\documentclass{gtart}

\def\ifplaintex{\expandafter\ifx\csname documentclass\endcsname\relax}

\def\gtp{{\mathsurround=0pt\it $\cal G\mskip-2mu$eometry \&\ 
$\cal T\!\!$opology $\cal P\!$ublications}}  

\def\recd{{\small Received:\qua\receiveddate\ifx\reviseddate\relax
\else\qquad Revised:\qua\reviseddate\fi\par}} 

\def\lognumber#1{\def\thelognumber{#1}}
\def\volumenumber#1{\def\thevolumenumber{#1}}
\def\volumeyear#1{\def\thevolumeyear{#1}}
\def\papernumber#1{\def\thepapernumber{#1}}
\def\pagenumbers#1#2{\def\startpage{#1}\def\finishpage{#2}}
\def\published#1{\def\publishdate{#1}}

\def\received#1{\def\receiveddate{#1}}
\def\revised#1{\def\reviseddate{#1}}
\def\accepted#1{\def\accepteddate{#1}}

\long\def\asciiabstract#1{\long\def\theasciiabstract{#1}}

\let\\\par\let\thelognumber\relax\let\thevolumenumber\relax
\let\thepapernumber\relax\let\thevolumeyear\relax\let\startpage\relax
\let\finishpage\relax\let\publishdate\relax\let\receiveddate\relax
\let\reviseddate\relax\let\accepteddate\relax\let\theasciititle\relax
\let\theasciiauthors\relax
\let\theasciiabstract\relax

\let\theasciiemail\relax

\ifplaintex
\font\logobig=cmssbx10 scaled 3836
\font\logomed=cmssbx10 scaled 2557
\else
\font\logobig=cmssbx10 scaled 4200
\font\logomed=cmssbx10 scaled 2800
\fi

\long\def\makeagttitle{   
\count0=\startpage
\agt\hfill      
\hbox to 45truept{\vbox to 0pt{\vglue -13truept{\logomed A\kern -.37em{\logobig 
T}\kern -.38em G}\vss}\hss}
\break
{\small Volume \thevolumenumber\ (\thevolumeyear)
\startpage--\finishpage\nl
Published: \publishdate}

\vglue .25truein

{\parskip=0pt\leftskip 0pt plus
1fil\def\\{\par\smallskip}{\Large\bf\thetitle}\par\medskip} \vglue
0.05truein

{\parskip=0pt\leftskip 0pt plus 1fil\def\\{\par}{\sc\theauthors}
\par\medskip}%
 
\vglue 0.03truein 

{\small\leftskip 25truept\rightskip 25truept{\bf Abstract}\stdspace\theabstract

{\bf AMS Classification}\stdspace\theprimaryclass
\ifx\thesecondaryclass\relax\else; \thesecondaryclass\fi\par
{\bf Keywords}\stdspace \thekeywords\par}\vglue 7truept

}   

\ifplaintex
\font\phead=cmsl9 scaled 950
\font\pnum=cmbx10 scaled 913
\font\pfoot=cmsl9 scaled 950
\headline{\vbox to 0pt{\vskip -4.5mm\line{\small\phead\ifnum
\count0=\startpage ISSN 1472-2739 (on-line) 1472-2747 (printed)
\hfill {\pnum\folio}\else\ifodd\count0\def\\{ }%
\ifx\theshorttitle\relax\thetitle\else\theshorttitle\fi\hfill{\pnum\folio}
\else\def\\{ and }{\pnum\folio}\hfill\ifx\theshortauthors\relax\theauthors
\else\theshortauthors\fi\fi\fi}\vss}}
\footline{\vbox to 0pt{\vglue 0mm\line{\small\pfoot\ifnum\count0=\startpage
\copyright\ \gtp\hfill\else
\agt, Volume \thevolumenumber\ (\thevolumeyear)\hfill\fi}\vss}}
\else
\font=cmsl9 scaled 1050
\font\lnum=cmbx10 
\font=cmsl9 scaled 1050
\makeatletter
\def\@oddhead{{\small\lhead\ifnum\count0=\startpage ISSN 1472-2739 
(on-line) 1472-2747 (printed)\hfill {\lnum\number\count0}\else\ifodd\count0
\def\\{ }\ifx\theshorttitle\relax \thetitle \else\theshorttitle\fi\hfill
{\lnum\number\count0}\else\def\\{ and }{\lnum\number\count0}
\hfill\ifx\theshortauthors\relax 
\theauthors\else\theshortauthors\fi\fi\fi}}\def\@evenhead{\@oddhead}
\def\@oddfoot{\small\lfoot\ifnum\count0=\startpage\copyright\ \gtp\hfill\else
\agt, Volume \thevolumenumber\ (\thevolumeyear)\hfill\fi}
\def\@evenfoot{\@oddfoot}
\makeatother
\fi
\let\maketitlepage\makeagttitle
\let\makeshorttitle\maketitlepage
\let\maketitle\maketitlepage

\newwrite\gtoutfile
\long\gdef\makeheadfile{  
{\def\\{, }\def\s{ }
\immediate\openout\gtoutfile head.xxx
\immediate\write\gtoutfile{To: math@arxiv.org}
\immediate\write\gtoutfile{Subject: put OR rep NNNNN:ppppp}
\immediate\write\gtoutfile{--text follows this line--}
\immediate\write\gtoutfile{Proxy-for: \ifx\theasciiauthors\relax
\theauthors\else\theasciiauthors\fi\s<\ifx\theasciiemail\relax\theemail\else\theasciiemail\fi>}
\immediate\write\gtoutfile{\noexpand\\}
\immediate\write\gtoutfile{Authors: \ifx\theasciiauthors\relax
\theauthors\else\theasciiauthors\fi}
{\def\\{ }\immediate\write\gtoutfile{Title: \ifx\theasciititle\relax
\thetitle\else\theasciititle\fi}}
\immediate\write\gtoutfile{Subj-class: GT or SG, GR etc}
\immediate\write\gtoutfile{MSC-class: \theprimaryclass\ifx\thesecondaryclass\relax\else, \thesecondaryclass\fi}
\immediate\write\gtoutfile{Journal-ref: Algebr. Geom. Topol. \thevolumenumber\s
(\thevolumeyear) \startpage-\finishpage}
\immediate\write\gtoutfile{Comments: Published by Algebraic and
Geometric Topology at}
\immediate\write\gtoutfile{\s\s\s  http://www.maths.warwick.ac.uk/agt/AGTVol\thevolumenumber/agt-\thevolumenumber-\thepapernumber.abs.html}
\immediate\write\gtoutfile{\noexpand\\}
\immediate\write\gtoutfile{}
\ifx\theasciiabstract\relax
\immediate\write\gtoutfile{\theabstract}\else
\immediate\write\gtoutfile{\theasciiabstract}\fi
\immediate\write\gtoutfile{}
\immediate\write\gtoutfile{\noexpand\\}
\immediate\write\gtoutfile{}
\immediate\closeout\gtoutfile}}  

\def\maketitlepage{\makeagttitle\makeheadfile}
\let\makeshorttitle\maketitlepage
\let\maketitle\maketitlepage

\def\ifplaintex{\expandafter\ifx\csname documentclass\endcsname\relax}

\def\gtp{{\mathsurround=0pt\it $\cal G\mskip-2mu$eometry \&\ 
$\cal T\!\!$opology $\cal P\!$ublications}}  

\def\recd{{\small Received:\qua\receiveddate\ifx\reviseddate\relax
\else\qquad Revised:\qua\reviseddate\fi\par}} 

\def\lognumber#1{\def\thelognumber{#1}}
\def\volumenumber#1{\def\thevolumenumber{#1}}
\def\volumeyear#1{\def\thevolumeyear{#1}}
\def\papernumber#1{\def\thepapernumber{#1}}
\def\pagenumbers#1#2{\def\startpage{#1}\def\finishpage{#2}}
\def\published#1{\def\publishdate{#1}}

\def\received#1{\def\receiveddate{#1}}
\def\revised#1{\def\reviseddate{#1}}
\def\accepted#1{\def\accepteddate{#1}}

\long\def\asciiabstract#1{\long\def\theasciiabstract{#1}}

\let\\\par\let\thelognumber\relax\let\thevolumenumber\relax
\let\thepapernumber\relax\let\thevolumeyear\relax\let\startpage\relax
\let\finishpage\relax\let\publishdate\relax\let\receiveddate\relax
\let\reviseddate\relax\let\accepteddate\relax\let\theasciititle\relax
\let\theasciiauthors\relax
\let\theasciiabstract\relax

\let\theasciiemail\relax

\ifplaintex
\font\logobig=cmssbx10 scaled 3836
\font\logomed=cmssbx10 scaled 2557
\else
\font\logobig=cmssbx10 scaled 4200
\font\logomed=cmssbx10 scaled 2800
\fi

\long\def\makeagttitle{   
\count0=\startpage
\agt\hfill      
\hbox to 45truept{\vbox to 0pt{\vglue -13truept{\logomed A\kern -.37em{\logobig 
T}\kern -.38em G}\vss}\hss}
\break
{\small Volume \thevolumenumber\ (\thevolumeyear)
\startpage--\finishpage\nl
Published: \publishdate}

\vglue .25truein

{\parskip=0pt\leftskip 0pt plus
1fil\def\\{\par\smallskip}{\Large\bf\thetitle}\par\medskip} \vglue
0.05truein

{\parskip=0pt\leftskip 0pt plus 1fil\def\\{\par}{\sc\theauthors}
\par\medskip}%
 
\vglue 0.03truein 

{\small\leftskip 25truept\rightskip 25truept{\bf Abstract}\stdspace\theabstract

{\bf AMS Classification}\stdspace\theprimaryclass
\ifx\thesecondaryclass\relax\else; \thesecondaryclass\fi\par
{\bf Keywords}\stdspace \thekeywords\par}\vglue 7truept

}   

\ifplaintex
\font\phead=cmsl9 scaled 950
\font\pnum=cmbx10 scaled 913
\font\pfoot=cmsl9 scaled 950
\headline{\vbox to 0pt{\vskip -4.5mm\line{\small\phead\ifnum
\count0=\startpage ISSN 1472-2739 (on-line) 1472-2747 (printed)
\hfill {\pnum\folio}\else\ifodd\count0\def\\{ }%
\ifx\theshorttitle\relax\thetitle\else\theshorttitle\fi\hfill{\pnum\folio}
\else\def\\{ and }{\pnum\folio}\hfill\ifx\theshortauthors\relax\theauthors
\else\theshortauthors\fi\fi\fi}\vss}}
\footline{\vbox to 0pt{\vglue 0mm\line{\small\pfoot\ifnum\count0=\startpage
\copyright\ \gtp\hfill\else
\agt, Volume \thevolumenumber\ (\thevolumeyear)\hfill\fi}\vss}}
\else
\font=cmsl9 scaled 1050
\font\lnum=cmbx10 
\font=cmsl9 scaled 1050
\makeatletter
\def\@oddhead{{\small\lhead\ifnum\count0=\startpage ISSN 1472-2739 
(on-line) 1472-2747 (printed)\hfill {\lnum\number\count0}\else\ifodd\count0
\def\\{ }\ifx\theshorttitle\relax \thetitle \else\theshorttitle\fi\hfill
{\lnum\number\count0}\else\def\\{ and }{\lnum\number\count0}
\hfill\ifx\theshortauthors\relax 
\theauthors\else\theshortauthors\fi\fi\fi}}\def\@evenhead{\@oddhead}
\def\@oddfoot{\small\lfoot\ifnum\count0=\startpage\copyright\ \gtp\hfill\else
\agt, Volume \thevolumenumber\ (\thevolumeyear)\hfill\fi}
\def\@evenfoot{\@oddfoot}
\makeatother
\fi
\let\maketitlepage\makeagttitle
\let\makeshorttitle\maketitlepage
\let\maketitle\maketitlepage

\newwrite\gtoutfile
\long\gdef\makeheadfile{  
{\def\\{, }\def\s{ }
\immediate\openout\gtoutfile head.xxx
\immediate\write\gtoutfile{To: math@arxiv.org}
\immediate\write\gtoutfile{Subject: put OR rep NNNNN:ppppp}
\immediate\write\gtoutfile{--text follows this line--}
\immediate\write\gtoutfile{Proxy-for: \ifx\theasciiauthors\relax
\theauthors\else\theasciiauthors\fi\s<\ifx\theasciiemail\relax\theemail\else\theasciiemail\fi>}
\immediate\write\gtoutfile{\noexpand\\}
\immediate\write\gtoutfile{Authors: \ifx\theasciiauthors\relax
\theauthors\else\theasciiauthors\fi}
{\def\\{ }\immediate\write\gtoutfile{Title: \ifx\theasciititle\relax
\thetitle\else\theasciititle\fi}}
\immediate\write\gtoutfile{Subj-class: GT or SG, GR etc}
\immediate\write\gtoutfile{MSC-class: \theprimaryclass\ifx\thesecondaryclass\relax\else, \thesecondaryclass\fi}
\immediate\write\gtoutfile{Journal-ref: Algebr. Geom. Topol. \thevolumenumber\s
(\thevolumeyear) \startpage-\finishpage}
\immediate\write\gtoutfile{Comments: Published by Algebraic and
Geometric Topology at}
\immediate\write\gtoutfile{\s\s\s  http://www.maths.warwick.ac.uk/agt/AGTVol\thevolumenumber/agt-\thevolumenumber-\thepapernumber.abs.html}
\immediate\write\gtoutfile{\noexpand\\}
\immediate\write\gtoutfile{}
\ifx\theasciiabstract\relax
\immediate\write\gtoutfile{\theabstract}\else
\immediate\write\gtoutfile{\theasciiabstract}\fi
\immediate\write\gtoutfile{}
\immediate\write\gtoutfile{\noexpand\\}
\immediate\write\gtoutfile{}
\immediate\closeout\gtoutfile}}  

\def\maketitlepage{\makeagttitle\makeheadfile}
\let\makeshorttitle\maketitlepage
\let\maketitle\maketitlepage

\def\ifplaintex{\expandafter\ifx\csname documentclass\endcsname\relax}

\def\gtp{{\mathsurround=0pt\it $\cal G\mskip-2mu$eometry \&\ 
$\cal T\!\!$opology $\cal P\!$ublications}}  

\def\recd{{\small Received:\qua\receiveddate\ifx\reviseddate\relax
\else\qquad Revised:\qua\reviseddate\fi\par}} 

\def\lognumber#1{\def\thelognumber{#1}}
\def\volumenumber#1{\def\thevolumenumber{#1}}
\def\volumeyear#1{\def\thevolumeyear{#1}}
\def\papernumber#1{\def\thepapernumber{#1}}
\def\pagenumbers#1#2{\def\startpage{#1}\def\finishpage{#2}}
\def\published#1{\def\publishdate{#1}}

\def\received#1{\def\receiveddate{#1}}
\def\revised#1{\def\reviseddate{#1}}
\def\accepted#1{\def\accepteddate{#1}}

\long\def\asciiabstract#1{\long\def\theasciiabstract{#1}}

\let\\\par\let\thelognumber\relax\let\thevolumenumber\relax
\let\thepapernumber\relax\let\thevolumeyear\relax\let\startpage\relax
\let\finishpage\relax\let\publishdate\relax\let\receiveddate\relax
\let\reviseddate\relax\let\accepteddate\relax\let\theasciititle\relax
\let\theasciiauthors\relax
\let\theasciiabstract\relax

\let\theasciiemail\relax

\ifplaintex
\font\logobig=cmssbx10 scaled 3836
\font\logomed=cmssbx10 scaled 2557
\else
\font\logobig=cmssbx10 scaled 4200
\font\logomed=cmssbx10 scaled 2800
\fi

\long\def\makeagttitle{   
\count0=\startpage
\agt\hfill      
\hbox to 45truept{\vbox to 0pt{\vglue -13truept{\logomed A\kern -.37em{\logobig 
T}\kern -.38em G}\vss}\hss}
\break
{\small Volume \thevolumenumber\ (\thevolumeyear)
\startpage--\finishpage\nl
Published: \publishdate}

\vglue .25truein

{\parskip=0pt\leftskip 0pt plus
1fil\def\\{\par\smallskip}{\Large\bf\thetitle}\par\medskip} \vglue
0.05truein

{\parskip=0pt\leftskip 0pt plus 1fil\def\\{\par}{\sc\theauthors}
\par\medskip}%
 
\vglue 0.03truein 

{\small\leftskip 25truept\rightskip 25truept{\bf Abstract}\stdspace\theabstract

{\bf AMS Classification}\stdspace\theprimaryclass
\ifx\thesecondaryclass\relax\else; \thesecondaryclass\fi\par
{\bf Keywords}\stdspace \thekeywords\par}\vglue 7truept

}   

\ifplaintex
\font\phead=cmsl9 scaled 950
\font\pnum=cmbx10 scaled 913
\font\pfoot=cmsl9 scaled 950
\headline{\vbox to 0pt{\vskip -4.5mm\line{\small\phead\ifnum
\count0=\startpage ISSN 1472-2739 (on-line) 1472-2747 (printed)
\hfill {\pnum\folio}\else\ifodd\count0\def\\{ }%
\ifx\theshorttitle\relax\thetitle\else\theshorttitle\fi\hfill{\pnum\folio}
\else\def\\{ and }{\pnum\folio}\hfill\ifx\theshortauthors\relax\theauthors
\else\theshortauthors\fi\fi\fi}\vss}}
\footline{\vbox to 0pt{\vglue 0mm\line{\small\pfoot\ifnum\count0=\startpage
\copyright\ \gtp\hfill\else
\agt, Volume \thevolumenumber\ (\thevolumeyear)\hfill\fi}\vss}}
\else
\font=cmsl9 scaled 1050
\font\lnum=cmbx10 
\font=cmsl9 scaled 1050
\makeatletter
\def\@oddhead{{\small\lhead\ifnum\count0=\startpage ISSN 1472-2739 
(on-line) 1472-2747 (printed)\hfill {\lnum\number\count0}\else\ifodd\count0
\def\\{ }\ifx\theshorttitle\relax \thetitle \else\theshorttitle\fi\hfill
{\lnum\number\count0}\else\def\\{ and }{\lnum\number\count0}
\hfill\ifx\theshortauthors\relax 
\theauthors\else\theshortauthors\fi\fi\fi}}\def\@evenhead{\@oddhead}
\def\@oddfoot{\small\lfoot\ifnum\count0=\startpage\copyright\ \gtp\hfill\else
\agt, Volume \thevolumenumber\ (\thevolumeyear)\hfill\fi}
\def\@evenfoot{\@oddfoot}
\makeatother
\fi
\let\maketitlepage\makeagttitle
\let\makeshorttitle\maketitlepage
\let\maketitle\maketitlepage

\newwrite\gtoutfile
\long\gdef\makeheadfile{  
{\def\\{, }\def\s{ }
\immediate\openout\gtoutfile head.xxx
\immediate\write\gtoutfile{To: math@arxiv.org}
\immediate\write\gtoutfile{Subject: put OR rep NNNNN:ppppp}
\immediate\write\gtoutfile{--text follows this line--}
\immediate\write\gtoutfile{Proxy-for: \ifx\theasciiauthors\relax
\theauthors\else\theasciiauthors\fi\s<\ifx\theasciiemail\relax\theemail\else\theasciiemail\fi>}
\immediate\write\gtoutfile{\noexpand\\}
\immediate\write\gtoutfile{Authors: \ifx\theasciiauthors\relax
\theauthors\else\theasciiauthors\fi}
{\def\\{ }\immediate\write\gtoutfile{Title: \ifx\theasciititle\relax
\thetitle\else\theasciititle\fi}}
\immediate\write\gtoutfile{Subj-class: GT or SG, GR etc}
\immediate\write\gtoutfile{MSC-class: \theprimaryclass\ifx\thesecondaryclass\relax\else, \thesecondaryclass\fi}
\immediate\write\gtoutfile{Journal-ref: Algebr. Geom. Topol. \thevolumenumber\s
(\thevolumeyear) \startpage-\finishpage}
\immediate\write\gtoutfile{Comments: Published by Algebraic and
Geometric Topology at}
\immediate\write\gtoutfile{\s\s\s  http://www.maths.warwick.ac.uk/agt/AGTVol\thevolumenumber/agt-\thevolumenumber-\thepapernumber.abs.html}
\immediate\write\gtoutfile{\noexpand\\}
\immediate\write\gtoutfile{}
\ifx\theasciiabstract\relax
\immediate\write\gtoutfile{\theabstract}\else
\immediate\write\gtoutfile{\theasciiabstract}\fi
\immediate\write\gtoutfile{}
\immediate\write\gtoutfile{\noexpand\\}
\immediate\write\gtoutfile{}
\immediate\closeout\gtoutfile}}  

\def\maketitlepage{\makeagttitle\makeheadfile}
\let\makeshorttitle\maketitlepage
\let\maketitle\maketitlepage

\def\ifplaintex{\expandafter\ifx\csname documentclass\endcsname\relax}

\def\gtp{{\mathsurround=0pt\it $\cal G\mskip-2mu$eometry \&\ 
$\cal T\!\!$opology $\cal P\!$ublications}}  

\def\recd{{\small Received:\qua\receiveddate\ifx\reviseddate\relax
\else\qquad Revised:\qua\reviseddate\fi\par}} 

\def\lognumber#1{\def\thelognumber{#1}}
\def\volumenumber#1{\def\thevolumenumber{#1}}
\def\volumeyear#1{\def\thevolumeyear{#1}}
\def\papernumber#1{\def\thepapernumber{#1}}
\def\pagenumbers#1#2{\def\startpage{#1}\def\finishpage{#2}}
\def\published#1{\def\publishdate{#1}}

\def\received#1{\def\receiveddate{#1}}
\def\revised#1{\def\reviseddate{#1}}
\def\accepted#1{\def\accepteddate{#1}}

\long\def\asciiabstract#1{\long\def\theasciiabstract{#1}}

\let\\\par\let\thelognumber\relax\let\thevolumenumber\relax
\let\thepapernumber\relax\let\thevolumeyear\relax\let\startpage\relax
\let\finishpage\relax\let\publishdate\relax\let\receiveddate\relax
\let\reviseddate\relax\let\accepteddate\relax\let\theasciititle\relax
\let\theasciiauthors\relax
\let\theasciiabstract\relax

\let\theasciiemail\relax

\ifplaintex
\font\logobig=cmssbx10 scaled 3836
\font\logomed=cmssbx10 scaled 2557
\else
\font\logobig=cmssbx10 scaled 4200
\font\logomed=cmssbx10 scaled 2800
\fi

\long\def\makeagttitle{   
\count0=\startpage
\agt\hfill      
\hbox to 45truept{\vbox to 0pt{\vglue -13truept{\logomed A\kern -.37em{\logobig 
T}\kern -.38em G}\vss}\hss}
\break
{\small Volume \thevolumenumber\ (\thevolumeyear)
\startpage--\finishpage\nl
Published: \publishdate}

\vglue .25truein

{\parskip=0pt\leftskip 0pt plus
1fil\def\\{\par\smallskip}{\Large\bf\thetitle}\par\medskip} \vglue
0.05truein

{\parskip=0pt\leftskip 0pt plus 1fil\def\\{\par}{\sc\theauthors}
\par\medskip}%
 
\vglue 0.03truein 

{\small\leftskip 25truept\rightskip 25truept{\bf Abstract}\stdspace\theabstract

{\bf AMS Classification}\stdspace\theprimaryclass
\ifx\thesecondaryclass\relax\else; \thesecondaryclass\fi\par
{\bf Keywords}\stdspace \thekeywords\par}\vglue 7truept

}   

\ifplaintex
\font\phead=cmsl9 scaled 950
\font\pnum=cmbx10 scaled 913
\font\pfoot=cmsl9 scaled 950
\headline{\vbox to 0pt{\vskip -4.5mm\line{\small\phead\ifnum
\count0=\startpage ISSN 1472-2739 (on-line) 1472-2747 (printed)
\hfill {\pnum\folio}\else\ifodd\count0\def\\{ }%
\ifx\theshorttitle\relax\thetitle\else\theshorttitle\fi\hfill{\pnum\folio}
\else\def\\{ and }{\pnum\folio}\hfill\ifx\theshortauthors\relax\theauthors
\else\theshortauthors\fi\fi\fi}\vss}}
\footline{\vbox to 0pt{\vglue 0mm\line{\small\pfoot\ifnum\count0=\startpage
\copyright\ \gtp\hfill\else
\agt, Volume \thevolumenumber\ (\thevolumeyear)\hfill\fi}\vss}}
\else
\font=cmsl9 scaled 1050
\font\lnum=cmbx10 
\font=cmsl9 scaled 1050
\makeatletter
\def\@oddhead{{\small\lhead\ifnum\count0=\startpage ISSN 1472-2739 
(on-line) 1472-2747 (printed)\hfill {\lnum\number\count0}\else\ifodd\count0
\def\\{ }\ifx\theshorttitle\relax \thetitle \else\theshorttitle\fi\hfill
{\lnum\number\count0}\else\def\\{ and }{\lnum\number\count0}
\hfill\ifx\theshortauthors\relax 
\theauthors\else\theshortauthors\fi\fi\fi}}\def\@evenhead{\@oddhead}
\def\@oddfoot{\small\lfoot\ifnum\count0=\startpage\copyright\ \gtp\hfill\else
\agt, Volume \thevolumenumber\ (\thevolumeyear)\hfill\fi}
\def\@evenfoot{\@oddfoot}
\makeatother
\fi
\let\maketitlepage\makeagttitle
\let\makeshorttitle\maketitlepage
\let\maketitle\maketitlepage

\newwrite\gtoutfile
\long\gdef\makeheadfile{  
{\def\\{, }\def\s{ }
\immediate\openout\gtoutfile head.xxx
\immediate\write\gtoutfile{To: math@arxiv.org}
\immediate\write\gtoutfile{Subject: put OR rep NNNNN:ppppp}
\immediate\write\gtoutfile{--text follows this line--}
\immediate\write\gtoutfile{Proxy-for: \ifx\theasciiauthors\relax
\theauthors\else\theasciiauthors\fi\s<\ifx\theasciiemail\relax\theemail\else\theasciiemail\fi>}
\immediate\write\gtoutfile{\noexpand\\}
\immediate\write\gtoutfile{Authors: \ifx\theasciiauthors\relax
\theauthors\else\theasciiauthors\fi}
{\def\\{ }\immediate\write\gtoutfile{Title: \ifx\theasciititle\relax
\thetitle\else\theasciititle\fi}}
\immediate\write\gtoutfile{Subj-class: GT or SG, GR etc}
\immediate\write\gtoutfile{MSC-class: \theprimaryclass\ifx\thesecondaryclass\relax\else, \thesecondaryclass\fi}
\immediate\write\gtoutfile{Journal-ref: Algebr. Geom. Topol. \thevolumenumber\s
(\thevolumeyear) \startpage-\finishpage}
\immediate\write\gtoutfile{Comments: Published by Algebraic and
Geometric Topology at}
\immediate\write\gtoutfile{\s\s\s  http://www.maths.warwick.ac.uk/agt/AGTVol\thevolumenumber/agt-\thevolumenumber-\thepapernumber.abs.html}
\immediate\write\gtoutfile{\noexpand\\}
\immediate\write\gtoutfile{}
\ifx\theasciiabstract\relax
\immediate\write\gtoutfile{\theabstract}\else
\immediate\write\gtoutfile{\theasciiabstract}\fi
\immediate\write\gtoutfile{}
\immediate\write\gtoutfile{\noexpand\\}
\immediate\write\gtoutfile{}
\immediate\closeout\gtoutfile}}  

\def\maketitlepage{\makeagttitle\makeheadfile}
\let\makeshorttitle\maketitlepage
\let\maketitle\maketitlepage

\lognumber{42}
\volumenumber{2}
\volumeyear{2002}
\papernumber{42}
\pagenumbers{1061}{1118}
\received{3 April 2001}
\revised{5 November 2002}
\accepted{19 November 2002}
\published{25 November 2002}

\usepackage{amsmath,amssymb} 
\usepackage{pb-diagram,lamsarrow,pb-lams}

\newcommand{\Div}       {\operatorname{Div}}
\newcommand{\End}       {\operatorname{End}}
\newcommand{\Hom}       {\operatorname{Hom}}
\newcommand{\Ind}       {\operatorname{Ind}}
\newcommand{\Int}       {\operatorname{Int}}
\newcommand{\Prim}      {\operatorname{Prim}}
\newcommand{\Sub}       {\operatorname{Sub}}
\newcommand{\Sym}       {\operatorname{Sym}}
\newcommand{\Tor}       {\operatorname{Tor}}

\newcommand{\alt}       {\operatorname{alt}}
\newcommand{\codim}     {\operatorname{codim}}
\newcommand{\cok}       {\operatorname{cok}}
\newcommand{\img}       {\operatorname{image}}
\renewcommand{\int}     {\operatorname{int}}
\newcommand{\rank}      {\operatorname{rank}}
\newcommand{\sgn}       {\operatorname{sgn}}
\newcommand{\spec}      {\operatorname{spec}}
\newcommand{\spf}       {\operatorname{spf}}
\newcommand{\tr}        {\operatorname{tr}}

\newcommand{\Z}         {{\mathbb{Z}}}

\newcommand{\R}         {{\mathbb{R}}}
\newcommand{\C}         {{\mathbb{C}}}
\renewcommand{\H}       {{\mathbb{H}}}

\newcommand{\GG}        {{\mathbb{G}}}
\renewcommand{\O}       {\mathcal{O}}
\newcommand{\OG}        {{\mathcal{O}_\GG}}
\newcommand{\OS}        {{\mathcal{O}_S}}
\newcommand{\OT}        {{\mathcal{O}_T}}
\newcommand{\OD}        {{\mathcal{O}_D}}

\newcommand{\CM}        {{\mathcal{M}}}

\newcommand{\al}        {\alpha}
\newcommand{\bt}        {\beta} 
\newcommand{\gm}        {\gamma}
\newcommand{\dl}        {\delta}
\newcommand{\ep}        {\epsilon}
\newcommand{\lm}        {\lambda}
\newcommand{\sg}        {\sigma}
\newcommand{\zt}        {\zeta}
\newcommand{\tht}       {\theta}
\newcommand{\og}        {\omega}

\newcommand{\Gm}        {\Gamma}
\newcommand{\Sg}        {\Sigma}
\newcommand{\Om}        {\Omega}
\newcommand{\Oml}       {\Omega^{\text{lau}}}
\newcommand{\Omp}       {\Omega^{\text{pol}}}

\newcommand{\CP}        {{\mathbb{C}P}}
\newcommand{\CPi}       {{\mathbb{C}P^\infty}}
\newcommand{\Ci}        {{\mathbb{C}}^\infty}
\newcommand{\HP}        {{\mathbb{H}P}}
\newcommand{\op}        {\oplus}
\newcommand{\om}        {\ominus}
\newcommand{\psb}[1]    {[\![#1]\!]}
\newcommand{\lsb}[1]    {(\!(#1)\!)}
\newcommand{\st}        {\;|\;}
\newcommand{\tE}        {\widetilde{E}}
\newcommand{\tm}        {\times}
\newcommand{\xra}       {\xrightarrow}
\newcommand{\xla}       {\xleftarrow}
\newcommand{\iffa}      {\Leftrightarrow}
\newcommand{\ot}        {\otimes}
\newcommand{\hot}       {\widehat{\otimes}}
\newcommand{\sse}       {\subseteq}
\newcommand{\Smash}     {\wedge}
\newcommand{\ov}[1]     {\overline{#1}}
\newcommand{\un}[1]     {\underline{#1}}
\newcommand{\um}        {{\underline{n}}}
\newcommand{\pnt}       {\text{point}}

\newcommand{\convto}    {\Longrightarrow}
\newcommand{\tH}        {\widetilde{H}}

\newcommand{\colim}  {\operatornamewithlimits{\underset{\longrightarrow}{lim}}}
\newcommand{\invlim} {\operatornamewithlimits{\underset{\longleftarrow}{lim}}}
\newcommand{\ip}[1]     {\langle #1\rangle}
\newcommand{\uu}        {{\mathfrak{u}}}

\newcommand{\Wedge}     {\vee}
\newcommand{\bcf}[2]{\left(\begin{array}{c}{#1}\\{#2}\end{array}\right)}
\newcommand{\sm}        {\setminus}
\newcommand{\phb}       {\overline{\phi}}
\newcommand{\jh}        {\hat{\jmath}}

\makeatletter
\@namedef{dgo@=>}{\let\dg@VECTOR=\dg@forkvector}%

\def\dg@forkvector(#1,#2)#3{%
   \begingroup
   \dg@XTEMP=#1\relax \dg@YTEMP=#2\relax
   \dg@XEND=1
   \multiply\dg@XEND\dg@XTEMP 
   \multiply\dg@XEND  50 \relax 
   \dg@WEND=1
   \multiply\dg@WEND\dg@YTEMP 
   \multiply\dg@WEND -50 \relax
   \begin{picture}(0,0)
      \put(\dg@WEND,\dg@XEND){\lamsvector(\dg@XTEMP,\dg@YTEMP){#3}}
      \multiply\dg@XEND -1 \relax
      \multiply\dg@WEND -1 \relax
      \put(\dg@WEND,\dg@XEND){\lamsvector(\dg@XTEMP,\dg@YTEMP){#3}}
   \end{picture}%
   \endgroup}%
\makeatother

\renewcommand{\:}{\colon\thinspace}

\newtheorem{theorem}{Theorem}[section]

\newtheorem{lemma}[theorem]{Lemma}
\newtheorem{proposition}[theorem]{Proposition}
\newtheorem{corollary}[theorem]{Corollary}

\theoremstyle{definition}
\newtheorem{construction}[theorem]{Construction}
\newtheorem{remark}[theorem]{Remark}
\newtheorem{definition}[theorem]{Definition}

\newenvironment{diag}{
 \renewcommand{\typeout}[1]{}
 \begin{displaymath}
 \begin{diagram}}{
 \end{diagram}
 \end{displaymath}} 

\begin{document}
\title{Common subbundles and intersections of divisors}
\author{N. P. Strickland}
\address{Department of Mathematics,
University of Sheffield\\Western Bank, Sheffield, S10 2TN,UK}
\email{N.P.Strickland@sheffield.ac.uk}
\begin{abstract}
 Let $V_0$ and $V_1$ be complex vector bundles over a space $X$.  We
 use the theory of divisors on formal groups to give obstructions in
 generalised cohomology that vanish when $V_0$ and $V_1$ can be
 embedded in a bundle $U$ in such a way that $V_0\cap V_1$ has
 dimension at least $k$ everywhere.  We study various algebraic
 universal examples related to this question, and show that they arise
 from the generalised cohomology of corresponding topological
 universal examples.  This extends and reinterprets earlier work on
 degeneracy classes in ordinary cohomology or intersection theory.
\end{abstract}

\asciiabstract{Let V_0 and V_1 be complex vector bundles over a space
 X.  We use the theory of divisors on formal groups to give
 obstructions in generalised cohomology that vanish when V_0 and V_1
 can be embedded in a bundle U in such a way that V_0\cap V_1 has
 dimension at least k everywhere.  We study various algebraic
 universal examples related to this question, and show that they arise
 from the generalised cohomology of corresponding topological
 universal examples.  This extends and reinterprets earlier work on
 degeneracy classes in ordinary cohomology or intersection theory.}



\maketitle

\section{Introduction}
\label{sec-intro}

There are a number of different motivations for the theory developed
here, but perhaps the most obvious is as follows.  Suppose we have a
space $X$ with vector bundles $V_0$ and $V_1$.  (Throughout this
paper, the term ``vector space'' refers to finite-dimensional complex
vector spaces equipped with Hermitian inner products, and similarly
for ``vector bundle''.)  We define the \emph{intersection index} of
$V_0$ and $V_1$ to be the largest $k$ such that $V_0$ and $V_1$ can be
embedded isometrically in some bundle $U$ in such a way that
$\dim(V_{0x}\cap V_{1x})\geq k$ for all $x\in X$.  We write
$\int(V_0,V_1)$ for this intersection index.  Our aim is to use
invariants from generalised cohomology theory to estimate
$\int(V_0,V_1)$, and to investigate the topology of certain universal
examples related to this question.

We will show in Proposition~\ref{prop-int-equiv} that $\int(V_0,V_1)$
is also the largest $k$ such that there is a linear map $V_0\xra{}V_1$
of rank at least $k$ everywhere.  This creates a link with the theory
of degeneracy loci and the corresponding classes in the cohomology of
manifolds or Chow rings of varieties, which are given by the
determinantal formula of Thom and Porteous.  The paper~\cite{pr:egd}
by Pragacz is a convenient reference for comparison with the present
work.  The relevant theory is based strongly on Schubert calculus, and
could presumably be transferred to complex cobordism (and thus to
other complex-orientable theories) by the methods of Bressler and
Evens~\cite{brev:scc}.

However, our approach will be different in a number of ways.  Firstly,
we use the language of formal groups, as discussed in~\cite{st:fsfg}
(for example).  We fix an even periodic cohomology theory $E$ with a
complex orientation $x\in\tE^0\CPi$.  For any space $X$ we have a
formal scheme $X_E=\spf(E^0X)$, the basic examples being $S:=(\pnt)_E$
and $\GG:=\CP^\infty_E=\spf(E^0\psb{x})$, which is a formal group over
$S$.  If $V$ is a complex vector bundle over $X$, we write $PV$ for
the associated bundle of projective spaces.  It is standard that
$E^0(PV)=E^0(X)\psb{x}/f_V(x)$, where
$f_V(x)=\sum_{i+j=\dim(V)}c_ix^j$, where $c_i$ is the $i$'th Chern
class of $V$.  In geometric terms, this means that the formal scheme
$D(V):=(PV)_E$ is naturally embedded as a divisor in $\GG\tm_SX_E$.
Most of our algebraic constructions will have a very natural
interpretation in terms of such divisors.  We will also consider the
bundle $U(V)=\coprod_{x\in X}U(V_x)$ of unitary groups associated to
$V$.  A key point is that $E^*U(V)$ is the exterior algebra over
$E^*X$ generated by $E^{*-1}PV$.  This provides a very natural link
with exterior algebra, and could be regarded as the ``real reason''
for the appearance of determinantal formulae, which seem rather
accidental in other approaches.  Our divisorial approach also leads to
descriptions of various cohomology rings that are manifestly
independent of the choice of complex orientation, and depend
functorially on $\GG$.  This functorality implicitly encodes the
action of stable cohomology operations and thus gives a tighter link
with the underlying homotopy theory.  

We were also influenced by work of Kitchloo~\cite{ki:css}, in which he
investigates the cohomological effect of Miller's stable splitting of
$U(n)$, and draws a link with the theory of Schur functions.

In Section~\ref{sec-int-div} we use the theory of Fitting ideals to
define an intersection index $\int(D_0,D_1)$, where $D_0$ and $D_1$
are divisors on $\GG$.  In Section~\ref{sec-unitary-bundles} we
identify $E^*U(V)$ with the exterior algebra generated by $E^{*-1}PV$,
and show that this identification is an isomorphism of Hopf algebras.
In Section~\ref{sec-int-bundles} we use this to prove our first main
theorem, that $\int(V_0,V_1)\leq\int(D(V_0),D(V_1))$; this implicitly
gives all the relations among Chern classes that are universally
satisfied when $\int(V_0,V_1)\geq k$ for some given integer $k$.
Next, in Section~\ref{sec-alg-universal} we study the universal
examples of our various algebraic questions, focusing on the scheme
$\Int_r(d_0,d_1)$ which classifies pairs $(D_0,D_1)$ of divisors of
degrees $d_0$ and $d_1$ such that $\int(D_0,D_1)\geq k$.  Our next
task is to construct spaces whose associated schemes are these
algebraic universal examples.  In Section~\ref{sec-flags} we warm up
by giving a divisorial account of the generalised cohomology of
Grassmannians and flag spaces, and then in
Section~\ref{sec-top-universal} we show that the space
\[ I'_r(d_0,d_1) :=
    \{(V_0,V_1)\in G_{d_0}(\Ci)\tm G_{d_1}(\Ci) 
               \st \dim(V_0\cap V_1) \geq k
    \}
\]
satisfies $I'_r(d_0,d_1)_E=\Int_r(d_0,d_1)$.  (The origin of the
present work is that the author needed to compute the cohomology of
certain spaces similar to $I'_r(d_0,d_1)$ as input to another project;
it would take us too far afield to discuss the background.)  This
completes the main work of the paper, but we have added three more
sections exploring the isomorphism $E^*U(V)\simeq\lm^*E^{*-1}PV$ in
more detail.  Section~\ref{sec-PkD} treats some purely algebraic
questions related to this situation, and in Sections~\ref{sec-adjoint}
and~\ref{sec-loops} we translate all the algebra into homotopy
theory.  In particular, this gives a divisorial interpretation of the
work of Mitchell, Richter and others on filtrations of $\Om U(n)$: the
scheme associated to the $k$'th stage in the filtration of $\Om_XU(V)$
is $D(V)^k/\Sg_k$, and the scheme associated to $\Om_XU(V)$ is the
free formal group over $X_E$ generated by $D(V)$.

Appendix~\ref{apx-functional} gives a brief treatment of the
functional calculus for normal operators, which is used in a number of
places in the text.

\begin{remark}
 There is a theory of degeneracy loci for morphisms with symmetries,
 where the formulae involve Pfaffians instead of determinants.  It
 would clearly be a natural project to reexamine this theory from the
 point of view of the present paper, but so far we have nothing to say
 about this. 
\end{remark}

\section{Notation and conventions}

\subsection{Spheres}\label{subsec-spheres}

We take $\R^n\cup\{\infty\}$ as our definition of $S^n$, with $\infty$
as the basepoint; we distinguish $S^1$ from the homeomorphic space
$U(1):=\{z\in\C\st|z|=1\}$.  Where necessary, we use the homeomorphism
$\gm\:U(1)\xra{}S^1$ given by 
\begin{align*}
 \gm(z)      &= (z+1)(z-1)^{-1}/i \\
 \gm^{-1}(t) &= (it+1)/(it-1). 
\end{align*}
One checks that $\gm(e^{i\tht})=\cot(-\tht/2)$, which is a strictly
increasing function of $\tht$ for $0<\tht<2\pi$.

\subsection{Fibrewise spaces}

We will use various elementary concepts from fibrewise topology; the
book of Crabb and James~\cite{crja:fht} is a convenient reference.
Very few topological technicalities arise, as our fibrewise spaces are
all fibre bundles, and the fibres are usually finite complexes.
 
In particular, given spaces $U$ and $V$ over a space $X$, we write
$U\tm_XV$ for the fibre product, and $U^n_X$ for the fibre power
$U\tm_X\ldots\tm_XU$.  If $U$ is pointed (in other words, it has a
specified section $s\:X\xra{}U$) and $E$ is any cohomology theory we
write $\tE_X^*U=E^*(U,sX)$.  We also write $\Sg_XU$ for the fibrewise
suspension of $U$, which is the quotient of $S^1\tm U$ in which
$\{\infty\}\tm U\cup S^1\tm sX$ is collapsed to a copy of $X$.  This
satisfies $\tE_X^*\Sg_XU=\tE^{*-1}_XU$.  We also write $\Om_XU$ for
the fibrewise loop space of $U$, which is the space of maps
$\og\:S^1\xra{}U$ such that the composite $S^1\xra{}U\xra{}X$ is
constant and $\og(\infty)\in sX$.  If $V$ is another pointed space
over $X$, we write $U\Smash_XV$ for the fibrewise smash product.  If
$W$ is an unpointed space over $X$ then we write $W_{+X}=W\amalg X$,
which is a pointed space over $X$ in an obvious way.

\subsection{Tensor products over schemes}

If $T$ is a scheme and $M$, $N$ are modules over the ring $\O_T$, we
will write $M\ot_TN$ for $M\ot_{\O_T}N$.  Similarly, we write
$\lm^k_TM$ for $\lm^k_{\O_T}M$, the $k$'th exterior power of $M$ over
$\O_T$. 

\subsection{Free modules}

Given a ring $R$ and a set $T$, we write $R\{T\}$ for the free
$R$-module generated by $T$.

\section{Intersections of divisors}
\label{sec-int-div}

Let $\GG$ be a commutative, one-dimensional formal group over a scheme
$S$.  Choose a coordinate $x$ so that $\OG=\OS\psb{x}$.  Let $D_0$ and
$D_1$ be divisors on $\GG$ defined over $S$, with degrees $d_0$ and
$d_1$ respectively.  This means that
$\O_{D_i}=\OG/f_i=\OS\psb{x}/f_i(x)$ for some monic polynomial
$f_i(x)$ of degree $d_i$ such that $f_i(x)=x^{d_i}$ modulo nilpotents.
It follows that $\O_{D_i}$ is a free module of rank $d_i$ over $\OS$,
with basis $\{x^j\st 0\leq j<d_i\}$.

As $D_0$ and $D_1$ are closed subschemes of $\GG$ we can form their
intersection, so that 
\[ \O_{D_0\cap D_1}=\OG/(f_0,f_1)=\OS\psb{x}/(f_0(x),f_1(x)). \] 
Typically this will not be a projective module over $\OS$, so some
thought is required to give a useful notion of its size.  We will use
a measure coming from the theory of Fitting ideals, which we now
recall briefly.

Let $R$ be a commutative Noetherian ring, and let $M$ be a finitely
generated $R$-module.  We can then find a presentation
$P_1\xra{\phi_1}P_0\xra{\phi_0}M$, where $P_0$ and $P_1$ are finitely
generated projective modules of ranks $p_0$ and $p_1$ say, and
$M=\cok(\phi_1)$.  The exterior powers $\lm^jP_i$ are again finitely
generated projective modules.  We define $I_j(\phi_1)$ to be the
smallest ideal in $R$ modulo which we have $\lm^j(\phi_1)=0$.  More
concretely, if $P_0$ and $P_1$ are free then $\phi_1$ can be
represented by a matrix $A$ and $I_j(\phi_1)$ is generated by the
determinants of all $j\tm j$ submatrices of $A$.  We then define
$I_j(M)=I_{p_0-j}(\phi_1)$; this is called the $j$'th Fitting ideal of
$M$.  It is a fundamental fact that this is well-defined; this was
already known to Fitting (see~\cite[Chapter 3]{no:ffr}, for example),
but we give a proof for the convenience of the reader.

\begin{proposition}\label{prop-fitting}
 The ideal $I_j(M)$ is independent of the choice of presentation of
 $M$. 
\end{proposition}
\begin{proof}
 We temporarily write $I_j(M,P_*,\phi_*)$ for the ideal called
 $I_j(M)$ above.
 
 Put $N=\ker(\phi_0)$ and let $\bt\:N\xra{}P_0$ be the inclusion.
 Then $\phi_1$ factors as $P_1\xra{\al}N\xra{\bt}P_0$, where $\al$ is
 surjective.  For any ideal $J\leq R$ we see that $\lm^k\al$ is
 surjective mod $J$, so $\lm^k\phi_1$ is zero mod $J$ iff $\lm^k\bt$
 is zero mod $J$.  This condition depends only on the map
 $\phi_0\:P_0\xra{}M$, so we can legitimately define
 $I_j(M,P_0,\phi_0):=I_j(M,P_*,\phi_*)$.

 Now suppose we have another presentation
 $Q_1\xra{\psi_1}Q_0\xra{\psi_0}M$, where $Q_i$ has rank $q_i$.
 Define $\chi_0\:P_0\op Q_0\xra{}M$ by
 $(u,v)\mapsto\phi_0(u)+\psi_0(v)$.  It will suffice to prove that
 \[ I_j(M,P_0,\phi_0) = I_j(M,P_0\op Q_0,\chi_0) = I_j(M,Q_0,\psi_0),
 \] 
 and by symmetry we need only check the first of these.  By
 projectivity we can choose a map $\tht\:Q_0\xra{}P_0$ with
 $\phi_0\tht=\psi_0$, and define $\chi_1\:P_1\op Q_0\xra{}P_0\op Q_0$
 by $(u,v)\mapsto(\phi_1(u)-\tht(v),v)$.  It is easy to check that
 this gives another presentation
 \[ P_1\op Q_0 \xra{\chi_1} P_0\op Q_0 \xra{\chi_0} M. \]

 If $k\leq q_0$ then $\lm^k\chi_1$ is certainly nonzero, because the
 composite
 \[ \lm^k Q_0 \xra{} \lm^k(P_1\op Q_0) \xra{\lm^k\chi}
    \lm^k(P_0\op Q_0) \xra{} \lm^kQ_0
 \]
 is the identity, and $\lm^kQ_0\neq 0$.  If $k>q_0$ and
 $\lm^k\chi_1=0$ then (by restricting to
 $\lm^{k-q_0}P_1\ot\lm^{q_0}Q_0$) we see that $\lm^{k-q_0}\phi_1=0$.

 For the converse, notice that $\lm^*N$ is a graded ring for any
 module $N$, and that $\lm^*\al$ is a ring map for any homomorphism
 $\al$ of $R$-modules.  One can check that $\lm^{j+q_0}(P_1\op Q_0)$
 is contained in the ideal in $\lm^*(P_1\op Q_0)$ generated by
 $\lm^jP_1$.  It follows that if $\lm^j\phi_1=0$ then
 $\lm^{j+q_0}\chi_1=0$. 

 This shows that $I_r(\phi_1)=I_{r+q_0}(\chi_1)$, and thus that
 $I_r(M,P_0,\phi_0)=I_r(M,P_0\op Q_0,\chi_0)$, as required.
\end{proof}

It is clear that
\[ I_0(M) \leq \ldots\leq I_m(M) = R, \]
and we define
\[ \rank(M) = \rank_R(M) = \min\{r \st I_r(M)\neq 0\}. \]
We call $\rank(M)$ the \emph{Fitting rank} of $M$.  For example, if
$R$ is a principal ideal domain with fraction field $K$, one can check
that $\rank(M)=\dim_K(K\ot_R M)$ for all $M$.  However, we will mostly
be interested in rings $R$ with many nilpotents, for which there is no
such simple formula.

The following lemma is easily checked from the definitions.
\begin{lemma}\label{lem-rank-omni}
 \begin{itemize}
 \item[\rm(a)] The Fitting rank is the same as the ordinary rank for
  projective modules.
 \item[\rm(b)] If $N$ is a quotient of $M$ then $\rank(N)\leq\rank(M)$.
 \item[\rm(c)] If there is a presentation $P\xra{}Q\xra{}M$ then
  $\rank(Q)-\rank(P)\leq\rank(M)\leq\rank(Q)$. \qed
 \end{itemize}
\end{lemma}
(It is not true, however, that $\rank(M\op N)=\rank(M)+\rank(N)$;
indeed, if $a\neq 0$ and $a^2=0$ then $\rank(R/a)=0$ but 
$\rank(R/a\op R/a)=1$.)

\begin{definition}\label{defn-int-D}
 The \emph{intersection multiplicity} of $D_0$ and $D_1$ is the
 integer 
 \[ \int(D_0,D_1):=\rank_{\OS}(\O_{D_0\cap D_1}). \]
 We also put
 \[ \Int_r(D_0,D_1)=\spec(\OS/I_{r-1}(\O_{D_0\cap D_1})), \]
 which is the largest subscheme of $S$ over which we have
 $\int(D_0,D_1)\geq r$.
\end{definition}
\begin{remark}
 Let $S'$ be a scheme over $S$, so that $\GG':=\GG\tm_SS'$ is a formal
 group over $S'$.  We refer to divisors on $\GG'$ as divisors on $\GG$
 over $S'$.  Given two such divisors $D_0$ and $D_1$, we get a closed
 subscheme $\Int_r(D_0,D_1)\sse S'$.  We will use this kind of
 base-change construction throughout the paper without explicit
 comment. 
\end{remark}

To make the above definitions more explicit, we will describe several
different presentations of $\O_{D_0\cap D_1}$ that can be used to
determine its rank.

\begin{construction}
 First, recall that we can form the divisor
 \[ D_0+D_1=\spec(\OG/f_0f_1)=\spec(\OS\psb{x}/f_0(x)f_1(x)). \]
 This contains $D_0$ and $D_1$, so we have a pullback square of closed
 inclusions as shown on the left below.  This gives a pushout square
 of $\OS$-algebras as shown on the right.
 \begin{diag}
  \node{D_0\cap D_1}   \arrow{s,V} \arrow{e,V}
  \node{D_0}           \arrow{s,V}
  \node{\O_{D_0\cap D_1}} 
  \node{\O_{D_0}}                  \arrow{w,A} \\
  \node{D_1}                       \arrow{e,V}
  \node{D_0+D_1}
  \node{\O_{D_1}}     \arrow{n,A}
  \node{\O_{D_0+D_1},}\arrow{n,A} \arrow{w,A}
 \end{diag}%
which gives a presentation
 \[ \O_{D_0+D_1} \xra{} \O_{D_0}\op\O_{D_1} \xra{} \O_{D_0\cap D_1}. \]
 Explicitly, this is just the presentation
 \[ \OG/(f_0f_1) \xra{\phi} \OG/f_0\op\OG/f_1 \xra{\psi} \OG/(f_0,f_1) \]
 given by 
 \begin{align*}
  \phi(g\bmod{f_0f_1}) &= (g\bmod{f_0},-g\bmod{f_1}) \\
  \psi(g_0\bmod{f_0},g_1\bmod{f_1}) &= g_0+g_1\bmod{(f_0,f_1)}. 
 \end{align*}
\end{construction}

Although this is probably the most natural presentation, it is not
easy to write down the effect of $\phi$ on the obvious bases of
$\OG/(f_0f_1)$ and $\OG/f_i$.  To remedy this, we give an alternate
presentation.  

\begin{construction}
 Let $J_i$ be the ideal generated by $f_i$ and put
 $J=J_0J_1$.  Then $J_i/J$ is free over $\O_S$ with basis
 $\{x^jf_i(x)\st 0\leq j<d_{1-i}\}$ and the inclusion maps
 $J_i\xra{}\OG$ give rise to a presentation
 \[ J_0/J \op J_1/J \xra{\zt} \OG/J=\O_{D_0+D_1} \xra{\xi}
     \OG/(J_0+J_1) = \O_{D_0\cap D_1}.
 \]
 Let $c_{ij}$ be the coefficient of $x^{d_i-j}$ in $f_i(x)$, so that
 $c_{i0}=1$ and $f_i(x)=\sum_{d_i=j+k}c_{ij}x^k$.  Then 
 \begin{align*}
  \zt(x^jf_0(x),0) &= \sum_{k=j}^{d_0+j} c_{0,d_0+j-k} x^k
                      \qquad \text{ for } 0\leq j<d_1 \\
  \zt(0,x^jf_1(x)) &= \sum_{k=j}^{d_1+j} c_{1,d_1+j-k} x^k
                      \qquad \text{ for } 0\leq j<d_0, 
 \end{align*}
 and this tells us the matrix for $\zt$ in terms of the obvious bases
 of $J_0/J\op J_1/J$ and $\OG/J$.  For example, if $d_0=2$ and $d_1=3$
 the matrix is
 \[ \left(\begin{array}{ccc|cc}
  c_{02} & 0      & 0      & c_{13} & 0      \\
  c_{01} & c_{02} & 0      & c_{12} & c_{13} \\
  1      & c_{01} & c_{02} & c_{11} & c_{12} \\
  0      & 1      & c_{01} & 1      & c_{11} \\
  0      & 0      & 1      & 0      & 1
  \end{array}\right)
 \]
 In general, we have a square matrix with $d_0+d_1$ rows and columns.
 The left hand block consists of $d_1$ columns, each of which contains
 $d_1-1$ zeros.  The right hand block consists of $d_0$ columns, each
 of which contains $d_0-1$ zeros.  Clearly $\Int_r(D_0,D_1)$ is the
 closed subscheme defined by the vanishing of the minors of this matrix
 of size $d_0+d_1-r+1$.  In particular, $\Int_1(D_0,D_1)$ is defined by
 the vanishing of the determinant of the whole matrix, which is
 classically known as the \emph{resultant} of $f_0$ and $f_1$.  If
 $f_0(x)=\prod_i(x-a_i)$ and $f_1(x)=\prod_j(x-b_j)$ then the resultant
 is just $\prod_{i,j}(a_i-b_j)$.  We do not know of any similar formula
 for the other minors.
\end{construction}

\begin{construction}
 For a smaller but less symmetrical presentation, we can just use the
 sequence $J_1/J\xra{}\OG/J_0\xra{}\OG/(J_0+J_1)$ induced by the
 inclusion of $J_1$ in $\OG$.  This is isomorphic to the presentation
 $\OG/J_0\xra{\mu_1}\OG/J_0\xra{}\OG/(J_0+J_1)$, where
 $\mu_1(g)=f_1g$.  However, the isomorphism depends on a choice of
 coordinate on $\GG$ (because the element $f_1$ does), so the previous
 presentation is sometimes preferable.  There is of course a similar
 presentation $\OG/J_1\xra{\mu_0}\OG/J_1\xra{}\OG/(J_0+J_1)$.
\end{construction}

Finally, we give a presentation that depends only on the formal
Laurent series $f_0/f_1$ and thus makes direct contact with the
classical Thom-Porteous formula.  
\begin{construction}
 Write $\CM_G=R\lsb{x}=\OG[x^{-1}]$.  Note that $f_1(x)/x^{d_1}$ is a
 polynomial in $x^{-1}$ whose constant term is $1$ and whose other
 coefficients are nilpotent, so it is a unit in $R[x^{-1}]$.  It
 follows that $f_1$ is a unit in $R\lsb{x}$.  Put
 $Q=x^{-1}R[x^{-1}]\sse R\lsb{x}$, so that $R\lsb{x}=R\psb{x}\op Q$.
 Multiplication by the series $x^{d_1}f_0/f_1$ induces a map
 \[ P_1 = \frac{R\psb{x}}{f_1R\psb{x}} \xra{\phi}
    P_0 = \frac{R\lsb{x}}{x^{d_1}R\psb{x}\op Q}. 
 \]
 We claim that the cokernel of $\phi$ is isomorphic to
 $R\psb{x}/(f_0,f_1)=\O_{D_0\cap D_1}$, so we have yet another
 presentation of this ring.  Indeed, the cokernel of $\phi$ is clearly
 given by $R\lsb{x}/(x^{d_1}f_0f_1^{-1}R\psb{x}+x^{d_1}R\psb{x}+Q)$.
 The element $f_1/x^{d_1}$ is invertible in $R[x^{-1}]$ so it is
 invertible in $R\lsb{x}$ and satisfies $(f_1/x^{d_1})Q=Q$.  Thus,
 multiplication by this element gives an isomorphism 
 \[ \frac{R\lsb{x}}{x^{d_1}f_0f_1^{-1}R\psb{x}+x^{d_1}R\psb{x}+Q}
     \simeq
    \frac{R\lsb{x}}{f_0R\psb{x}+f_1R\psb{x}+Q}.
 \]
 As $R\lsb{x}=R\psb{x}\op Q$, we see that the right hand side is just
 $R\psb{x}/(f_0,f_1)$ as claimed.
 
 The elements $\{1,x,\ldots,x^{d_1-1}\}$ give a basis for both $P_0$
 and $P_1$, and the matrix elements of $\phi$ with respect to these
 bases are just the coefficients of $f_0/f_1$ (suitably indexed).
 More precisely, we have
 \[ f_0/f_1 = x^{d_0-d_1} \sum_{i\geq 0} c_ix^{-i}, \]
 where $c_0=1$ and $c_i$ is nilpotent for $i>0$.  We take $c_i=0$ for
 $i<0$ by convention.  The matrix elements $\Phi_{ij}$ of $\phi$ are
 then given by $\Phi_{ij}=c_{d_0+i-j}$ for $0\leq i,j<d_1$.  For
 example, if $d_0=3$ and $d_1=5$ then the matrix is
 \[ \Phi =
   \left(\begin{array}{ccccc}
     c_3 & c_4 & c_5 & c_6 & c_7 \\
     c_2 & c_3 & c_4 & c_5 & c_6 \\
     c_1 & c_2 & c_3 & c_4 & c_5 \\
      1  & c_1 & c_2 & c_3 & c_4 \\
      0  &  1  & c_1 & c_2 & c_3
    \end{array}\right).
 \]
 Now suppose that our divisors $D_i$ arise in the usual way from
 vector bundles $V_i$ over a stably complex manifold $X$, and we have
 a generic linear map $g\:V_0\xra{}V_1$.  Let $Z_r$ be the locus where
 the rank of $g$ is at most $r$, and let $i\:Z_r\xra{}X$ be the
 inclusion.  Generically, this will be a smooth stably complex
 submanifold of $X$, so we have a class $z_r=i_*[Z_r]\in E^0X$.  The
 Thom-Porteous formula says that $z_r=\det(\Psi_r)$, where $\Psi_r$ is
 the square block of size $d_1-r$ taken from the bottom left of
 $\Phi$.  More explicitly, the matrix elements are
 $(\Psi_r)_{ij}=c_{d_0-k+i-j}$ for $0\leq i,j<d_1-r$.  Clearly
 $\det(\Psi_r)\in I_{d_1-r}(\phi)=I_r(\O_{D_0\cap D_1})$.  If $Z_r$ is
 empty then $z_r=0$.  On the other hand,
 Proposition~\ref{prop-int-equiv} will tell us that $\int(D_0,D_1)>r$
 and so $I_r(\O_{D_0\cap D_1})=0$, so $\det(\Psi_r)=0$, which is
 consistent with the Thom-Porteous formula.  It is doubtless possible
 to prove the formula using the methods of this paper, but we have not
 yet done so.
\end{construction}

\begin{proposition}\label{prop-included}
 We always have $\int(D_0,D_1)\leq\min(d_0,d_1)$ (unless the base
 scheme $S$ is empty).  If $\int(D_0,D_1)=d_0$ then $D_0\leq D_1$, and
 if $\int(D_0,D_1)=d_1$ then $D_1\leq D_0$.
\end{proposition}
\begin{proof}
 The presentation $\O_{D_1}\xra{\mu_0}\O_{D_1}\xra{}\O_{D_0\cap D_1}$
 shows that 
 \[ \int(D_0,D_1)=\rank(\O_{D_0\cap D_1})\leq\rank(\O_{D_1})=d_1. \]
 If this is actually an equality we must have
 $\lm^{d_1-d_1+1}(\mu_0)=0$ or in other words $\mu_0=0$, so
 $f_0=0\pmod{f_1}$, so $D_1\leq D_0$.  The remaining claims follow by
 symmetry. 
\end{proof}

\begin{proposition}\label{prop-dec-int}
 If there is a divisor $D$ of degree $k$ such that $D\leq D_0$
 and $D\leq D_1$, then $\int(D_0,D_1)\geq k$.
\end{proposition}
\begin{proof}
 Clearly $\OD$ is a quotient of the ring $\O_{D_0\cap D_1}$, and it
 is free of rank $k$, so
 $\int(D_0,D_1)=\rank(\O_{D_0\cap D_1})\geq k$.
\end{proof}

\begin{definition}\label{defn-sub}
 Given two divisors $D_0,D_1$, we write $\Sub_r(D_0,D_1)$ for the
 scheme of divisors $D$ of degree $r$ such that $D\leq D_0$ and
 $D\leq D_1$.  The proposition shows that the projection
 $\pi\:\Sub_r(D_0,D_1)\xra{}S$ factors through the closed subscheme
 $\Int_r(D_0,D_1)$.  
\end{definition}

\begin{remark}\label{rem-included}
 Proposition~\ref{prop-included} implies that $\Int_{d_0}(D_0,D_1)$ is
 just the largest closed subscheme of $S$ over which we have
 $D_0\leq D_1$.  From this it is easy to see that
 $\Sub_{d_0}(D_0,D_1)=\Int_{d_0}(D_0,D_1)$.  
\end{remark}

It is natural to expect that the map
$\pi\:\Sub_r(D_0,D_1)\xra{}\Int_r(D_0,D_1)$ should be surjective in
some suitable sense.  Unfortunately this does not work as well as one
might hope: the map $\pi$ is not faithfully flat or even dominant, so
the corresponding ring map $\pi^*$ need not be injective.  However, it
is injective in a certain universal case, as we shall show in
Section~\ref{sec-alg-universal}.

We conclude this section with an example where $\pi^*$ is not
injective.  Let $\GG$ be the additive formal group over the scheme
$S=\spec(\Z[a]/a^2)$.  Let $D_0$ and $D_1$ be the divisors with
equations $x^2-a$ and $x^2$, respectively.  Then 
$\O_{D_0\cap D_1}=\OS[x]/(x^2-a,x^2)=\OS[x]/(a,x^2)$, which is the
cokernel of the map $\mu\:\OS[x]/x^2\xra{}\OS[x]/x^2$ given by
$\mu(t)=at$.  The matrix of $\mu$ is
$\left(\begin{array}{cc}a&0\\0&a\end{array}\right)$ which is clearly
nonzero, but $\lm^2(\mu)=a^2=0$.  It follows that $\int(D_0,D_1)=1$,
so $\Int_1(D_0,D_1)=S$.  However, $\Sub_1(D_0,D_1)$ is just the scheme
$D_0\cap D_1=\spec(\OS[x]/(a,x^2))$, so $\pi^*(a)=0$.

For a topological interpretation, let $V_0$ be the tautological bundle
over $\HP^1=S^4$, and let $V_1$ be the trivial rank two complex
bundle.  If we use the cohomology theory $E^*Y=(H^*Y)[u,u^{-1}]$ (with
$|u|=2$) and let $a$ be the second Chern class of $V_0$ we find that
$E^0X=\Z[a]/a^2$, and the equations of $D(V_0)$ and $D(V_1)$ are $x^2-a$
and $x^2$.  Using the theory to be developed in
Section~\ref{sec-int-bundles} and the calculations of the previous
paragraph, we deduce that $V_0$ and $V_1$ cannot have a common
subbundle of rank one, but there is no cohomological obstruction to
finding a map $f\:V_0\xra{}V_1$ with rank at least $1$ everywhere.
To see that such a map does in fact exist, choose a subspace $W<\H^2$
which is a complex vector space of dimension $2$, but not an
$\H$-submodule.  We can then take the constant bundle with fibre
$\H^2/W$ as a model for $V_1$.  The bundle $V_0$ is by definition a
subbundle of the constant bundle with fibre $\H^2$, so there is an
evident projection map $f\:V_0\xra{}V_1$.  As $W$ is not an
$\H$-submodule, we see that $f$ is nowhere zero and thus has rank at
least one everywhere, as claimed.

\section{Unitary bundles}
\label{sec-unitary-bundles}

In order to compare the constructions of the previous section with
phenomena in topology, we need a topological interpretation of the
exterior powers $\lm^k\OD$ when $D$ is the divisor associated to a
vector bundle.

Let $V$ be a complex vector bundle of dimension $d$ over a space $X$.
We can thus form a bundle $U(V)$ of unitary groups in the evident way
(so $U(V)=\{(x,g)\st x\in X \text{ and } g\in U(V_x)\}$).  The key
point is that $E^*U(V)$ can be naturally identified with
$\lm^*_{E^*X}E^{*-1}PV$ (the exterior algebra over the ring $E^*X$
generated by the module $E^*PV$).  Moreover, we can use the group
structure on $U(V)$ to make $E^*U(V)$ into a Hopf algebra over $E^*X$,
and we can make $\lm^*_{E^*X}E^{*-1}PV$ into a Hopf algebra by
declaring $E^*PV$ to be primitive.  We will need to know that our
isomorphism respects these structures.  All this is of course
well-known when $X$ is a point and $E$ represents ordinary cohomology.
Kitchloo~\cite{ki:css} has shown that if one chooses the right proof
then the restriction on $E$ can be removed.  With just a few more
words, we will be able to remove the restriction on $X$ as well.

We start by comparing $U(V)$ with a suitable classifying space.  First
let $V$ be a vector space rather than a bundle.  We let $EU(V)$ denote
the geometric realisation of the simplicial space
$\{U(V)^{n+1}\}_{n\geq 0}$ and we put $BU(V)=EU(V)/U(V)$, which is the
usual simplicial model for the classifying space of $U(V)$.  There is
a well-known map $\eta\:U(V)\xra{}\Om BU(V)$, which is a weak
equivalence of $H$-spaces.  By adjunction we have a map 
$\zt\:\Sg U(V)\xra{}BU(V)$, which gives a map
\[ \zt^* \: \tE^*BU(V) \xra{} \tE^*\Sg U(V)=\tE^{*-1}U(V). \]
The fact that $\eta$ is an $H$-map means
that $\zt$ is primitive, or in other words that
\[ \zt\circ\mu=\zt\circ(\pi_0+\pi_1)\in [\Sg U(V)^2,BU(V)]. \]
We can also construct a tautological bundle $T=EU(V)\tm_{U(V)}V$ over
$BU(V)$.

We now revert to the case where $V$ is a vector bundle over a space
$X$, and perform all the above constructions fibrewise.  Firstly, we
construct the bundle 
$BU(V)=\{(x,e)\st x\in X\text{ and }e\in BU(V_x)\}$.  Note that each
space $BU(V_x)$ has a canonical basepoint, and using these we get an
inclusion $X\xra{}BU(V)$.  

A slightly surprising point is that there is a canonical homotopy
equivalence $BU(V)\xra{}X\tm BU(d)$.  Indeed, we can certainly perform
the definition of $T$ fibrewise to get a tautological bundle over
$BU(V)$, which is classified by a map $q\:BU(V)\xra{}BU(d)$, which is
unique up to homotopy.  We can combine this with the projection
$p\:BU(V)\xra{}X$ to get a map $f=(p,q)\:BU(V)\xra{}X\tm BU(d)$.  The
map $p$ is a fibre bundle projection, and the restriction of $q$ to
each fibre of $p$ is easily seen to be an equivalence.  It is now an
easy exercise with the homotopy long exact sequence of $p$ to see that
$f$ is a weak equivalence.  (Nothing untoward happens with $\pi_0$ and
$\pi_1$ because $BU(d)$ is simply connected.)  

\begin{remark}\label{rem-U-twisted}
 Let $q_0\:X\xra{}BU(d)$ be the restriction of $q$.  Then $q_0$
 classifies the bundle $T|_X\simeq V$, so in general it will be an
 essential map.  Thus, if we just use the basepoint of $BU(d)$ to make
 $X\tm BU(d)$ into a based space over $X$, then our equivalence
 $f\:BU(V)\simeq X\tm BU(d)$ does not preserve basepoints, and cannot
 be deformed to do so.  If it did preserve basepoints we could apply
 the fibrewise loop functor $\Om_X$ and deduce that
 $U(V)\simeq X\tm U(d)$, but this is false in general.
\end{remark}

It follows from the above that $E^*BU(V)$ is a formal power series
algebra over $E^*X$, generated by the Chern classes of $T$.  It will
be convenient for us to modify this description slightly by
considering the virtual bundle $T-V$ (where $V$ is implicitly pulled
back to $BU(V)$ by the map $p\:BU(V)\xra{}X$).  We have
$f_T(t)=t^d\sum_{k=0}^da_kt^{-k}$ and
$f_V(t)=t^d\sum_{k=0}^db_kt^{-k}$ for some coefficients $a_k\in
E^0BU(V)$ and $b_k\in E^0X$ so 
$f_{T-V}(t)=f_T(t)/f_V(t)=\sum_{k\geq 0}c_kt^{-k}$ for some $c_k\in
E^0BU(V)$.  For $k\leq d$ we have $c_k=a_k\pmod{b_1,\ldots,b_d}$ and
it follows easily that
\[ E^*BU(V) = (E^*X)\psb{c_1,\ldots,c_d}. \]
Note that the restriction of $T-V$ to $X\subset BU(V)$ is trivial, so
the classes $c_k$ restrict to zero on $X$.

Next, consider the fibrewise suspension $\Sg_XU(V)$.  By dividing each
fibre into two cones we obtain a decomposition $\Sg_XU(V)=C_0\cup C_1$
where the inclusion of $X$ in each $C_i$ is a homotopy equivalence,
and $C_0\cap C_1=U(V)$.  Using a Mayer-Vietoris sequence we deduce
that $\tE_X^*\Sg_XU(V)\simeq\tE^{*-1}U(V)$ and that this can be
regarded as an ideal in $E^*\Sg_XU(V)$ whose square is zero.
Moreover, the construction of $\zt$ can be carried out fibrewise to
get a map $\Sg_XU(V)\xra{}BU(V)$ which is again primitive.  It follows
that $\zt$ induces a map
\[ \zt^*\:\Ind(E^*BU(V)) \xra{} \Prim(E^{*-1}U(V)). \]
(Here $\Ind$ and $\Prim$ denote indecomposables and primitives over
$E^*X$.)  Note also that $\Ind(E^*BU(V))$ is a free module over $E^*X$
generated by $\{c_1,\ldots,c_d\}$.

To prove that $\zt^*$ is injective, we need to consider the complex
reflection map $\rho\:\Sg_XPV_{+X}\xra{}U(V)$, which we define as
follows.  For $t\in S^1=\R\cup\{\infty\}$ and $x\in X$ and 
$L\in PV_x$, the map $\rho(t,x,L)$ is the endomorphism of $V_x$ that
has eigenvalue $\gm^{-1}(t)$ on the line $L$, and eigenvalue $1$ on
$L^\perp$.  Here $\gm^{-1}(t)=(it+1)/(it-1)\in U(1)$, as in
Section~\ref{subsec-spheres}.  Using this we obtain a map
$\xi=\zt\circ\Sg_X\rho\:\Sg^2_XPV_{+X}\xra{}BU(V)$.

Our next problem is to identify the virtual bundle $\xi^*(T-V)$ over
$\Sg_X^2PV_{+X}$.  For this it is convenient to identify $S^2$ with
$\CP^1$ and thus $\Sg^2PV_{+X}$ with a quotient of $\CP^1\tm PV$.  We
have tautological bundles $H$ and $L$ over $\CP^1$ and $PV$, whose
Euler classes we denote by $y$ and $x$.
\begin{lemma}
 We have $\xi^*(T-V)\simeq (H-1)\ot L$.  Moreover, there is a power
 series $g(s)\in E^0\psb{s}$ with $g(0)=1$ such that
 $\xi^*c_k=-yx^{k-1}g(x)$ for $k=1,\ldots,d$.  (If $E^0$ is
 torsion-free then $g(s)=1/\log'_F(x)$.)
\end{lemma}
\begin{proof}
 In the proof it will be convenient to write $T_V$ and $L_V$ instead
 of $T$ and $L$, to display the dependence on $V$.

 First consider the case where $X$ is a point and $V=\C$.  Then
 $\rho\:S^1\xra{}U(1)=U(\C)$ is a homeomorphism and
 $BU(\C)\simeq\CPi$.  It is a standard fact that
 $\xi\:S^2\xra{}BU(\C)$ can be identified with the inclusion
 $\CP^1\xra{}\CPi$, and thus that $\xi^*T_\C=H$.

 In the general case, note that we have a map
 $\xi_L\:\CP^1\tm PL\xra{}BU(L)$ of spaces over $PV$.  The projection
 $PL\xra{}PV$ is a homeomorphism which we regard as the identity.  If
 we let $\pi\:PV\xra{}X$ be the projection, we have a splitting
 $\pi^*V=L\op(\pi^*V\om L)$.  The inclusion $L\xra{}\pi^*V$ gives an
 inclusion $U(L)\xra{}\pi^*U(V)$ and thus an inclusion
 $BU(L)\xra{}\pi^*BU(V)$, or equivalently a map
 $\phi\:BU(L)\xra{}BU(V)$ covering $\pi$.  As $T_V=V\tm_{U(V)}EU(V)$
 and $U(L)$ acts trivially on $\pi^*V\om L$ we see that
 $\phi^*T_V=T_L\op(\pi^*V\om L)$.

 Next, we note that tensoring with $L$ gives an isomorphism
 $\tau\:U(\C)\tm PV\xra{}U(L)$ and thus an isomorphism
 $B\tau\:BU(\C)\tm PV\xra{}BU(L)$ with $(B\tau)^*T_L=T_\C\ot L$.

 One can check that the following diagram commutes:
 \begin{diag}
  \node{\CP^1\tm PV}
   \arrow{e,t}{1}
   \arrow{s,l}{\xi_\C\tm 1}
  \node{\CP^1\tm PV}
   \arrow{e,t}{1}
   \arrow{s,l}{\xi_L}
  \node{\CP^1\tm PV}
   \arrow{s,r}{\xi_V} \\
  \node{BU(\C)\tm PV}
   \arrow{e,tb}{\simeq}{B\tau}
  \node{BU(L)}
   \arrow{e,b}{\phi}
  \node{BU(V).}
 \end{diag}%
It follows that $\xi_V^*T_V\simeq(\xi_\C\tm 1)^*(B\tau)^*\phi^*T_V$, 
 and the previous discussion identifies this with
 $(H\ot L)\op(\pi^*V\om L)$.  It follows that
 $\xi_V^*(T_V-V)\simeq(H\ot L)-L=(H-1)\ot L$, as claimed.

 Now let $g(s)$ be the partial derivative of $t+_Fs$ with respect to
 $t$ evaluated at $t=0$.  This is characterised by the equation
 $t+_Fs=tg(s)+s\pmod{t^2}$; it is clear that $g(0)=1$, and by applying
 $\log_F$ we see that $g(s)=1/\log'_F(s)$ in the torsion-free case.
 As $y^2=0$ we see that the Euler class of $H\ot L$ is
 $x+_Fy=x+yg(x)$.  Thus, we have
 \begin{align*}
  f_{H\ot L - L}(t) &= (t-x-yg(x))/(t-x) \\
                    &= 1 - yg(x)t^{-1}/(1-x/t) \\
                    &= 1 - \sum_{k>0} y g(x) x^{k-1} t^{-k}.
 \end{align*}
 The $k$'th Chern class of $(H-1)\ot L$ is the coefficient of $t^{-k}$
 in this series, which is $-yg(x)x^{k-1}$ as claimed.
\end{proof}
\begin{corollary}\label{cor-xi-iso}
 The induced map
 $\xi^*\:\Ind(E^*BU(V))\xra{}E^*(\Sg^2_XPV_{+X},X)=E^{*-2}PV$ is an
 isomorphism.  \qed
\end{corollary}

\begin{theorem}\label{thm-EUV}
 There is a natural isomorphism $\lm^*E^{*-1}PV\xra{}E^*U(V)$ of Hopf
 algebras over $E^*X$.
\end{theorem}
\begin{proof}
 Put $a_i=\xi^*c_i\in\Prim(E^*U(V))$ for $i=1,\ldots,d$.  Given a
 sequence $I=(i_1,\ldots,i_r)$ with $1\leq i_1<\ldots<i_r\leq d$, put
 $a_I=\prod_ja_{i_j}$.  We first claim that the elements $a_I$ form a
 basis for $E^*U(V)$ over $E^*X$.  This is very well-known in the case
 where $X$ is a point (so $U(V)\simeq U(d)$) and $E$ represents
 ordinary cohomology; it can proved using the Serre spectral sequence
 of the fibration $U(d-1)\xra{}U(d)\xra{}S^{2d-1}$.  For a more
 general theory $E$ we still have an Atiyah-Hirzebruch-Serre spectral
 sequence $H^p(S^{2d-1};E^qU(d-1))\convto E^{p+q}U(d)$.  It follows
 easily that the elements $a_I$ form a basis whenever $X$ is a point.
 A standard argument now shows that they form a basis for any $X$.
 Indeed, it follows easily from the above that they form a basis
 whenever $V$ is trivialisable.  We can give $X$ a cell structure such
 that $V$ is trivialisable over each cell, and then use Mayer-Vietoris
 sequences to check that the elements $a_I$ form a basis whenever $X$
 is a finite complex.  Finally, we can use the Milnor exact sequence
 to show that the elements $a_I$ form a basis for all $X$.

 The ring $E^*U(V)$ is graded-commutative so we certainly have
 $a_ia_j=-a_ja_i$ and in particular $2a_i^2=0$ for all $i$.  Suppose
 we can show that $a_i^2=0$.  Then $\zt^*$ extends to give a map
 $\lm_{E^*X}^*\Ind(E^*BU(V))\xra{}E^{*-1}U(V)$ of Hopf algebras, and
 from the previous paragraph we see that this is an isomorphism.
 Combining this with the isomorphism of Corollary~\ref{cor-xi-iso}
 gives the required isomorphism $\lm^*E^{*-1}PV\xra{}E^*U(V)$.

 All that is left is to check that $a_i^2=0$.  For this we consider
 the case of the tautological bundle $T$ over $BU(d)$, and take
 $E=MP=MU[u,u^{-1}]$.  (We use this $2$-periodic version of $MU$
 simply to comply with our standing assumptions on $E$; we could
 equally well use $MU$ itself.)  Here it is standard that $MP^*BU(d)$
 is a formal power series algebra over $MP^*$ and thus is
 torsion-free.  The ring $MP^*U(T)$ is a free module over $MP^*BU(d)$
 and thus is also torsion-free.  As $2a_i^2=0$ we must have $a_i^2=0$
 as required.  More generally, for an arbitrary bundle $V$ over a
 space $X$ we have a classifying map $X\xra{}BU(d)$ giving rise to a
 map $U(V)\xra{}U(T)$.  Moreover, for any $E$ we can choose an
 orientation in degree zero and thus a ring map $MP\xra{}E$.  Together
 these give a ring map $MP^*U(T)\xra{}E^*U(V)$, which carries $a_i$ to
 $a_i$.  As $a_i^2=0$ in $MP^*U(T)$, the same must hold in $E^*U(V)$.
\end{proof}

We will need to extend the above result slightly to give a topological
interpretation of the quotient rings 
\[ \lm^{\leq r}E^{*-1}PV=\lm^*E^{*-1}PV/\lm^{>r}E^{*-1}PV. \]
For this we recall Miller's filtration of $U(V)$:
\begin{align*}
 F_kU(V) &= \{g\in U(V) \st \codim(\ker(g-1))\leq k\} \\
         &= \{g\in U(V) \st \rank(g-1)\leq k\}.
\end{align*}
More precisely, this is supposed to be interpreted fibrewise, so
\[ F_kU(V)=\{(x,g)\st x\in X\text{ and } g\in U(V_x) \text{ and }
              \rank(g-1)\leq k\}.
\]
It is not hard to see that $\rho$ gives a homeomorphism
$\Sg_XPV_{+X}\xra{}F_1U(V)$.  It is known from work of
Miller~\cite{mi:sss} that when $X$ is a point, the filtration is
stably split.  Crabb showed in~\cite{cr:ssu} that the splitting works
fibrewise; our outline of related material essentially follows his
account.  

We will need to recall the basic facts about the quotients in Miller's
filtration.  Consider the space
\[ G_k(V)=\{(x,W)\st x\in X\;,\; W\leq V_x \;,\; \dim(W)=k\}. \]
For each point $(x,W)\in G_k(V)$ we have a Lie group $U(W)$ and its
associated Lie algebra $\uu(W)=\{\al\in\End(W)\st\al+\al^*=0\}$.
These fit together to form a bundle over $G_k(V)$ which we denote by
$\uu$.  Given a point $(x,W,\al)$ in the total space of this bundle
one checks that $\al-1$ is invertible and that
$g:=(\al+1)(\al-1)^{-1}$ is a unitary automorphism of $W$ without
fixed points, so $g\op 1_{W^\perp}\in F_kU(V_x)\sm F_{k-1}U(V_x)$.  It
is not hard to show that this construction gives a homeomorphism of
the total space of $\uu$ with $F_kU(V)\sm F_{k-1}U(V)$ and thus a
homeomorphism of the Thom space $G_k(V)^\uu$ with
$F_kU(V)/F_{k-1}U(V)$.

If $g\in F_jU(V_x)$ and $h\in F_kU(V_x)$ then $\ker(g-1)\cap\ker(h-1)$
has codimension at most $j+k$, so $gh\in F_{j+k}U(V)$, so the
filtration is multiplicative.  A less obvious argument shows that it
is also comultiplicative, up to homotopy:
\begin{lemma}\label{lem-diagonal}
 The diagonal map $\dl\:U(V)\xra{}U(V)\tm_XU(V)$ is homotopic to a
 filtration-preserving map.
\end{lemma}
\begin{proof}
 For notational convenience, we will give the proof for a vector
 space; it can clearly be done fibrewise for vector bundles.

 We regard $U(1)$ as the set of unit complex numbers and define
 $p_0,p_1\:U(1)\xra{}U(1)$ as follows:
 \begin{align*}
  p_0(z) &= \begin{cases} z^2 & \text{ if } \text{Im}(z)\geq 0 \\
                           1   & \text{ otherwise }
            \end{cases} \\
  p_1(z) &= \begin{cases} z^2 & \text{ if } \text{Im}(z)\leq 0 \\
                           1   & \text{ otherwise. }
            \end{cases} \\
 \end{align*}
 Thus $(p_0,p_1)\:U(1)\xra{}U(1)\tm U(1)$ is just the usual pinch map
 $U(1)\xra{}U(1)\Wedge U(1)\subset U(1)\tm U(1)$.

 Note that if $g\in U(V)$ and $r\in\{0,1\}$ then the eigenvalues of
 $g$ lie in $U(1)$ so we can interpret $p_r(g)$ as an endomorphism of
 $V$ as in Appendix~\ref{apx-functional}.  As $p_r(U(1))\sse U(1)$ we
 see that $\ov{p_r(z)}=p_r(z)^{-1}$ for all $z\in U(1)$ and thus that
 $p_r(g)^*=p_r(g)^{-1}$, so $p_r$ gives a map from $U(V)$ to itself.

 We now define $\dl'\:U(V)\xra{}U(V)\tm U(V)$ by
 $\dl'(g)=(p_0(g),p_1(g))$.  It is clear that the
 filtration of $p_0(g)$ is the number of eigenvalues of $g$ (counted
 with multiplicity) lying in the open upper half-circle, and the
 filtration of $p_1(g)$ is the number in the open lower half-circle.
 Thus, the filtration of $\dl'(g)$ is the number of eigenvalues not
 equal to $\pm 1$, which is less than or equal to the filtration of
 $g$. 

 On the other hand, each map $p_r\:U(1)\xra{}U(1)$ has degree $1$ and
 thus is homotopic to the identity, so $\dl'$ is homotopic to $\dl$.
\end{proof}

\begin{theorem}\label{thm-E-Fk}
 There is a natural isomorphism
 $\lm_{E^*X}^{<k}E^{*-1}PV\xra{}E^*F_{k-1}U(V)$. 
\end{theorem}
\begin{proof}
 For brevity we write $\lm^k=\lm^k_{E^*X}E^{*-1}PV$.  We also write
 $\lm^*=\bigoplus_k\lm^k$ and $\lm^{\geq k}=\bigoplus_{j\geq k}\lm^j$
 and $\lm^{<k}=\lm^*/\lm^{\geq k}=\bigoplus_{j<k}\lm^j$.

 Because the filtration of $U(V)$ is stably split, the restriction map
 $\lm^*E^{*-1}PV=E^*U(V)\xra{}E^*F_{k-1}U(V)$ is a split surjection,
 with kernel $J_k$ say.  Note that $\lm^*/J_k$ and $J_k$ are both
 projective over $E^*X$.  We need to show that $J_k=\lm^{\geq k}$.

 First, we have $F_0U(V)=X$ and it follows easily that
 $J_1=\lm^{\geq 1}$. 

 We next claim that $J_jJ_k\leq J_{j+k}$ for all $j,k$.  Indeed, $J_j$
 is the image in cohomology of the map $U(V)\xra{}U(V)/F_{j-1}$, and
 so $J_jJ_k$ is contained in the image in cohomology of the map
 \[ \phi = (U(V) \xra{\dl} U(V)\tm_XU(V) \xra{}
     U(V)/F_{j-1} \Smash_X U(V)/F_{k-1}).
 \]
 Note that $\dl$ is homotopic to the map $\dl'$, which sends
 $F_{j+k-1}$ into $F_{j-1}\tm_XU(V)\cup U(V)\tm_XF_{k-1}$.  It follows
 that the restriction of $\phi$ to $F_{j+k-1}$ is null, and thus that
 $J_jJ_k\leq J_{j+k}$ as claimed.  It follows inductively that
 $\lm^{\geq k}\leq J_k$ for all $k$.  This gives us a natural
 surjective map $\lm^{<k}\xra{}E^*F_{k-1}U(V)$.

 We previously gave a natural basis $\{a_I\}$ for $\lm^*$, and it is
 clear that the subset $\{a_I\st\; |I|<k\}$ is a basis for $\lm^{<k}$.
 It will be enough to prove that the images of these form a basis for
 $E^*F_{k-1}U(V)$.  The argument of Theorem~\ref{thm-EUV} allows us to
 reduce to the case where $X$ is a point, $V=\C^d$, and $E$ represents
 ordinary cohomology.  A proof in this case has been given by
 Kitchloo~\cite{ki:css} (and possibly by others) but we will sketch an
 alternate proof for completeness.  As the map
 $\lm^{<k}\xra{}H^*F_{k-1}U(d)$ is surjective, it will suffice to show
 that the source and target have the same rank as free Abelian groups.
 For this, it will suffice to show that $\lm^jH^*\CP^{d-1}$ has the
 same rank as $\tH^*(F_jU(d)/F_{j-1}U(d))$ for $0\leq j\leq d$.  As
 $H^*\CP^{d-1}$ has rank $d$, it is clear that $\lm^jH^*\CP^{d-1}$ has
 rank $\bcf{d}{j}$.  On the other hand $F_jU(d)/F_{j-1}U(d)$ is the
 Thom space $G_j(\C^d)^\uu$.  Note that although $\uu$ is not a
 complex bundle, it is necessarily orientable because $G_j(\C^d)$ is
 simply connected.  Thus, the Thom isomorphism theorem tells us that
 the rank of $\tH^*G_j(\C^d)^\uu$ is the same as that of
 $H^*G_j(\C^d)$.  By counting Schubert cells we see that this is again
 $\bcf{d}{j}$, as required.  (This will also follow from
 Proposition~\ref{prop-Gr}.)
\end{proof}

\section{Intersections of bundles}
\label{sec-int-bundles}

Let $X$ be a space, and let $V_0$ and $V_1$ be complex vector bundles
over $X$.  In Section~\ref{sec-int-div} we defined divisors
$D(V_i)=(PV_i)_E$ on $\GG$ over $X_E$, and we also defined the
intersection index $\int(V_0,V_1)$.

\begin{theorem}\label{thm-int}
 We have $\int(V_0,V_1)\leq\int(D(V_0),D(V_1))$.
\end{theorem}
\begin{proof}
 Suppose we have isometric linear embeddings
 $V_0\xra{j_0}W\xla{j_1}V_1$ such that
 $\dim((j_0V_{0x})\cap(j_1V_{1x}))\geq r$ for all $x$.  We must show
 that $\rank(\OG/(f_{V_0},f_{V_1}))\geq r$.  Put $d_i=\dim(V_i)$ and
 $e=\dim(W)$.  Recall that $E^0PV_i=\OG/f_{V_i}$ and that
 $E^0PW=\OG/f_W$.  As each $V_i$ embeds in $W$ we see that $f_{V_i}$
 divides $f_W$ and there is a natural surjection $E^0PW\xra{}E^0PV_i$.
 By combining these maps we get a map $\phi\:E^0PW\xra{}E^0PV_0\op
 E^0PV_1$, whose cokernel is $\OG/(f_{V_0},f_{V_1})$.  From the
 definition of the Fitting rank, we must prove that
 $\lm^{d_0+d_1-r+1}\phi=0$.

 For this, we first note that an isometric embedding $j\:V\xra{}W$ of
 vector spaces gives rise to a homomorphism $j_*\:U(V)\xra{}U(W)$ by 
 \[ j_*(g) = jgj^{-1}\op 1_{jV^\perp} \: W=jV\op jV^\perp\xra{}W. \]
 The alternative description $j_*(g)=jgj^*+1-jj^*$ makes it clear that
 $j_*(g)$ depends continuously on $j$ and $g$.

 We now extend this definition fibrewise, and define
 $\gm\:U(V_0)\tm_X U(V_1)\xra{}U(W)$ by
 $\gm(g_0,g_1)=(j_{0*}g_0)(j_{1*}g_1)$.  We have
 $E^*U(W)=\lm^*E^{*-1}PW$ and 
 \begin{align*}
  E^*U(V_0)\tm_XU(V_1) &= E^*U(V_0)\ot_{E^*X}E^*U(V_1) \\
                       &= \lm^*E^{*-1}PV_0\ot_{E^*X}\lm^*E^{*-1}PV_1\\
                       &= \lm^*(E^{*-1}PV_0\op E^{*-1}PV_1).
 \end{align*}
 Using the fact that $E^{*-1}PW$ is primitive in $E^*U(W)$, we find
 that $\gm^*=\lm^*\phi$.  Next, observe that if $g_i\in U(V_{ix})$ for
 $i=0,1$ we have 
 \[ \gm(g_0,g_1)\in U(j_0V_{0x}+j_1V_{1x})\sse U(W) \]
 and $\dim(j_0V_{0x}+j_1V_{1x})\leq d_0+d_1-r$ so
 $\gm(g_0,g_1)\in F_{d_0+d_1-r}U(W)$.  Thus $\gm$ factors through
 $F_{d_0+d_1-r}U(W)$, and it follows that $\lm^{d_0+d_1-r+1}E^{*-1}PW$
 is mapped to zero by $\gm^*$, as required.
\end{proof}

As an addendum, we show that some natural variations of the definition
of intersection index do not actually make a difference.
\begin{lemma}\label{lem-isom}
 Let $V$ and $W$ be vector bundles over a space $X$, and let
 $j\:V\xra{}W$ be a linear embedding.  Then $j$ is an isometric
 embedding if and only if $j^*j=1$ (where $j^*$ is the adjoint of
 $j$).  In any case, there is a canonical isometric embedding
 $\jh\:V\xra{}W$ with the same image as $j$.
\end{lemma}
\begin{proof}
 If $j^*j=1$ then $\|jv\|^2=\ip{jv,jv}=\ip{v,j^*jv}=\ip{v,v}=\|v\|^2$,
 so $j$ is an isometry.  Conversely, if $j$ is an isometry then it
 preserves inner products so $\ip{v',j^*jv}=\ip{jv',jv}=\ip{v',v}$ for
 all $v,v'$ which means that $j^*jv=v$.

 Even if $j$ is not an isometry we have $\ip{v,j^*jv}=\|jv\|^2$ which
 implies that $j^*j$ is injective.  It is thus a strictly positive
 self-adjoint operator on $V$, so we can define $(j^*j)^{-1/2}$ by
 functional calculus (as in Appendix~\ref{apx-functional}).  We then
 define $\jh=j\circ(j^*j)^{-1/2}$.  This is the composite of $j$
 with an automorphism of $V$, so it has the same image as $j$.  It
 also satisfies $\jh^*\jh=1$, so it is an isometric
 embedding. 
\end{proof}

\begin{proposition}\label{prop-int-equiv}
 Let $V_0$ and $V_1$ be bundles over a space $X$.  Consider the
 following statements:
 \begin{itemize}
  \item[\rm(a)] There exists a bundle $V$ of dimension $k$ and linear
   isometric embeddings $V_0\xla{i_0}V\xra{i_1}V_1$.
  \item[\rm(a$'$)] There exists a bundle $V$ of dimension $k$ and linear
   embeddings $V_0\xla{i_0}V\xra{i_1}V_1$.
  \item[\rm(b)] There exist a bundle $W$ and isometric linear embeddings
   $V_0\xra{j_0}W\xla{j_1}V_1$ such that
   $\dim((j_0V_{0x})\cap(j_1V_{1x}))\geq k$ for all $x\in X$.
  \item[\rm(b$'$)] There exist a bundle $W$ and linear embeddings
   $V_0\xra{j_0}W\xla{j_1}V_1$ such that
   $\dim((j_0V_{0x})\cap(j_1V_{1x}))\geq k$ for all $x\in X$.
  \item[\rm(c)] There is a linear map $f\:V_0\xra{}V_1$ such that
   $\rank(f_x)\geq k$ for all $x\in X$.
 \end{itemize}
 Then (a)$\iffa$(a$'$)$\Rightarrow$(b)$\iffa$(b$'$)$\iffa$(c).
\end{proposition}
\begin{proof}
 It follows immediately from Lemma~\ref{lem-isom} that
 (a)$\iffa$(a$'$) and (b)$\iffa$(b$'$).

 (a)$\Rightarrow$(b): Define $W$, $j_0$ and $j_1$ by the following
 pushout square:
 \begin{diag}
  \node{V}   \arrow{e,t,V}{i_0} \arrow{s,l,V}{i_1} 
  \node{V_0}                    \arrow{s,r,V}{j_0} \\
  \node{V_1} \arrow{e,b,V}{j_1} \node{W.}
 \end{diag}%
Equivalently, we can write $V'_t$ for the orthogonal complement of
 $i_tV$ in $V_t$ and then $W=V\op V'_0\op V'_1$.

 (b)$\Rightarrow$(c): Put $f=j_1^*j_0\:V_0\xra{}V_1$.  By hypothesis,
 for each $x$ we can choose an orthonormal sequence $u_1,\ldots,u_k$
 in $(j_0V_{0x})\cap(j_1V_{1x})$.  We can then choose elements
 $v_p\in V_{0x}$ and $w_p\in V_{1x}$ such that $u_p=j_0v_p=j_1w_p$.
 We find that
 $\ip{fv_p,w_q}=\ip{j_0v_p,j_1w_q}=\ip{u_p,u_q}=\dl_{pq}$.  This
 implies that the elements $fv_1,\ldots,fv_k$ are linearly
 independent, so $\rank(f)\geq k$ as required.

 (c)$\Rightarrow$(b): Note that $f_x^*f_x\:V_{0x}\xra{}V_{0x}$ is a
 nonnegative self-adjoint operator with the same kernel as $f_x$, and
 thus the same rank as $f_x$.  Similarly, $f_xf_x^*$ is a nonnegative
 self-adjoint operator on $V_{1x}$ with the same rank as $f_x$.  More
 basic facts about these operators are recorded in
 Proposition~\ref{prop-square}. 

 As in Definition~\ref{defn-ek} we let $\lm_j=e_j(f_x^*f_x)$ be the
 $j$'th eigenvalue of $f_x^*f_x$ (listed in descending order and
 repeated according to multiplicity).  We see from
 Proposition~\ref{prop-ek-cont} that $\lm_j$ is a continuous function
 of $x$.  Moreover, as $f_x^*f_x$ has rank at least $k$ we see that
 $\lm_k>0$.  Now define $\tau_x\:[0,\infty)\xra{}[0,\infty)$ by
 $\tau_x(t)=\max(\lm_k,t)$, and define $\mu_x=\tau_x(f_x^*f_x)$ and
 $\nu_x=\tau(f_xf_x^*)$.  (Here we are using functional calculus as in
 Appendix~\ref{apx-functional} again.)  One checks that
 $f_x\mu_x=\nu_xf_x$ and $\mu_xf_x^*=f_x^*\nu_x$.  We now have maps
 \begin{align*}
          \mu^{1/2} \: V_0 &\xra{} V_0 \\
                  f \: V_0 &\xra{} V_1 \\
  (\mu+f^*f)^{-1/2} \: V_0 &\xra{} V_0, 
 \end{align*}
 which we combine to get a map
 \[ j_0 = (\mu^{1/2},f)\circ(\mu+f^*f)^{-1/2} \: V_0\xra{}V_0\op V_1.
 \] 
 Similarly, we define
 \[ j_1 = (f^*,\nu^{1/2})\circ(\nu+ff^*)^{-1/2} \: V_1\xra{}V_0\op V_1.
 \] 
 It is easy to check that $j_0^*j_0=1$ and $j_1^*j_1=1$, so $j_0$ and
 $j_1$ are isometric embeddings.  
  
 Now choose an orthonormal sequence $v_1,\ldots,v_k$ of eigenvectors
 of $f_x^*f_x$, with eigenvalues $\lm_1,\ldots,\lm_k$.  Put
 $v'_i=f_x(v_i)/\sqrt{\lm_i}\in V_1$; these vectors form an
 orthonormal sequence of eigenvectors of $f_xf_x^*$, with the same
 eigenvalues. 

 For $i\leq k$ we have $\lm_i\geq\lm_k>0$ so $\tau_x(\lm_i)=\lm_i$ so
 $(\mu+f^*f)^{-1/2}(v_i)=v_i/\sqrt{2\lm_i}$ and
 $\mu^{1/2}(v_i)=\sqrt{\lm_i}v_i$ so $j_0(v_i)=(v_i,v'_i)/\sqrt{2}$.
 This is the same as $j_1(v'_i)$, so it lies in
 $(j_0V_{0x})\cap(j_1V_{1x})$.  Thus, this intersection has dimension
 at least $k$, as required.
\end{proof}

We conclude this section with a topological interpretation of the
scheme $D(V_0)\cap D(V_1)$ itself.
\begin{proposition}\label{prop-sphere-bundle}
 Let $V_0$ and $V_1$ be vector bundles over a space $X$, and let $L_0$
 and $L_1$ be the tautological bundles of the two factors in
 $PV_0\tm_XPV_1$.  Then there is a natural map
 $S(\Hom(L_0,L_1))_E\xra{}D(V_0)\cap D(V_1)$, which is an isomorphism
 if the map $E^*P(V_0\op V_1)\xra{}E^*PV_0\op E^*PV_1$ is injective.
\end{proposition}
\begin{proof}
 We divide the sphere bundle $S(V_0\op V_1)$ into two pieces, which
 are preserved by the evident action of $U(1)$:
 \begin{align*}
  C_0 &= \{ (v_0,v_1)\in S(V_0\op V_1) \st \; \|v_0\|\geq\|v_1\|\} \\
  C_1 &= \{ (v_0,v_1)\in S(V_0\op V_1) \st \; \|v_1\|\geq\|v_0\|\}.
 \end{align*}
 The inclusions $V_i\xra{}V_0\op V_1$ give inclusions
 $S(V_i)\xra{}C_i$ which are easily seen to be homotopy equivalences.
 It follows that $C_i/U(1)\simeq PV_i$.  We also have 
 \[ C_0\cap C_1 = \{(v_0,v_1)\st \|v_0\|=\|v_1\|=2^{-1/2}\}
                  \simeq S(V_0)\tm S(V_1).
 \]
 Given a point in this space we have a map $\al\:\C v_0\xra{}\C v_1$
 sending $v_0$ to $v_1$.  This has norm $1$ and is unchanged if we
 multiply $(v_0,v_1)$ by an element of $U(1)$.  Using this we see that
 $(C_0\cap C_1)/U(1)=S(\Hom(L_0,L_1))$.  Of course, we also have
 $(C_0\cup C_1)/U(1)=P(V_0\op V_1)$.  We therefore have a homotopy
 pushout square as shown on the left below, giving rise to a
 commutative square of formal schemes as shown on the right.
 \begin{diag}
  \node{S(\Hom(L_0,L_1))}
   \arrow{e,t}{i_0}
   \arrow{s,l}{i_1}
  \node{PV_0}
   \arrow{s,r}{j_0} 
  \node{S(\Hom(L_0,L_1))_E} 
   \arrow{e}
   \arrow{s}
  \node{D(V_0)} 
   \arrow{s,V} \\
  \node{PV_1}
   \arrow{e,b}{j_1} 
  \node{P(V_0\op V_1)}
  \node{D(V_1)}
   \arrow{e,V}
  \node{D(V_0\op V_1).}
 \end{diag}%
This evidently gives us a map
 $S(\Hom(L_0,L_1))_E\xra{}D(V_0)\cap D(V_1)$.

 To be more precise, we use the Mayer-Vietoris sequence associated to
 our pushout square.  This gives a short exact sequence
 \[ \cok(f^0) \xra{p} E^0S(\Hom(L_0,L_1)) \xra{q} \ker(f^{-1}), \]
 where 
 \[ f^k=(j_0^*,j_1^*)\:E^kP(V_0\op V_1)\xra{}E^kPV_0\op E^kPV_1. \]
 We have seen that $\cok(f^0)=\O_{D(V_0)\cap D(V_1)}$, and the map $p$
 just corresponds to our map
 \[ S(\Hom(L_0,L_1))_E\xra{}D(V_0)\cap D(V_1). \]
 This map will thus be an isomorphism if $f^*$ is injective, as
 claimed.
\end{proof}

\section{Algebraic universal examples}
\label{sec-alg-universal}

Let $\GG$ be a formal group over a formal scheme $S$.  Later we will
work with bundles over a space $X$, and we will take $S=X_E$ and
$\GG=(\CPi\tm X)_E$.  We write
$\Div_d^+=\Div_d^+(\GG)\simeq\GG^d/\Sg_d$, so
$\O_{\Div_d^+}=\OS\psb{c_1,\ldots,c_d}$.

Fix integers $d_0,d_1,r\geq 0$.  We write $\Int_r(d_0,d_1)$ for the
scheme of pairs $(D_0,D_1)$ where $D_0$ and $D_1$ are divisors of
degrees $d_0$ and $d_1$ on $\GG$, and $\int(D_0,D_1)\geq r$.  In other
words, if $D_i$ is the evident tautological divisor over
$\Div_{d_0}^+\tm\Div_{d_1}^+$ then $\Int_r(d_0,d_1)=\Int_r(D_0,D_1)$.
We will assume that $r\leq\min(d_0,d_1)$ (otherwise we would have
$\Int_r(d_0,d_1)=\emptyset$.)

For a more concrete description, put 
\[ R=\O_{\Div_{d_0}^+\tm\Div_{d_1}^+} = 
   \OS\psb{c_{0j}\st j<d_0}\psb{c_{1j}\st j<d_1}.
\]
Let $A$ be the matrix of $\zt$ over $R$ as in
Section~\ref{sec-int-div}, and let $I$ be the ideal in $R$ generated
by the minors of $A$ of size $d_0+d_1-r+1$.  Then
$\Int_r(d_0,d_1)=\spf(R/I)$.

We will also consider a ``semi-universal'' case.  Suppose we have a
divisor $D_1$ on $\GG$ over $S$, with degree $d_1$.  Let $D_0$ be the
tautological divisor over $\Div_{d_0}^+$.  We can regard $D_0$ and
$D_1$ as divisors on $\GG$ over $\Div_{d_0}^+$ and thus form the
closed subscheme $\Int_r(D_0,D_1)\sse\Div_{d_0}^+$.  We denote this
scheme by $\Int_r(d_0,D_1)$.

We can also define schemes $\Sub_r(d_0,d_1)$ and $\Sub_r(d_0,D_1)$ in
a parallel way.  

\begin{remark}\label{rem-sub-univ}
 $\Sub_r(d_0,d_1)$ is just the scheme of triples $(D,D_0,D_1)$ for
 which $D\leq D_0$ and $D\leq D_1$.  This is isomorphic to the scheme
 of triples
 $(D,D'_0,D'_1)\in\Div_r^+\tm\Div_{d_0-r}^+\tm\Div_{d_1-r}^+$, by the
 map $(D,D'_0,D'_1)\mapsto(D,D+D'_0,D+D'_1)$.
\end{remark}
\begin{definition}\label{defn-subdiv}
 We write $\Sub_r(D)$ for the scheme of divisors $D'$ of degree $r$
 such that $D'\leq D$.  Using Remark~\ref{rem-included} we see that
 $\Sub_r(D)=\Sub_r(r,D)=\Int_r(r,D)$.  
\end{definition}

\begin{theorem}\label{thm-int-dec}
 The ring $\O_{\Int_r(d_0,d_1)}$ is freely generated over 
 \[ \OS\psb{c_{0i}\st 0<i\leq d_0-r}\psb{c_{1j}\st 0<j\leq d_1} \]
 by the monomials 
 \[ c_0^\al:=\prod_{i=d_0-r+1}^{d_0}c_{0i}^{\al_i} \]
 for which $\sum_i\al_i\leq d_1-r$.  Moreover, if we let 
 $\pi\:\Sub_r(d_0,d_1)\xra{}\Int_r(d_0,d_1)$ be the usual projection,
 then the corresponding ring map $\pi^*$ is a split monomorphism of
 modules over $\O_{\Div_{d_1}^+}$ (so $\pi$ itself is dominant).
\end{theorem}

The proof will be given after a number of intermediate results.  It
seems likely that the injectivity of $\pi^*$ could be extracted from
work of Pragacz~\cite[Section 3]{pr:egd}.  He works with Chow groups
of varieties rather than generalised cohomology rings of spaces, and
his methods and language are rather different; we have not attempted a
detailed comparison.

We start by setting up some streamlined notation.  We put $n=d_0-r$
and $m=d_1-r$.  We use the following names for the coordinate rings of
various schemes of divisors, and the standard generators of these
rings:
\begin{align*}
 C_0 &= \O_{\Div_{d_0}^+} = \OS\psb{u_1,\ldots,u_{n+r}} \\
 C_1 &= \O_{\Div_{d_1}^+} = \OS\psb{v_1,\ldots,v_{m+r}} \\
 A   &= \O_{\Div_n^+}     = \OS\psb{a_1,\ldots,a_n}   \\
 B   &= \O_{\Div_m^+}     = \OS\psb{b_1,\ldots,b_m}   \\
 C   &= \O_{\Div_r^+}     = \OS\psb{c_1,\ldots,c_r}.
\end{align*}
(In particular, we have renamed $c_{0i}$ and $c_{1i}$ as $u_i$ and
$v_i$.)  We put $u_0=v_0=a_0=b_0=c_0=1$.  We define $u_i=0$ for $i<0$
or $i>n+r$, and similarly for $v_i$, $a_i$, $b_i$ and $c_i$.  The
equations of the various tautological divisors are as follows:
\begin{align*}
 f_0(x) &= \sum_i u_i x^{n+r-i} \in C_0[x] \\
 f_1(x) &= \sum_i v_i x^{m+r-i} \in C_1[x] \\
 f(x)   &= \sum_i a_i x^{n-i}   \in A[x]   \\
 g(x)   &= \sum_i b_i x^{m-i}   \in B[x]   \\
 h(x)   &= \sum_i c_i x^{r-i}   \in C[x].
\end{align*}
We write $T_0$ for the set of monomials of weight at most $m$ in
$u_{n+1},\ldots,u_{n+r}$, and $T$ for the set of monomials 
of weight at most $m$ in $c_1,\ldots,c_r$.  We also
introduce the subrings
\begin{align*}
 C'_0  &= \OS\psb{u_1,\ldots,u_n} \sse C_0 \\
 C''_0 &= \OS\psb{u_1,\ldots,u_{n-1}} \sse C'_0.
\end{align*}
We note that the ring $Q:=\O_{\Int_r(d_0,d_1)}$ has the form
$(C_0\hot C_1)/I$ for a certain ideal $I$.  The theorem claims that
$Q$ is freely generated as a module over $C'_0\hot C_1$ by $T_0$.

The map 
\[ \pi^*\:C_0\hot C_1\xra{}A\hot B\hot C \]
sends $f_0(x)$ to $f(x)h(x)$ and $f_1(x)$ to $g(x)h(x)$.  This induces
a map $\pi^*\:Q\xra{}A\hot B\hot C$, and the theorem also claims that
this is a split injection.

We will need to approximate certain determinants by calculating their
lowest terms with respect to a certain ordering.  More precisely, we
consider monomials of the form $u^\al=\prod_{i=1}^{n+r}u_i^{\al_i}$,
and we order these by $u^\al<u^\bt$ if there exists $i$ such that
$\al_i>\bt_i$ and $\al_j=\bt_j$ for $j>i$.  The mnemonic is that
$u_1\ll\ldots\ll u_{n+r}$, so any difference in the exponent of $u_i$
overwhelms any difference in the exponents of $u_1,\ldots,u_{i-1}$.
\begin{lemma}\label{lem-leading}
 Suppose we have integers $\gm_i$ satisfying
 $0\leq\gm_0<\ldots<\gm_m<m+r$, and we put $M_{ij}=u_{n+r+i-\gm_j}$
 for $0\leq i,j\leq m$, where $u_k$ is interpreted as $0$ if $k<0$ or
 $k>n+r$.  Then the lowest term in $\det(M)$ is the product of the
 diagonal entries, so
 \[ \det(M)=\prod_{i=0}^m u_{n+r+i-\gm_i}+\text{ higher terms }. \]
\end{lemma}
\begin{remark}
 Determinants of this type are known as \emph{Schur functions}.
\end{remark}
\begin{proof}
 Put $\dl=\prod_{i=0}^m u_{n+r+i-\gm_i}$.  Let $M'_i$ be obtained from
 $M$ by removing the $0$'th row and $i$'th column.  The matrix $M'_0$
 has the same general form as $M$ so by induction we have
 $\det(M'_0)=\prod_{i=1}^m u_{n+r+i-\gm_i}+\text{ higher terms }$.  If
 we expand $\det(M)$ along the top row then the $0$'th term is
 $u_{n+r-\gm_0}\det(M'_0)=\dl+\text{ higher terms }$.  As
 $0\leq\gm_0<\ldots<\gm_m$ we have $\gm_i\geq i+\gm_0$ and so $\dl$
 only involves variables $u_j$ with $j\leq n+r-\gm_0$.  The remaining
 terms in the row expansion of $\det(M)$ have the form
 $(-1)^iu_{n+r-\gm_0+i}\det(M'_i)$ for $i>0$, and $u_{n+r-\gm_0+i}$ is
 either zero (if $i>\gm_0$) or a variable strictly higher than all
 those appearing in $\dl$.  The lemma follows easily.
\end{proof}

\begin{lemma}\label{lem-gens}
 The ring $Q$ is generated by $T_0$ as a module over $C'_0\hot C_1$.
\end{lemma}
\begin{proof}
 Let $J$ be the ideal in $C'_0\hot C_1$ generated by $u_1,\ldots,u_n$
 and $v_1,\ldots,v_{m+r}$, so $(C'_0\hot C_1)/J=\OS$.  We also put
 $C''_0=(C_0\hot C_1)/J=\OS\psb{u_{n+1},\ldots,u_{n+r}}$.  As $J$
 is topologically nilpotent, it will suffice to prove the result
 modulo $J$.  We will thus work modulo $J$ throughout the proof, so
 that $f_1=x^{m+r}$, and we must show that $Q/J$ is generated over
 $\OS$ by $T_0$.  

 Let $\mu\:C''_0\psb{x}/x^{m+r}\xra{}C''_0\psb{x}/x^{m+r}$ be defined
 by $\mu(t)=f_0t$, and let $M$ be the matrix of $\mu$ with respect to
 the obvious bases.  It is then easy to see that $Q/J=C''_0/I$, where
 $I$ is generated by the minors of $M$ of size $m+1$.  The entries in
 $M$ are $M_{ij}=u_{n+r+i-j}$.

 We next claim that all the generators $u_k$ are nilpotent mod $I$, or
 equivalently that $u_k=0$ in the ring $R=C''_0/\sqrt{I}$ for all $k$.   
 By downward induction we may assume that $u_l=0$ in $R$ for
 $k<l\leq n+r$.  We consider the submatrix $M'$ of $M$ given by
 $M'_{ij}=M_{i,n+r-k+j}=u_{i+k-j}$ for $0\leq i,j\leq m$.  By the
 definition of $I$ we have $\det(M')\in I$ and thus $\det(M')=0$ in
 $R$.  On the other hand, we have $u_l=0$ for $l>k$ so $M'$ is lower
 triangular so $\det(M')=\prod_iM'_{ii}=u_k^{m+1}$.  Thus $u_k$ is
 nilpotent in $R$ but clearly $\text{Nil}(R)=0$ so $u_k=0$ in $R$ as
 required.  It follows that $Q/J$ is a quotient of the polynomial ring
 $\OS[u_{n+1},\ldots,u_{n+r}]\subset\OS\psb{u_{n+1},\ldots,u_{n+r}}$.  

 Now let $W$ be the submodule of $Q/J$ spanned over $\OS$ by $T_0$; we
 must prove that this is all of $Q/J$.  As $1\in W$, it will suffice
 to show that $W$ is an ideal.  In the light of the previous
 paragraph, it will suffice to show that $W$ is closed under
 multiplication by the elements $u_{n+1},\ldots,u_{n+r}$, or
 equivalently that $W$ contains all monomials of weight $m+1$.

 We thus let $\al=(\al_{n+1},\ldots,\al_{n+r})$ be a multiindex of
 weight $m+1$.  There is then a unique sequence $(\bt_0,\ldots,\bt_m)$
 with $n+r\geq\bt_0\geq\ldots\geq\bt_m>n$ and
 $u^\al=\prod_iu_{\bt_i}$.  Put $\gm_i=n+r+i-\bt_i$, so that
 $0\leq\gm_0<\ldots<\gm_m<m+r$.  Let $M_\al$ be the submatrix of $M$
 consisting of the first $m+1$ columns of the rows of indices
 $\gm_0,\ldots,\gm_m$, so the $(i,j)$'th entry of $M_\al$ is
 $u_{n+r+i-\gm_j}$.  Note that the elements $r_\al:=\det(M_\al)$ lie
 in $I$.

 Lemma~\ref{lem-leading} tells us that the lowest term in $r_\al$ is
 $\prod_iu_{n+r+i-\gm_i}=\prod_iu_{\bt_i}=u^\al$.  It is clear that
 the weight of the remaining terms is at most the size of $M_\al$,
 which is $m+1$.  By an evident induction, we may assume that their
 images in $C''_0/I$ lie in $W$.  As $r_\al\in I$ we deduce that
 $u^\al\in W$ as well.
\end{proof}

\begin{corollary}\label{cor-gens}
 Let $D_1$ be a divisor of degree $d_1$ on $\GG$ over $S'$, for some
 scheme $S'$ over $S$.  Then $\O_{\Int_r(d_0,D_1)}$ is generated over 
 $\O_{S'}\psb{c_{01},\ldots,c_{0,d_0-r}}$ by the monomials
 $c_0^\al=\prod_{i=d_0-r+1}^{d_0}c_{0i}^{\al_i}$ for which
 $|\al|\leq d_1-r$.
\end{corollary}
\begin{proof}
 The previous lemma is the universal case.
\end{proof}

We next treat the special case of Theorem~\ref{thm-int-dec} where
$n=0$ and so $r=d_0$.  As remarked in Definition~\ref{defn-subdiv},
the map $\pi\:\Sub_r(d_1)=\Sub_r(r,d_1)\xra{}\Int_r(r,d_1)$ is an
isomorphism in this case.

\begin{lemma}\label{lem-points}
 Let $D$ be a divisor of degree $d$ on $\GG$ over $S$.  For any
 $r\leq d$ we let $P_r(D)$ denote the scheme of tuples
 $(u_1,\ldots,u_r)\in\GG^r$ such that $\sum_{i=1}^r[u_r]\leq D$.  Then
 $\O_{P_r(D)}$ is free of rank $d!/(d-r)!$ over $\O_S$.
\end{lemma}
\begin{proof}
 There is an evident projection $P_r(D)\xra{}P_{r-1}(D)$, which
 identifies $P_r(D)$ with the divisor $D-[u_1]-\ldots-[u_{r-1}]$ on
 $\GG$ over $P_{r-1}(D)$.  This divisor has degree $d-r+1$, so
 $\O_{P_r(D)}$ is free of rank $d-r+1$ over $\O_{P_{r-1}(D)}$.  It
 follows by an evident induction that $\O_{P_r(D)}$ is free over
 $\OS$, with rank $d(d-1)\ldots(d-r+1)=d!/(d-r)!$.
\end{proof}

\begin{lemma}\label{lem-int-dec}
 Let $D$ be a divisor of degree $d$ on $\GG$ over $S$, let $D'$ be
 the tautological divisor of degree $r$ over $\Sub_r(D)$, and let
 $f(x)=\sum_{i=0}^rc_ix^{r-i}$ be the equation of $D$.  Then the set
 $T$ of monomials of degree at most $d-r$ in $c_1,\ldots,c_r$ is a basis
 for $\O_{\Sub_r(D)}$ over $\OS$.
\end{lemma}
\begin{proof}
 Put $K=|T|$; by elementary combinatorics we find that $K=\bcf{d}{r}$.
 Put $R=\O_{\Sub_r(D)}$.  Using $T$ we obtain an $\OS$-linear map
 $\bt\:\O_S^K\xra{}R$, which is surjective by Lemma~\ref{lem-gens}; we
 must prove that it is actually an isomorphism.

 Now consider the scheme $P_r(D)$; Lemma~\ref{lem-points} tells us
 that the ring $R':=\O_{P_r(D)}$ is a free module over $\OS$ of rank
 $d!/(d-r)!=r!K$.  On the other hand, $P_r(D)$ can be identified with
 the scheme of tuples $(D',u_1,\ldots,u_r)$ where $D'\in\Sub_r(D)$ and
 $D'=[u_1]+\ldots+[u_r]$.  In other words, if we change base to
 $\Sub_r(D)$ we can regard $P_r(D)$ as $P_r(D')$, and now
 Lemma~\ref{lem-points} tells us that $R'$ is free of rank $r!$ over
 $R$.

 Now choose a basis $e_1,\ldots,e_{r!}$ for $R'$ over $R$.  We can
 combine this with $\bt$ to get a map $\gm\:\OS^{r!K}\xra{}R'$.  This
 is a direct sum of copies of $\bt$, so it is surjective.  Both source
 and target of $\gm$ are free of rank $r!K$ over $\OS$.  Any
 epimorphism between free modules of the same finite rank is an
 isomorphism (choose a splitting and then take determinants).  Thus
 $\gm$ is an isomorphism, and it follows that $\bt$ is an isomorphism
 as required. 
\end{proof}

\begin{corollary}
 The set $T$ is a basis for $B\hot C$ over $C_1$.
\end{corollary}
\begin{proof}
 This is the universal case of the lemma.
\end{proof}

\begin{corollary}\label{cor-T-basis}
 The set $T$ is a basis for $A\hot B\hot C$ over $C'_0\hot C_1$.
\end{corollary}
\begin{proof}
 Note that $A\hot B\hot C=(B\hot C)\psb{a_1,\ldots,a_n}$.  For
 $0<i\leq n$ we have
 \[ \pi^*u_i=\sum_{i=j+k}a_jc_k=a_i+c_i
     \text{ mod decomposables, }
 \]
 where $c_i$ may be zero, but $a_i$ is nonzero.  It follows that our
 ring $A\hot B\hot C$ can also be described as
 $(B\hot C)\psb{u_1,\ldots,u_n}$, or equivalently as
 $C'_0\hot B\hot C$.  The claim now follows easily from the previous
 corollary.
\end{proof}

Now let $T_1$ be the set of monomials of the form
$u_n^ic_1^{\al_1}\ldots c_n^{\al_n}$ for which $0\leq i<|\al|\leq m$.
These monomials can be regarded as elements of $A\hot B\hot C$, giving
a map $(C''_0\hot C_1)\{T_1\}\xra{}A\hot B\hot C$.  The map
$\pi^*\:C_0\hot C_1\xra{}A\hot B\hot C$ also gives us a map 
$(C'_0\hot C_1)\{T_0\}\xra{}A\hot B\hot C$, and by combining these we
get a map
\[ \phi\:(C'_0\hot C_1)\{T_0\}\op (C''_0\hot C_1)\{T_1\}
   \xra{} A\hot B\hot C \simeq (C'_0\hot C_1)\{T\}
\]
of modules over $C''_0\hot C_1$.  Our main task will be to prove that
this is an isomorphism.  The proof will use the following lemma.

\begin{lemma}
 Let $R$ be a ring, and let $\al\:M\xra{}N$ be a homomorphism of
 modules over $R\psb{x}$.  Suppose that $M$ can be written as a
 product of copies of $R\psb{x}$, and similarly for $N$.  Suppose also
 that the induced map $M/xM\xra{}N/xN$ is an isomorphism.  Then $\al$
 is also an isomorphism.
\end{lemma}
\begin{proof}
 We have diagrams as shown below, in which the rows are easily seen to
 be exact:
 \begin{diag}
  \node{M/x^kM}     \arrow{e,t,V}{x} \arrow{s,l}{\al} 
  \node{M/x^{k+1}M} \arrow{e,A}      \arrow{s,l}{\al} 
  \node{M/xM}                        \arrow{s,r}{\al} \\
  \node{N/x^kN}     \arrow{e,b,V}{x} 
  \node{N/x^{k+1}N} \arrow{e,A}      
  \node{N/xN}                        
 \end{diag}%
We see by induction that the maps $M/x^kM\xra{}N/x^kN$ are all
 isomorphisms, and the claim follows by taking inverse limits.
\end{proof}

Our map $\phi$ is a map of modules over the ring 
\[ C''_0\hot C_1=\OS\psb{u_1,\ldots,u_{n-1},v_1,\ldots,v_{m+r}}. \]
Moreover, we have
$C'_0\hot C_1=(C''_0\hot C_1)\psb{u_n}\simeq
 \prod_{k=0}^\infty C''_0\hot C_1$.  Now let $J$ be the ideal in
$C''_0\hot C_1$ generated by
$\{u_1,\ldots,u_{n-1},v_1,\ldots,v_{m+r}\}$, so 
$(C''_0\hot C_1)/J=\OS$ and $\phi$ induces a map
\[ \phb\:\OS\psb{u_n}\{T_0\}\op\OS\{T_1\} \xra{}
    \OS\psb{u_n}\{T\}.
\]
Note also that $\OS\{T\}$ is the image of $C$ in $(A\hot B\hot C)/J$
and is thus a subring of $\OS\psb{u_n}\{T\}$.  By an evident inductive
extension of the lemma, it will suffice to show that $\phb$ is an
isomorphism.

\begin{lemma}\label{lem-ucw}
 We have $u_{n+j}=u_nc_j+w_j\pmod{J}$ for some polynomial $w_j$ in
 $c_1,\ldots,c_r$.
\end{lemma}
\begin{proof}
 For any monic polynomial $p(x)$ of degree $d$ we write
 $\hat{p}(y)=y^dp(1/y)$.  If $p(x)=\sum_ir_ix^{d-i}$ then
 $\hat{p}(y)=\sum_ir_iy^i$.  Note that $\widehat{pq}=\hat{p}\hat{q}$,
 and that $\hat{p}(0)=1$.  As we work mod $(u_i\st i<n)$ we have
 $\hat{f}_0=1\pmod{y^n}$.  As we work mod $(v_j\st j\leq m+r)$ we have
 $\hat{f}_1=1$.  We also have $fh=f_0$ and $gh=f_1$, so
 $\hat{f}\hat{h}=\hat{f}_0=1\pmod{y^n}$ and
 $\hat{g}\hat{h}=\hat{f}_1=1$.  It follows easily that
 $\hat{f}=\hat{g}\pmod{y^n}$, so $a_i=b_i$ for $i<n$.

 We now have to distinguish between the case $m<n$ and the case
 $m\geq n$.  First suppose that $m<n$.  Then for $i>n$ we have
 $a_i=b_i=0$, and also $b_n=0$, and $a_i=b_i$ for $i<n$ by the
 previous paragraph.  This implies that $\hat{f}-\hat{g}=a_ny^n$.  We
 also have $(\hat{f}-\hat{g})\hat{h}=\hat{f}_0-1$, and by comparing
 coefficients we deduce that $a_nc_i=u_{n+i}$ for $i=0,\ldots,r$.  The
 case $i=0$ gives $u_n=a_n$, so $u_{n+i}=u_nc_i$ for $i=1,\ldots,r$,
 so the lemma is true with $w_i=0$.

 Now suppose instead that $m\geq n$.  As $a_i=0$ for $i>n$ we have
 \[ \hat{f}-\hat{g}-(a_n-b_n)y^n=-\sum_{i=n+1}^m b_iy^i\in C[y]. \]
 We now multiply this by $\hat{h}$ and use the fact that
 $(\hat{f}-\hat{g})\hat{h}=\hat{f}_0-1$.  By comparing coefficients of
 $y^n$ we find that $u_n=a_n-b_n$.  In view of this, our equation
 reads 
 \[ \hat{f}_0-1-u_ny^n\hat{h}=
     -(\sum_{i=n+1}^m b_iy^i)\hat{h}\in C[y].
 \]
 The right hand side has the form $\sum_{j>0}w_jy^{n+j}$ with
 $w_j\in C$, and by comparing coefficients we see that
 $u_{n+j}=u_nc_j+w_j$ as claimed.
\end{proof}

\begin{proof}[Proof of Theorem~\ref{thm-int-dec}]
 Lemma~\ref{lem-ucw} tells us that $\phb(u^\al)$ is $u_n^{|\al|}c^\al$
 plus terms involving lower powers of $u_n$.  It follows easily that
 if we filter the source and target of $\phb$ by powers of $u_n$, then
 the resulting map of associated graded modules is a isomorphism.  It
 follows that $\phb$ is an isomorphism, and thus that $\phi$ is an
 isomorphism.  It follows that the map
 $(C'_0\hot C_1)\{T_0\}\xra{}A\hot B\hot C$ is a split monomorphism of
 modules over $C''_0\hot C_1$ (and thus certainly of modules over
 $C_1$).  We have seen that this map factors as
 \[ (C'_0\hot C_1)\{T_0\} \xra{\psi} Q \xra{\pi^*} A\hot B\hot C, \]
 where $\psi$ is surjective by Lemma~\ref{lem-gens}.  It follows that
 $\psi$ is an isomorphism and that $\pi^*$ is a split monomorphism, as
 required. 
\end{proof}

\section{Flag spaces}
\label{sec-flags}

In the next section we will (in good cases) construct spaces whose
associated formal schemes are the schemes $\Sub_r(D(V_0),D(V_1))$ and
$\Int_r(D(V_0),D(V_1))$ considered previously.  As a warm-up, and also
as technical input, we will first consider the schemes associated to
Grassmannian bundles and flag bundles.  The results discussed are
essentially due to Grothendieck~\cite{gr:qpf}; we have merely adjusted
the language and technical framework.

Let $V$ be a bundle of dimension $d$ over a space $X$.  We write
$P_r(V)$ for the space of tuples $(x,L_1,\ldots,L_r)$ where $x\in X$
and $L_1,\ldots,L_k\in PV_x$ and $L_i$ is orthogonal to $L_j$ for
$i\neq j$.  Recall also that in Lemma~\ref{lem-points} we defined
$P_r(D(V))$ to be the scheme over $X_E$ of tuples
$(u_1,\ldots,u_r)\in\GG^r$ for which $[u_1]+\ldots+[u_r]\leq D(V)$.

\begin{proposition}\label{prop-Pr}
 There is a natural isomorphism $P_r(V)_E=P_r(D(V))$.
\end{proposition}
\begin{proof}
 For each $i$ we have a line bundle over $P_r(V)$ whose fibre over the
 point $(x,L_1,\ldots,L_r)$ is $L_i$.  This is classified by a map
 $P_r(V)\xra{}\CPi$, which gives rise to a map
 $u_i\:P_r(V)_E\xra{}\GG$.  The direct sum of these line bundles
 corresponds to the divisor $[u_1]+\ldots+[u_r]$.  This direct sum is
 a subbundle of $V$, so $[u_1]+\ldots+[u_r]\leq D(V)$.  This
 construction therefore gives us a map $P_r(V)_E\xra{}P_r(D(V))$.

 In the case $r=1$ we have $P_1(V)=PV$ and $P_1(D(V))=D(V)$ so the
 claim is that $(PV)_E=D(V)$, which is true by definition.  In
 general, suppose we know that $P_{r-1}(V)_E=P_{r-1}(D(V))$.  We can
 regard $P_r(V)$ as the projective space of the bundle over
 $P_{r-1}(V)$ whose fibre over a point $(x,L_1,\ldots,L_{r-1})$ is the
 space $V_x\om(L_1\op\ldots\op L_{r-1})$.  It follows that $P_r(V)_E$
 is just the divisor $D(V)-([u_1]+\ldots+[u_{r-1}])$ over
 $P_{r-1}(D(V))$, which is easily identified with $P_r(D(V))$.  The
 proposition follows by induction.
\end{proof}
\begin{remark}\label{rem-Pr}
 One can easily recover the following more concrete statement.  The
 ring $E^0P_r(V)=\O_{P_r(D(V))}$ is the largest quotient ring of
 $(E^0X)\psb{x_1,\ldots,x_r}$ in which the polynomial $f_V(t)$ is
 divisible by $\prod_{i=1}^k(t-x_i)$.  It is a free module over $E^0X$
 with rank $d!/(d-r)!$, and the monomials $x^\al$ with
 $0\leq\al_i\leq d-i$ (for $i=1,\ldots,r$) form a basis.  More details
 about the multiplicative relations are given in
 Section~\ref{sec-PkD}.  
\end{remark}

We next consider the Grassmannian bundle 
\[ G_r(V)=\{(x,W)\st x\in X\;,\; W\leq V_x \text{ and }\dim(W)=r\}.
\]
\begin{proposition}\label{prop-Gr}
 There is a natural isomorphism $G_r(V)_E=\Sub_r(D(V))$.
\end{proposition}
\begin{proof}
 Let $T$ denote the tautological bundle over $G_r(V)$.  This is a rank
 $r$ subbundle of the pullback of $V$ so we have a degree $r$
 subdivisor $D(T)$ of the pullback of $D(V)$ over $G_r(V)_E$.  This
 gives rise to a map $G_r(V)\xra{}\Sub_r(D(V))$.  

 Next, consider the space $P_r(V)$.  There is a map
 $P_r(V)\xra{}G_r(V)$ given by
 $(x,L_1,\ldots,L_r)\mapsto(x,L_1\op\ldots\op L_r)$.  This lifts in an
 evident way to give a homeomorphism $P_r(V)\simeq P_r(T)$.  Of
 course, this is exactly parallel to the proof of
 Lemma~\ref{lem-int-dec}.  Over $P_r(D(T))$ we have points
 $a_1,\ldots,a_r$ of $\GG$ with coordinate values
 $x_1,\ldots,x_r\in\O_{P_r(D(T))}$ say.  Let $B$ be the set of monomials
 $x^\al$ with $0\leq\al_i\leq r-i$ for $i=1,\ldots,r$.  From our
 earlier analysis of $\Sub_r(D(V))$ and $P_r(D(T))$ we see that $B$ is
 a basis for $\O_{P_r(D(T))}$ over $\O_{\Sub_r(D(V))}$.  We also see
 from Remark~\ref{rem-Pr} (applied to the bundle $T$) that $B$ is a
 basis for $E^0P_r(T)$ over $E^0G_r(V)$.  This means that our
 isomorphism $f\:\O_{P_r(D(V))}\xra{}E^0P_r(V)$ is a direct sum
 (indexed by $B$) of copies of our map
 $g\:\O_{\Sub_r(D(V))}\xra{}E^0G_r(V)$.  It follows that $g$ must also
 be an isomorphism.
\end{proof}
\begin{remark}
 Lemma~\ref{lem-int-dec} now gives us an explicit basis for
 $E^0G_r(V)$ over $E^0X$, consisting of monomials in the Chern classes
 of the tautological bundle $T$.
\end{remark}

\section{Topological universal examples}
\label{sec-top-universal}

In this section we construct spaces whose associated formal schemes
are the algebraic universal examples considered in
Section~\ref{sec-alg-universal}. 

We first consider the easy case of the schemes $\Sub_r(D_0,D_1)$.
\begin{definition}\label{defn-Gr}
 Given vector bundles $V_0$ and $V_1$ over $X$, we define
 $G_r(V_0,V_1)$ to be the space of quadruples $(x,W_0,W_1,g)$ such
 that
 \begin{itemize}
  \item[\rm(a)] $x\in X$;
  \item[\rm(b)] $W_i$ is an $r$-dimensional subspace of $V_{ix}$ for
   $i=0,1$; and
  \item[\rm(c)] $g$ is an isometric isomorphism $W_0\xra{}W_1$. 
 \end{itemize}
 (We would obtain a homotopy equivalent space if we dropped the
 requirement that $g$ be an isometry.)

 If $V_i$ is the evident tautological bundle over $BU(d_i)$ we write
 $G_r(d_0,d_1)$ for $G_r(V_0,V_1)$.  More generally, if $V$ is a
 bundle over $X$ and $d_0\geq 0$ we can let $V_1$ be the pullback of
 $V$ to $BU(d_0)\tm X$, and let $V_0$ be the pullback of of the
 tautological bundle over $BU(d_0)$; in this context we write
 $G_r(d_0,V)$ for $G_r(V_0,V_1)$.
\end{definition}

\begin{theorem}\label{thm-Gr}
 There is a natural map
 $p\:G_r(V_0,V_1)_E\xra{}\Sub_r(D(V_0),D(V_1))$.  In the universal
 case this is an isomorphism, so
 \[ G_r(d_0,d_1)_E = \Sub_r(d_0,d_1). \]
 More generally, there is a spectral sequence
 \[ \Tor_{E^*BU(d_0)\tm BU(d_1)}^{**}(E^*X,E^*G_r(d_0,d_1))
     \convto E^*G_r(V_0,V_1),
 \]
 whose edge map in degree zero is the map
 \[ p^* \: \O_{\Sub_r(D(V_0),D(V_1))} \xra{} E^0G_r(V_0,V_1). \]
 The spectral sequence collapses in the universal case.  (We do not
 address the question of convergence in the general case.)
\end{theorem}
\begin{proof}
 First, we can pull back the bundles $V_i$ from $X$ to $G_r(V_0,V_1)$
 (without change of notation).  We also have a bundle over
 $G_r(V_0,V_1)$ whose fibre over a point $(x,W_0,W_1,g)$ is $W_0$; we
 denote this bundle by $W$, and note that there are natural inclusions
 $V_0\xla{1}W\xra{g}V_1$.  We then have divisors $D(W)$ and $D(V_i)$
 on $\GG$ over $G_r(V_0,V_1)_E$ with $D(W)\leq D(V_0)$ and
 $D(W)\leq D(V_1)$, so the triple $(D(V_0),D(V_1),D(W))$ is classified
 by a map $G_r(V_0,V_1)_E\xra{}\Sub_r(D(V_0),D(V_1))$.

 We next consider the universal case.  As our model of $EU(d)$ we use
 the space of orthonormal $d$-frames in $\Ci$, so $BU(d)$ is
 just the Grassmannian of $d$-planes in $\Ci$.  Given a point
 \[ (\un{u},\un{v})=(u_1,\ldots,u_{d_0},v_1,\ldots,v_{d_1})\in
     EU(d_0)\tm EU(d_1)
 \]
 we construct a point $((V_0,V_1),W_0,W_1,g)\in G_r(d_0,d_1)$ as
 follows:
 \begin{itemize}
  \item[\rm(a)] $V_0$ is the span of $u_1,\ldots,u_{d_0}$ 
  \item[\rm(b)] $V_1$ is the span of $v_1,\ldots,v_{d_1}$ 
  \item[\rm(c)] $W_0$ is the span of $u_1,\ldots,u_r$ 
  \item[\rm(d)] $W_1$ is the span of $v_1,\ldots,v_r$
  \item[\rm(e)] $g$ is the map $W_0\xra{}W_1$ that sends $u_i$ to $v_i$. 
 \end{itemize}
 This gives a map $f\:EU(d_0)\tm EU(d_1)\xra{}G_r(d_0,d_1)$.  Next,
 the group $U(d_0)\tm U(d_1)$ has a subgroup
 $U(r)\tm U(d_0-r)\tm U(r)\tm U(d_1-r)$, inside which we have the
 smaller subgroup $\Gm$ consisting of elements of the form
 $(h,k_0,h,k_1)$.  It is not hard to see that
 $\Gm\simeq U(r)\tm U(d_0-r)\tm U(d_1-r)$, and that $f$ gives a
 homeomorphism $(EU(d_0)\tm EU(d_1))/\Gm\xra{}G_r(d_0,d_1)$.
 Moreover, $EU(d_0)\tm EU(d_1)$ is contractible and $\Gm$ acts freely
 so $G_r(d_0,d_1)\simeq B\Gm=BU(r)\tm BU(d_0-r)\tm BU(d_1-r)$, so
 $G_r(d_0,d_1)_E=\Div^+_r\tm\Div^+_{d_0-r}\tm\Div^+_{d_1-r}=
  \Sub_r(d_0,d_1)$ as claimed.

 In the general case we can choose maps $f_i\:X\xra{}BU(d_i)$
 classifying $V_i$, and this gives rise to a pullback square as
 follows: 
 \begin{diag}
  \node{G_r(V_0,V_1)} \arrow{e}\arrow{s}
  \node{G_r(d_0,d_1)} \arrow{s} \\
  \node{X} \arrow{e} \node{BU(d_0)\tm BU(d_1).}
 \end{diag}%
The vertical maps are fibre bundle projections so this is actually a
 homotopy pullback square.  This give an Eilenberg-Moore spectral
 sequence as in the statement of the theorem.  On the edge we have
 \[ E^*X \ot_{E^*BU(d_0)\tm BU(d_1)} E^*G_r(d_0,d_1), \]
 which is the same as 
 \[ E^*X \ot_{E^0BU(d_0)\tm BU(d_1)} E^*G_r(d_0,d_1). \]
 We can now identify this as the tensor product of $E^*X$ with
 $\O_{\Sub_r(d_0,d_1)}$ over $\O_{\Div_{d_0}^+\tm\Div_{d_1}^+}$.  The
 part in degree zero is easily seen to be $\O_{\Sub_r(D(V_0),D(V_1))}$
 as claimed.
\end{proof}

We next show that our map $G_r(V_0,V_1)_E\xra{}\Sub_r(D(V_0),D(V_1))$
is an isomorphism in the semiuniversal case as well as the universal
case.  We start by analysing the semiuniversal spaces $G_r(d_0,V)$ in
more familiar terms.
\begin{proposition}\label{prop-Gr-semi}
 There are natural homotopy equivalences
 \[ G_r(d_0,V)\simeq G_r(V)\tm BU(d_0-r) \]
 (and in particular $G_r(r,V)\simeq G_r(V)$).
\end{proposition}
\begin{proof}
 A point of $G_r(d_0,V)$ is a tuple $(V_0,x,W_0,W_1,g)$ where
 $V_0\in G_{d_0}(\Ci)$, $x\in X$, $W_0\in G_r(V_0)$,
 $W_1\in G_r(V_x)$ and $g\:W_0\xra{}W_1$.  We can define a map
 $f\:G_r(d_0,V)\xra{}G_r(V)\tm BU(d_0-r)$ by
 $f(V_0,x,W_0,W_1,g)=(x,W_1,V_0\om W_0)$.  It is not hard to see
 that this is a fibre bundle projection, and that the fibre over a
 point $(x,W,V')$ is the space of linear isometric embeddings from $W$
 to $\Ci\om V'$.  This space is homeomorphic to the space of
 linear isometric embeddings of $\C^r$ in $\Ci$, which is
 well-known to be contractible.  Thus $f$ is a fibration with
 contractible fibres and thus is a weak equivalence.
\end{proof}

\begin{corollary}
 The map $G_r(d_0,V)_E\xra{}\Sub_r(d_0,D(V))$ is an isomorphism.
\end{corollary}
\begin{proof}
 Recall that $\Sub_r(d_0,D(V))$ is the scheme of pairs $(D_1,D)$ where
 $D_1$ is a divisor of degree $d_1$, $D$ is a divisor of degree $r$
 and $D\sse D_1\cap D(V)$.  There is an evident isomorphism
 $\Sub_r(D(V))\tm_S\Div_{d_0-r}^+\xra{}\Sub_r(d_0,D(V))$ sending
 $(D',D)$ to $(D+D',D)$.  The proposition tells us that
 $G_r(d_0,V)=G_r(V)\tm BU(d_0-r)$.  We already know that
 $BU(d_0-r)_E=\Div_{d_0-r}^+$, and Proposition~\ref{prop-Gr} tells us
 that $G_r(V)_E=\Sub_r(D(V))$.  We therefore have an isomorphism
 $G_r(d_0,V)_E=\Sub_r(D(V))\tm_S\Div_{d_0-r}^+=\Sub_r(d_0,D(V))$.
 (This involves an implicit K\"unneth isomorphism, which is valid
 because $BU(d_0-r)$ has only even-dimensional cells.)  We leave it to
 the reader to check that this isomorphism is the same as the map
 considered previously.
\end{proof}

We now turn to parallel results for the schemes
$\Int_r(D(V_0),D(V_1))$.  
\begin{definition}
 Given vector bundles $V_0$ and $V_1$ over a space $X$, we define
 $I_r(V_0,V_1)$ to be the space of pairs $(x,f)$ where
 $f\:V_{0x}\xra{}V_{1x}$ is a linear map of rank at least $r$.  We
 define the universal and semiuniversal spaces $I_r(d_0,d_1)$ and
 $I_r(d_0,V)$ by the evident analogue of Definition~\ref{defn-Gr}. 
\end{definition}
\begin{remark}
 There is a natural map 
 \[ G_r(V_0,V_1)\xra{}I_r(V_0,V_1), \]
 sending $(x,W_0,W_1,g)$ to $(x,f)$, where $f$ is the composite
 \[ V_0 \xra{\text{proj}} W_0 \xra{g} W_1 \xra{\text{inc}} V_1. \]
 This gives a homeomorphism of $G_r(V_0,V_1)$ with the subspace of
 $I_r(V_0,V_1)$ consisting of pairs $(x,f)$ for which $f^*f$ and
 $ff^*$ are idempotent.
\end{remark}

\begin{definition}\label{defn-I-Int}
 We define a natural map $q\:I_r(V_0,V_1)\xra{}\Int_r(D(V_0),D(V_1))$
 as follows.  If we let $\pi$ denote the projection
 $I_r(V_0,V_1)\xra{}X$ then we have a tautological map
 $f\:\pi^*V_0\xra{}\pi^*V_1$ which has rank at least $r$ everywhere.
 Proposition~\ref{prop-int-equiv} now tells us that
 $\int(\pi^*V_0,\pi^*V_1)\geq r$.  We can therefore apply
 Theorem~\ref{thm-int} and deduce that the map
 $I_r(V_0,V_1)_E\xra{}X_E$ factors through a map
 $q\:I_r(V_0,V_1)\xra{}\Int_r(D(V_0),D(V_1))\sse X_E$ as required.
\end{definition}

Later we will show that the map $q$ is an isomorphism in the universal
case.  For this, it will be convenient to have an alternative model
for the universal space $I_r(d_0,d_1)$.
\begin{proposition}
 Put 
 \[ I'_r(d_0,d_1)=
    \{(V_0,V_1)\in G_{d_0}(\Ci)\tm G_{d_1}(\Ci) 
               \st \dim(V_0\cap V_1) \geq k
    \}.
 \]
 Then $I'_r(d_0,d_1)$ is homotopy equivalent to $I_r(d_0,d_1)$.
\end{proposition}
\begin{proof}
 The basic idea is to refine the proof of
 Proposition~\ref{prop-int-equiv}.  We will take $G_d(\Ci)$ as
 our model for $BU(d)$.  We write $I=I_r(d_0,d_1)$ and
 $I'=I'_r(d_0,d_1)$ for brevity.

 We will need various isometries between infinite-dimensional vector
 spaces.  We define $\dl\:\Ci\xra{}\Ci\op\Ci$ by
 $\dl(v)=(v,v)/\sqrt{2}$, and we define $\tht\:\Ci\op\Ci\xra{}\Ci$ by
 $\tht(v,w)=(v_0,w_0,v_1,w_1,\ldots)$.  Next, it is well-known that
 the space of linear isometric embeddings of $\Ci$ in itself is
 contractible, so we can choose a continuous family of isometries
 $\phi_t$ with $\phi_0=\tht\dl$ and $\phi_1=1$.  Similarly, we can
 choose continuous families of isometric embeddings
 $\psi_0^t,\psi_1^t\:\Ci\xra{}\Ci\op\Ci$ with $\psi_0^0(v)=\tht(v,0)$ and
 $\psi_1^0(v)=\tht(0,v)$ and $\psi_0^1(v)=\psi_1^1(v)=v$.

 We now define a map $\al\:I'\xra{}I$ by $\al(V_0,V_1)=(V_0,V_1,f)$,
 where $f$ is the orthogonal projection map from $V_0$ to $V_1$.  This
 acts as the identity on $V_0\cap V_1$ and thus has rank at least $k$.
 If we choose $n$ large enough that $V_0+V_1\leq\C^n$ and let
 $V_0\xra{i_0}\C^n\xla{i_1}V_1$ be the inclusions, then $f=i_1^*i_0$.

 Next, we need to define a map $\bt\:I\xra{}I'$.   Given
 $(V_0,V_1,f)\in I$ we can construct maps
 \begin{align*}
  \mu \: V_0 & \xra{} V_0 \\
  \nu \: V_1 & \xra{} V_1 \\
  j_0 \: V_0 & \xra{} V_0\op V_1 < \Ci\op\Ci \\
  j_1 \: V_1 & \xra{} V_0\op V_1 < \Ci\op\Ci
 \end{align*}
 as in the proof of the implication (c)$\Rightarrow$(b) in
 Proposition~\ref{prop-int-equiv}, so
 $\dim(j_0V_0\cap j_1V_1)\geq k$.  We can thus define $\bt\:I\xra{}I'$
 by $\bt(V_0,V_1,f)=(\tht j_0V_0,\tht j_1V_1)$.

 Suppose we start with $(V_0,V_1)\in I'$, define $f\:V_0\xra{}V_1$ to
 be the orthogonal projection, and then define $j_0,j_1$ as above so
 that $\bt\al(V_0,V_1)=(\tht j_0V_0,\tht j_1V_1)$.  Observe that
 $f^*f\:V_0\xra{}V_0$ decreases distances, and acts as the identity on
 $V:=V_0\cap V_1$.  If we let $\lm_1,\ldots,\lm_{d_0}$ be the
 eigenvalues of $f^*f$ (listed in the usual way) we deduce that
 $\lm_1=\ldots=\lm_k=1$ and that $0\leq\lm_i\leq 1$ for all $i$.  It
 follows from this that $\mu$ and $\nu$ are the respective identity
 maps, so
 \begin{align*}
  j_0 &= (1,f)\circ(1+f^*f)^{-1/2} \\
  j_1 &= (f^*,1)\circ(1+ff^*)^{-1/2}.
 \end{align*}
 In particular, we have $j_0(v)=j_1(v)=(v,v)/\sqrt{2}$ for $v\in V$,
 so $j_0|_V=j_1|_V=\dl|_V$.

 Next, for $0\leq t\leq 1$ we define $j^t_0\:V_0\xra{}\Ci\op\Ci$ by
 \[ j^t_0 = (i_0,ti_0 + (1-t)f)\circ(1+t^2+(1-t^2)f^*f)^{-1/2}. \]
 One can check that this is an isometric embedding, with $j^0_0=j_0$
 and $j^1_0=\dl|_{V_0}$ and $j^t_0|_V=\dl|_V$ for all $t$.  Similarly,
 if we put
 \[ j^t_1=(i_1,ti_0 + (1-t)f^*)\circ(1+t^2+(1-t^2)ff^*)^{-1/2}, \]
 we find that this is an isometric embedding of $V_1$ in $\Ci\op\Ci$
 with $j^0_1=j_1$ and $j^1_1=\dl|_{V_1}$ and $j^t_1|_V=\dl|_V$ for all
 $t$.  It follows that $(\tht j^t_0V_0,\tht j^t_1V_1)\in I'$ for all
 $t$, and this gives a path from
 $\bt\al(V_0,V_1)=(\tht j_0V_0,\tht j_1V_1)$ to
 $(\tht\dl V_0,\tht\dl V_1)$.  Recall that we chose a path
 $\{\phi_t\}$ from $\tht\dl$ to $1$.  The pairs
 $(\phi_tV_0,\phi_tV_1)$ now give a path from
 $(\tht\dl V_0,\tht\dl V_1)$ to $(V_0,V_1)$ in $I'$.  Both of the
 paths considered above are easily seen to depend continuously on the
 point $(V_0,V_1)\in I'$ that we started with, so we have constructed
 a homotopy $\bt\al\simeq 1$.

 Now suppose instead that we start with a point $(V_0,V_1,f)\in I$; we
 need a path from $\al\bt(V_0,V_1,f)$ to $(V_0,V_1,f)$.  We have
 $\bt(V_0,V_1,f)=(\tht j_0V_0,\tht j_1V_1)$, so
 $\al\bt(V_0,V_1,f)=(\tht j_0V_0,\tht j_1V_1,f')$, where
 $f'\:\tht j_0V_0\xra{}\tht j_1V_1$ is the orthogonal projection.  One
 can check that this is characterised by
 $f'(\tht j_0(v))=\tht j_1(j_1^*j_0(v))$.  Next, for $0\leq t\leq 1$
 we define $k^t_0\:V_0\xra{}V_0\op V_1$ by 
 \[ k^t_0=
     (\sqrt{1-t^2+t^2\mu},tf)\circ(1-t^2+t^2\mu+t^2f^*f)^{-1/2}.
 \]
 This is an isometric embedding with $k^1_0=j_0$ and
 $k^0_0(v)=(v,0)$.  Similarly, we define $k^t_1\:V_1\xra{}V_0\op V_1$
 by  
 \[ k^t_1=
     (tf^*,\sqrt{1-t^2+t^2\nu})\circ(1-t^2+t^2\nu+t^2ff^*)^{-1/2},
 \]
 and we define $f'_t\:\tht k^t_0V_0\xra{}\tht k^t_1V_1$ by
 \[ f'_t(\tht k^t_0(v))=\tht k^t_1(j_1^*j_0(v)), \]
 so $f'_1=f'$.  The points $(k^t_0V_0,k^t_1V_1,f'_t)$ give a path from
 $\al\bt(V_0,V_1,f)$ to $(\tht(V_0\op 0),\tht(0\op V_1),f'_0)$ in $I$.

 Next, we define $f''_t\:\psi_0^tV_0\xra{}\psi_1^tV_1$ by
 $f''_t(\psi^t_0(v))=\psi^t_1(j_1^*j_0(v))$.  The points
 $(\psi^t_0V_0,\psi^t_1V_1,f''_t)$ give a path from 
 $(\tht(V_0\op 0),\tht(0\op V_1),f'_0)$ to $(V_0,V_1,j_1^*j_0)$ in
 $I$.  Using Proposition~\ref{prop-square} one can check that
 \begin{align*}
  j_1^*j_0 &=
   (\nu+ff^*)\circ(f\sqrt{\mu}+\sqrt{\nu}f)\circ(\mu+f^*f) \\ 
   &= f\circ(2\mu^{1/2}(\mu+f^*f)^{-1}).
 \end{align*}
 The map $\zt:=2\mu^{1/2}(\mu+f^*f)^{-1}$ is a strictly positive
 self-adjoint automorphism of $V_0$, so the same is true of $t+(1-t)\zt$
 for $0\leq t\leq 1$.  The points $(V_0,V_1,f\circ(t+(1-t)\zt))$ form
 a path from $(V_0,V_1,j_1^*j_0)$ to $(V_0,V_1,f)$.  All the paths
 considered depend continuously on the point $(V_0,V_1,f)$ that we
 started with, so we have defined a homotopy $\al\bt\simeq 1$. 
\end{proof}

\begin{theorem}
 The map $q\:I_r(d_0,d_1)_E\xra{}\Int_r(d_0,d_1)$ is an isomorphism.
\end{theorem}
\begin{proof}
 We first replace $I_r(d_0,d_1)$ by the homotopy-equivalent space
 $I'_r(d_0,d_1)$.  We write $I_r=I'_r(d_0,d_1)$ and $G_r=G_r(d_0,d_1)$
 for brevity, and similarly for $\Int_r$ and $\Sub_r$.  We first claim
 that there is a commutative diagram as follows.
 \begin{diag}
  \node{\O_{\Int_r}}
   \arrow{s,l}{q}
   \arrow{e,V}
  \node{\O_{\Sub_r}}
   \arrow{s,lr}{p}{\simeq} \\
  \node{E^0I_r}
   \arrow{e}
  \node{E^0G_r.}
 \end{diag}%
Indeed, the isomorphism
 \[ p\:\O_{\Sub_r(d_0,d_1)}\xra{}E^0G_r(d_0,d_1) \]
 comes from Theorem~\ref{thm-Gr}, and the map $q$ comes from
 Definition~\ref{defn-I-Int}.  It was proved in
 Theorem~\ref{thm-int-dec} that the top horizontal map is a split
 monomorphism of $\OS$-modules, and it follows that the same is true
 of the map $q\:\O_{\Int_r}\xra{}E^0I_r$.

 We now specialise to the case where $E$ is $H[u,u^{-1}]$, the
 two-periodic version of the integer Eilenberg-MacLane spectrum.  We
 then have $E^0X=\prod_kH^{2k}X$ for all spaces $X$.  This splits each
 of the rings on the bottom row of our diagram as a product of
 homogeneous pieces, and it is not hard to check that there is a
 unique compatible way to split the rings on the top row.  We know
 that $q$ is a split monomorphism; if we can show that the source and
 target have the same Poincar\'e series, it will follow that $q$ is an
 isomorphism.  If $r=\min(d_0,d_1)$ then $\Int_r=\Sub_r$ so the claim
 is certainly true.  To work downwards from here by induction, it will
 suffice to show that
 \[ PS(H^*I_{r+1})-PS(H^*I_r) = PS(\O_{\Int_{r+1}})-PS(\O_{\Int_r}) \]
 for all $r$.

 To evaluate the left hand side, we consider the space 
 \[ I_r\sm I_{r+1} =
     \{(V_0,V_1)\in G_{d_0}(\Ci)\tm G_{d_1}(\Ci) \st 
        \dim(V_0\cap V_1)=k \}.
 \]
 Let $G'_r$ be the space of triples $(V,V'_0,V'_1)$ of mutually
 orthogonal subspaces of $\Ci$ such that $\dim(V)=r$ and
 $\dim(V_i)=d_i-r$.  This is well-known to be a model of
 $BU(r)\tm BU(d_0-r)\tm BU(d_1-r)$ and thus homotopy-equivalent to
 $G_r$; the argument uses frames much as in the proof of
 Theorem~\ref{thm-Gr}.  Let $W$ be the bundle over $G'_r$ whose fibre
 over $(V,V'_0,V'_1)$ is $\Hom(V'_1,V'_0)$.  If
 $\al\in\Hom(V'_0,V'_1)$ and we put $V_0=V\op V'_1$ and
 $V_1=V\op\text{graph}(\al)$ then $V_0\cap V_1=V$ and so
 $(V_0,V_1)\in I_r$.  It is not hard to see that this construction
 gives a homeomorphism of the total space of $W$ with
 $I_r\sm I_{r+1}$.  This in turn gives a homeomorphism of the Thom
 space $(G'_r)^W$ with the quotient space $I_r/I_{r+1}$.  By induction
 we may assume that $H^*I_{r+1}$ is concentrated in even degrees, and
 it is clear from the Thom isomorphism theorem that the same is true of
 $\tH^*(G'_r)^W$.  This implies that $H^*I_r$ is in even degrees and
 that $PS(H^*I_r)-PS(H^*I_{r+1})=PS(\tH^*(G'_r)^W)$.  As $W$ has
 dimension $(d_0-r)(d_1-r)$, we see that
 $PS(\tH^*(G'_r)^W)=t^{2(d_0-r)(d_1-r)}PS(H^*G'_r)$.  We also know
 that $H^*G'_r\simeq\O_{\Sub_r}$.  The conclusion is that 
 \[ PS(H^*I_r)-PS(H^*I_{r+1})=t^{2(d_0-r)(d_1-r)}PS(\O_{\Sub_r}). \]

 We next evaluate $PS(\O_{\Int_{r+1}})-PS(\O_{\Int_r})$.  Put
 \[ R_r^*=\Z\psb{c_{01},\ldots,c_{0,d_0-r},c_{11},\ldots,c_{1,d_1}}.
 \] 
 We know from Theorem~\ref{thm-int-dec} that $\O_{\Int_r}$ is freely
 generated over $R_r^*$ by the monomials
 $\prod_{i=1}^rc_{0,d_0-r+i}^{\al_i}$ for which
 $\sum_{i=1}^r\al_i\leq d_1-r$.  It follows that the monomials
 $\prod_{i=0}^rc_{0,d_0-r+i}^{\al_i}$ for which
 $\sum_1^r\al_i\leq d_1-r$ form a basis for $\O_{\Int_r}$ over
 $R_{r+1}^*$.  Similarly, those for which $\sum_{i=0}^r\al_i<d_1-r$
 form a basis for $\O_{\Int_{r+1}}$ over $R_{r+1}^*$.  Thus, if we let
 $M^*$ be the module generated over $R_{r+1}^*$ by the monomials with
 $\sum_1^r\al_i\leq d_1-r\leq\sum_0^r\al_i$, we find that
 $PS(\O_{\Int_{r+1}})-PS(\O_{\Int_r})=PS(M^*)$.

 It is not hard to check that the monomials for which
 $\sum_0^r\al_i=d_1-r$ form a basis for $M^*$ over $R^*_r$.  Next, let
 $N^*$ be generated over $\Z$ by the monomials
 $\prod_{i=0}^rc_i^{\al_i}$ for which $\sum_0^r\al_i=d_1-r$; note that
 this involves the variables $1=c_0,\ldots,c_r$ rather than the
 variables $c_{d_0-r},\ldots,c_{d_0}$ used in $M^*$.  Because
 $\deg(c_{d_0-r+i})=\deg(c_i)+2(d_0-r)$ we have
 \[ \deg(\prod_ic_{d_0-r+i}^{\al_i})=
    \deg(\prod_ic_i^{\al_i})+2(d_0-r)\sum_i\al_i.
 \]
 Using this, we see that
 $PS(M^*)=t^{2(d_0-r)(d_1-r)}PS(N^*)PS(R_r^*)$.  However,
 Corollary~\ref{cor-T-basis} essentially says that
 $\O_{\Sub_r}\simeq R_r^*\ot N^*$ as graded Abelian groups, so
 $PS(N^*)PS(R_r^*)=PS(\O_{\Sub_r})$, so
 \begin{align*}
  PS(\O_{\Int_{r+1}})-PS(\O_{\Int_r}) &=
    t^{2(d_0-r)(d_1-r)}PS(\O_{\Sub_r}) \\
     & = PS(H^*I_r) - PS(H^*I_{r+1}).
 \end{align*}
 
 As explained previously, this implies that $q$ is an isomorphism in
 the case $E=H[u,u^{-1}]$.  We next consider the case
 $E=MU[u,u^{-1}]$.  Let $I$ be the kernel of the usual map
 $MU_*\xra{}\Z$.  Because $H^*I_r$ is free of finite type and
 concentrated in even degrees, we see that the Atiyah-Hirzebruch
 spectral sequence collapses and that the associated graded ring
 $\text{gr}_IMU^*I_r$ is isomorphic to
 $\text{gr}_I(MU^*)\hot H^*I_r$.  Using this it is not hard to check
 that $q$ is an isomorphism in the case $E=MU[u,u^{-1}]$ also.
 Finally, given an arbitrary even periodic ring spectrum $E$ we can
 choose a complex orientation in $\tE^0\CPi$ and thus a ring map
 $MU[u,u^{-1}]\xra{}E$.  Using this, we deduce that $q$ is an
 isomorphism for all $E$.
\end{proof}

\begin{corollary}\label{cor-Ir}
 Let $V_0$ and $V_1$ be bundles of dimensions $d_0$ and $d_1$ over a
 space $X$.  Then there is a spectral sequence
 \[ \Tor_{E^*BU(d_0)\tm BU(d_1)}^{**}(E^*X,E^*I_r(d_0,d_1))
     \convto E^*I_r(V_0,V_1),
 \]
 whose edge map in degree zero is the map
 \[ q^* \: \O_{\Int_r(D(V_0),D(V_1))} \xra{} E^0I_r(V_0,V_1). \]
 The spectral sequence collapses in the semiuniversal and universal
 cases.  (We do not address the question of convergence in the general
 case.)
\end{corollary}
\begin{proof}
 This is another Eilenberg-Moore spectral sequence.
\end{proof}

\section{The schemes $P_kD$}
\label{sec-PkD}

Let $D$ be a divisor of degree $d$ on $\GG$ over $S$, with equation
\[ f(t)=f_D(t)=\sum_{i=0}^dc_ix^{d-i}\in\OS[t], \]
say.  In this section we assemble some useful facts about the scheme
$P_kD$.  This is a closed subscheme of $\GG^k$, so
$\O_{P_kD}=\OS\psb{x_0,\ldots,x_{k-1}}/J_k$ for some ideal $J_k$; our
main task will be to find systems of generators for $J_k$.  We put
$p_i(t)=\prod_{j<i}(t-x_j)$, and we let $q_i(t)$ and $r_i(t)$ be the
quotient and remainder when $f(t)$ is divided by $p_i(t)$.  Thus
$f(t)=q_i(t)p_i(t)+r_i(t)$ and $r_i(t)$ has the form
$\sum_{j=0}^{i-1}a_{ij}t^j$ for some $a_{i0},\ldots,a_{i,i-1}\in\OS$.
From the definitions it is clear that $J_k$ is the smallest ideal
modulo which $f(t)$ becomes divisible by $p_k(t)$, or in other words
the smallest ideal modulo which $r_k(t)=0$, so $J_k$ is generated by
$a_{k0},\ldots,a_{k,k-1}$.  Now put $b_i=a_{i+1,i}$ for $0\leq i<k$;
we will show that these elements also generate $J_k$.

\begin{lemma}
 We have $b_i=q_i(x_i)$ and $r_{i+1}=b_ip_i+r_i$ for all $i$.
\end{lemma}
\begin{proof}
 The polynomial $q_i(t)-q_i(x_i)$ is evidently divisible by $t-x_i$,
 say $q_i(t)-q_i(x_i)=(t-x_i)q'_{i+1}(t)$.  If we put
 $r'_{i+1}(t)=q_i(x_i)p_i(t)+r_i(t)$ we find that $r'_{i+1}$ is a
 polynomial of degree at most $i$ and that
 $f(t)=q'_{i+1}(t)p_{i+1}(t)+r'_{i+1}(t)$, so we must have
 $q_{i+1}=q'_{i+1}$ and $r_{i+1}=r'_{i+1}$.  Thus $b_i$ is the
 coefficient of $t^i$ in $r'_{i+1}(t)$.  As $r_i$ has degree less than
 $i$ and $p_i$ is monic of degree $i$ we deduce that $b_i=q_i(x_i)$.
\end{proof}

\begin{corollary}
 The ideal $J_k$ is generated by $b_0,\ldots,b_{k-1}$.
\end{corollary}
\begin{proof}
 Put $J'_k=(b_0,\ldots,b_{k-1})$.  If we work modulo $J'_k$ then it is
 immediate from the lemma that $r_k=r_{k-1}=\ldots=r_0=0$; this shows
 that $J_k\sse J'_k$.  Conversely, if we work modulo $J_k$ then $f$ is
 divisible by $p_k$ and hence by $p_i$ for all $i\leq k$, so
 $r_0=\ldots=r_k=0$.  It follows from the lemma that $b_ip_i=0$ for
 all $i$, and $p_i$ is monic so $b_i=0$.  Thus $J'_k\sse J_k$.
\end{proof}

We now give a determinantal formula for the relators $b_j$.  Consider
the Vandermonde determinant
\[ v_k:=\det(x_i^j)_{0\leq i,j<k}=\prod_{0\leq i<j<k}(x_j-x_i). \]
We also define a matrix $B_k$ by 
\[ (B_k)_{ij} = 
    \begin{cases} x_i^j & \text{ if } 0\leq j<k-1 \\
                  f(x_i)& \text{ if } j=k-1.
    \end{cases}
\]

\begin{proposition}\label{prop-b-det}
 We have $b_j=\det(B_j)/v_j$ for all $j$.  (More precisely, we have
 $v_jb_j=\det(B_j)\in\OS\psb{x_0,\ldots,x_{j-1}}$, and $v_j$ is not
 a zero-divisor in this ring.)
\end{proposition}
\begin{proof}
 Define $\al\:\OS^j\xra{}\OS\psb{t}/p_j(t)$ by
 \[ \al(u_0,\ldots,u_{j-2},w)=\sum_iu_it^i+wf(t) \pmod{p_j}
                             =\sum_iu_it^i+wr_j(t) \pmod{p_j}.
 \]
 Next, define $\bt\:\OS\psb{t}/p_j(t)\xra{}\OS^j$ by
 $\bt(g)=(g(x_0),\ldots,g(x_{j-1}))$.  We identify $\OS\psb{t}/p_j(t)$
 with $\OS^j$ using the basis $\{t^i\st 0\leq i<j\}$.  It is easy to
 see that $\det(\bt)=v_j$ and $\det(\bt\al)=\det(B_j)$.  Moreover, the
 matrix of $\al$ has the form
 $\left(\begin{array}{c|c} I & * \\ \hline 0 & b_j\end{array}\right)$
 so $\det(\al)=b_j$.  It follows immediately that $v_jb_j=\det(B_j)$.
 It is easy to see that none of the polynomials $x_j-x_i$ (where
 $i<j$) are zero-divisors, so $v_j$ is not a zero-divisor.
\end{proof}

We next need to relate the schemes $P_kD$ to the exterior powers
$\lm^k\OD$. 
\begin{lemma}\label{lem-alt-Jk}
 The ideal $J_k$ maps to zero under the natural projection
 $\O_{D^k}=\O_D^{\ot k}\xra{}\lm^k\OD$.
\end{lemma}
\begin{proof}
 It is enough to prove the corresponding result in the universal case,
 where $D$ is the tautological divisor over $\Div_d^+$.  As the map
 $\GG^d\xra{}\GG^d/\Sg_d=\Div_d^+$ is faithfully flat, it is enough to
 prove the result after pulling back along this map.  In other words,
 we need only consider the divisor over the ring
 $R:=\O_{\GG^d}=\OS\psb{y_i\st i<d}$ with equation
 $f(t)=\prod_i(t-y_i)$.  Let $w$ be the discriminant of this
 polynomial, so $w=\prod_{i\neq j}(y_i-y_j)\in R$.  Put
 $N=\{0,\ldots,d-1\}$, and let $F(N,R)$ denote the ring of functions
 from $N$ to $R$, with pointwise operations.  We can define
 $\phi\:\OD\xra{}F(N,R)$ by $\phi(g)(i)=g(y_i)$, and the Chinese
 Remainder Theorem tells us that the resulting map
 $w^{-1}\OD\xra{}F(N,w^{-1}R)$ is an isomorphism, and it follows that
 $w^{-1}\O_{D^k}=F(N^k,w^{-1}R)$.  We also have
 $\O_{D^k}=R\psb{x_j\st j<k}/(f_D(x_j)\st j<k)$; the element $x_j$
 corresponds to the function $\um\mapsto y_{n_j}$.

 Now put
 \[ N_k = \{ (n_0,\ldots,n_{k-1})\in N^k
              \st n_i\neq n_j \text{ when } i\neq j \},
 \]
 and $N^c_k=N^k\setminus N_k$.  Let $r(t)$ be the remainder when the
 polynomial $f_D(t):=\prod_{i<d}(t-y_i)$ is divided by
 $f_{D'}(t):=\prod_{j<k}(t-x_j)$.  This corresponds to the function
 $\um\mapsto r_\um(t)$, where $r_\um(t)$ is the remainder of $f_D(t)$
 modulo $\prod_{j<k}(t-y_{n_j})$.  As the discriminant is invertible
 in $w^{-1}R$ we see that $r_\um(t)=0$ iff $\um\in N_k$, and otherwise
 some coefficient of $r_\um(t)$ is invertible.  Using this, we deduce
 that $w^{-1}\O_{P_kD}=F(N_k,w^{-1}R)$ and
 $w^{-1}J_k=F(N^c_k,w^{-1}R)$.  If we let $\{e_0,\ldots,e_{k-1}\}$ be
 the evident basis of $F(N,R)$ over $R$, this means that $w^{-1}J_k$
 is spanned over $w^{-1}R$ by the elements
 $e_{n_0}\ot\ldots\ot e_{n_{k-1}}$ for which $n_i=n_j$ for some
 $i\neq j$, and these elements satisfy
 $e_{n_0}\Smash\ldots\Smash e_{n_{k-1}}=0$ so the map
 $w^{-1}J_k\xra{}w^{-1}\lm^k\OD$ is zero.  As $w$ is not a
 zero-divisor we deduce that the map $J_k\xra{}\lm^k\OD$ is zero, as
 claimed.
\end{proof}

Next note that the symmetric group $\Sg_k$ acts on $D^k$ and $P_kD$
and thus on the corresponding rings.  In either case we define
$\alt_k(a)=\sum_{\sg\in\Sg_k}\sgn(\sg)\sg.a$.  We also let
$\mu_k\:\O_D^{\ot k}\xra{}\lm^k\OD$ be the usual projection, or
equivalently the restriction of the product map
$\mu_k\:(\lm^*\OD)^{\ot k}\xra{}\lm^*\OD$.  Dually, we let
$\psi_k\:\lm^k\OD\xra{}\O_D^{\ot k}$ be the component of the coproduct
map $\psi_k\:\lm^*\OD\xra{}(\lm^*\OD)^{\ot k}$.  We also let
$p^*\:\O_D^{\ot k}=\O_{D^k}\xra{}\O_{P_kD}$ denote the usual
projection, corresponding to the closed inclusion $P_kD\xra{}D^k$.

\begin{proposition}\label{prop-cross}
 There is a natural commutative diagram as follows.
 \begin{diag}
  \node{\O_{P_kD}}
  \node[2]{\O_{P_kD}}
   \arrow[2]{w,t}{\alt_k} 
   \arrow{sw,t,A}{\mu'_k} \\  
  \node[2]{\lm^k\OD} 
   \arrow{nw,t}{\psi'_k} 
   \arrow{sw,t,V}{\psi_k} \\
  \node{\O_D^{\ot k}} 
   \arrow[2]{n,l,A}{p^*}
  \node[2]{\O_D^{\ot k}} 
   \arrow[2]{n,r,A}{p^*}
   \arrow{nw,t,A}{\mu_k}
   \arrow[2]{w,b}{\alt_k}
 \end{diag}
\end{proposition}
\begin{proof}
 The main point to check is that
 $\psi_k\mu_k=\alt_k\:\O_D^{\ot k}\xra{}\O_D^{\ot k}$.  Consider an
 element $a=a_0\ot\ldots\ot a_{k-1}\in\O_D^{\ot k}$.  Let $a_i^j$
 denote the element $1^{\ot j}\ot a_i\ot 1^{\ot k-j-1}\in\O_D^{\ot k}$,
 so that $\psi_k(a_i)=\sum_{j=0}^{k-1}a_i^j$ and
 $\psi_k\mu_k(a)=\prod_i\sum_ja_i^j$.  We are interested in the
 component of this in $\O_D^{\ot k}\subset(\lm^*\OD)^{\ot k}$, which is
 easily seen to be $\sum_\sg\prod_ia_i^{\sg(i)}$.  Moreover, one
 checks that 
 \[ \prod_ia_i^{\sg(i)} =
     \sgn(\sg)a_{\sg^{-1}(0)}\ot\ldots\ot a_{\sg^{-1}(k-1)} = 
     \sgn(\sg) \sg.a,
 \]
 so the relevant component of $\psi_k\mu_k(a)$ is
 $\sum_\sg\sgn(\sg)\sg.a=\alt_k(a)$, as claimed.

 Let $A$ be the set of multiindices $\al=(\al_0,\ldots,\al_{k-1})$
 with $0\leq\al_i<d$ for all $i$, and let $A_0$ be the subset of those
 for which $\al_0>\ldots>\al_{k-1}$.  Put
 $x^\al:=x^{\al_0}\ot\ldots\ot x^{\al_{k-1}}\in\O_D^{\ot k}$.  Then
 $\{x^\al\st \al\in A\}$ is a basis for $\O_D^{\ot k}$, and
 $\{\mu_k(x^\al)\st \al\in A_0\}$ is a basis for $\lm^k\OD$.
 Moreover, if $\al\in A_0$ and we write
 $\psi_k\mu_k(x^\al)=\alt_k(x^\al)=\sum_{\bt\in A}c_{\al\bt}x^\bt$ we
 see that $c_{\al\al}=1$ and $c_{\al\bt}=0$ if $\bt\in A_0$ and
 $\bt\neq\al$.  It follows that $\mu_k$ is surjective and $\psi_k$ is
 a split injection of $\OS$-modules, as indicated in the diagram.  
 
 Lemma~\ref{lem-alt-Jk} tells us that $\mu_k$ factor as $\mu'_k p^*$
 for some $\mu'_k\:\O_{P_kD}\xra{}\lm^k\OD$, and a diagram chase shows
 that $\mu'_k$ is surjective.  This gives the right hand triangle of
 the diagram.  We simply define $\psi'_k=p^*\psi_k$ to get the left
 hand triangle.  As $p^*$ is $\Sg_k$-equivariant we have
 \[ \alt_kp^* = p^*\alt_k = p^*\psi_k\mu_k = \psi'_k \mu'_k p^*. \]
 As $p^*$ is surjective, this proves that $\psi'_k\mu'_k=\alt_k$, so
 the top triangle commutes. 
\end{proof}

We next study certain orbit schemes for actions of $\Sg_k$.  Recall
that $\O_{\GG^k}=\OS\psb{x_i\st i<k}$ has a topological basis
consisting of monomials in the variables $x_i$.  This basis is
permuted by $\Sg_k$, and the sums of the orbits form a topological
basis for the invariant subring
$\O_{\GG^k}^{\Sg_k}=\O_{\GG^k/\Sg_k}=\O_{\Div_k^+}$.  It is clear from
this analysis that our quotient construction commutes with base
change, in other words $(S'\tm_S\GG^k)/\Sg_k=S'\tm_S(\GG^k/\Sg_k)$ for
any scheme $S'$ over $S$.  Similarly, the set 
$\{x^\al\st\al_i<d\text{ for all }i\}$ is a basis for $\O_{D^k}$ that
is permuted by $\Sg_k$, so the orbit sums give a basis for
$\O_{D^k}^{\Sg_k}$ and we have a quotient scheme
$D^k/\Sg_k=\spf(\O_{D^k}^{\Sg_k})$ whose formation commutes with base
change.  By comparing our bases we see that the projection
$\OG\xra{}\OD=\OG/f_D$ induces a surjective map
$\O_{\GG^k/\Sg_k}\xra{}\O_{D^k/\Sg_k}$.  In other words, we have a
commutative square of schemes as shown, in which $j$ and $j'$ are
closed inclusions, and $q_2$ is a faithfully flat map of degree $d!$.
\begin{diag}
 \node{D^k}       \arrow{e,t,V}{j}  \arrow{s,l}{q_1}
 \node{\GG^k}                       \arrow{s,r,A}{q_2} \\
 \node{D^k/\Sg_k} \arrow{e,b,V}{j'}
 \node{\GG^k/\Sg_k=\Div_k^+} 
\end{diag}%
One might hope to show that $P_k(D)/\Sg_k=\Sub_k(D)$ in a similar
sense, but this is not quite correct.  For example if $D=3[0]$ (so
$f_D(t)=t^3$) and $k=2$ then
$\O_{P_2D}=\OS\psb{x,y}/(x^3,x^2+xy+y^2)$.  If we define a basis of
this ring by
\[ \{e_0,\ldots,e_5\} = \{1,x,y,x^2,-x^2-xy,x^2y\}, \]
we find that the generator of $\Sg_2$ has the effect
\[ e_0 \leftrightarrow e_0 \;,\;
   e_1 \leftrightarrow e_2 \;,\;
   e_3 \leftrightarrow e_4 \;,\;
   e_5 \leftrightarrow -e_5.
\]
If $\OS$ has no $2$-torsion we find that $\O_{P_2D}^{\Sg_2}$ is
spanned by $\{1,x+y,xy\}$ and thus is equal to $\O_{\Sub_2(D)}$.
However, if $2=0$ in $\OS$ we have an additional generator $x^2y$, so
$\O_{P_2D}^{\Sg_2}$ is strictly larger than $\O_{\Sub_2(D)}$.  This
example also shows that the formation of $\O_{P_2D}^{\Sg_2}$ is not
compatible with base change.

The following proposition provides a substitute for the hope described
above. 
\begin{proposition}\label{prop-invariant}
 There is a commutative diagram as follows, in which $i$, $i'$, $j$
 and $j'$ are closed inclusions, and $q_0$ and $q_2$ are faithfully
 flat of degree $k!$.  Moreover, the outer rectangle is a pullback,
 and if $J_k:=\ker(i^*)$ then
 $\ker((i')^*)=J_k^{\Sg_k}$.  
 \begin{diag}
  \node{P_kD}      \arrow{e,t,V}{i}  \arrow{s,l,A}{q_0}
  \node{D^k}       \arrow{e,t,V}{j}  \arrow{s,l}{q_1}
  \node{\GG^k}                       \arrow{s,r,A}{q_2} \\
  \node{\Sub_k(D)} \arrow{e,b,V}{i'}
  \node{D^k/\Sg_k} \arrow{e,b,V}{j'}
  \node{\GG^k/\Sg_k=\Div_k^+} 
 \end{diag}
\end{proposition}
\begin{proof}
 We have already produced the right hand square.  The map $i$ is just
 the obvious inclusion.  The map $q_0$ sends
 $(a_0,\ldots,a_{k-1})\in P_kD$ to
 $[a_0]+\ldots+[a_{k-1}]\in\Sub_k(D)$; it was observed in the proof of
 Lemma~\ref{lem-int-dec} that this makes $\O_{P_kD}$ into a free
 module of rank $k!$ over $\O_{\Sub_k(D)}$, so $q_0$ is faithfully
 flat of degree $k!$.

 The points of $\Sub_k(D)$ are the divisors of degree $k$ contained in
 $D$, so $\Sub_k(D)$ is a closed subscheme of $\Div_k^+$; we write
 $m'\:\Sub_k(D)\xra{}\Div_k^+$ for the inclusion, and note that
 $m'q_0=q_2ji$.  As $q_0$ is faithfully flat and $m'q_0$ factors
 through $D^k/\Sg_k$ we see that $m'$ factors through $D^k/\Sg_k$, so
 there is a unique map $i'\:\Sub_k(D)\xra{}D^k/\Sg_k$ such that
 $m'=j'i'$.  As $m'$ is a closed inclusion, the same is true of $i'$.
 A point of the pullback of $m$ and $q_2$ is a list
 $a=(a_0,\ldots,a_{k-1})$ of points of $\GG$ such that the divisor
 $q_2(a)=\sum_r[a_r]$ lies in $\Sub_k(D)$, and thus satisfies
 $\sum_r[a_r]\leq D$.  It follows from the definitions that this
 pullback is just $P_kD$ as claimed.

 As $q_0$ is faithfully flat we have
 $\ker((i')^*)=\ker(q_0^*(i')^*)=\ker(i^*q_1^*)$.  By construction,
 $q_1^*$ is just the inclusion of the $\Sg_k$-invariants in
 $\O_{D^k}$, so $\ker(i^*q_1^*)=\ker(i^*)^{\Sg_k}=J_k^{\Sg_k}$ as
 claimed. 
\end{proof}
\begin{corollary}\label{cor-module}
 $\lm^k\OD$ is naturally a module over $\O_{\Sub_k(D)}$.
\end{corollary}
\begin{proof}
 We can certainly regard $\O_{D^k}$ as a module over the subring
 $\O_{D^k/\Sg_k}=\O_{D^k}^{\Sg_k}$, and the map
 $\alt_k\:\O_{D^k}\xra{}\O_{D^k}$ respects this structure.  This makes
 $\lm^k\OD=\img(\alt_k)$ into a module over $\O_{D^k}^{\Sg_k}$.  If
 $a\in J_k^{\Sg_k}$ and $b\in\O_{D^k}$ then $a\alt_k(b)=\alt_k(ab)$
 but $ab\in J_k$ so $\alt_k(ab)=0$.  This shows that $\lm^k\O_D$ is
 annihilated by $J_k^{\Sg_k}$, so it is a module over
 $\O_{D^k}^{\Sg_k}/J_k^{\Sg_k}=\O_{\Sub_k(D)}$, as claimed.
\end{proof}

We next identify $\lm^k\OD$ as a module over $\O_{\Sub_k(D)}$.  Let
$D'$ be the tautological divisor of degree $k$ over $\Sub_k(D)$.  Then
$\O_{D'}$ is naturally a quotient of the ring
$\OD\ot_S\O_{\Sub_k(D)}$, which contains the subring $\OD=\OD\ot 1$.
This gives us a map $\OD\xra{}\O_{D'}$, which extends to give a map 
$\phi\:\lm^k_S\OD\xra{}\lm^k_{\Sub_k(D)}\O_{D'}$.

\begin{proposition}\label{prop-phi-iso}
 The map $\phi\:\lm^k_S\OD\xra{}\lm^k_{\Sub_k(D)}\O_{D'}$ is an
 isomorphism of free rank-one modules over $\O_{\Sub_k(D)}$. 
\end{proposition}
\begin{proof}
 Put $T=\Sub_k(D)$ for brevity.  Note that
 \begin{align*}
  \OD            &= \OS\{x^i\st i<d\} \\
  \lm^k_S\OD     &= \OS\{ x^{i_0}\Smash\ldots\Smash x^{i_{k-1}} \st
                          d>i_0>\ldots>i_{k-1} \} \\
  \O_{D'}        &= \OT\{x^i \st i<k\} \\
  \lm^k_T\O_{D'} &= \OT\{x^{k-1}\Smash\ldots\Smash x^0\}.
 \end{align*}
 In particular, we see that $\lm^k_S\OD$ is free of rank
 $K:=\bcf{d}{k}$ over $\OS$, and $\lm^k_T\O_{D'}$ is free of rank one
 over $\OT$.  We also know from Lemma~\ref{lem-int-dec} that $\OT$ is
 free of rank $K$ over $\OS$, so $\lm^k_T\O_{D'}$ is also free of rank
 $K$ over $\OS$.  

 Suppose for the moment that $\phi$ is a homomorphism of
 $\OT$-modules.  It is clear that 
 \[ \phi(x^{k-1}\Smash\ldots\Smash x^0)=x^{k-1}\Smash\ldots\Smash x^0,
 \]
 and this element generates $\lm^k_T\O_{D'}$, so $\phi$ is
 surjective.  As the source and target are free of the same finite
 rank over $\OS$, we deduce that $\phi$ is an isomorphism as claimed.

 We still need to prove that $\phi$ is linear over $\OT$.  By the
 argument of Lemma~\ref{lem-alt-Jk} we reduce to the case where $D$ is
 the divisor with equation $\prod_{i=0}^{d-1}(t-y_i)$ defined over the
 ring 
 \[ R := \O_{\GG^d}=\OS\psb{y_0,\ldots,y_{d-1}}, \]
 and we can invert the discriminant $w=\prod_{i\neq j}(y_i-y_j)$.  We
 reuse the notation in the proof of that lemma, so
 $w^{-1}\OD=F(N,w^{-1}R)$ and $w^{-1}\O_{D^k}=F(N^k,w^{-1}R)$ and
 $w^{-1}\O_{P_kD}=F(N_k,w^{-1}R)$.  We see from
 Proposition~\ref{prop-invariant} that $w^{-1}\OT$ is the image of
 $w^{-1}\O_{D^k}^{\Sg_k}$ in $w^{-1}\O_{P_kD}$, which is the ring
 $F(N_k,w^{-1}R)^{\Sg_k}$ of symmetric functions from $N_k$ to
 $w^{-1}R$.  If we write $N_k^+=\{\um\in N^k\st n_0>\ldots>n_{k-1}\}$
 then $N_k=\Sg_k\tm N_k^+$ as $\Sg_k$-sets so
 $w^{-1}\O_{P_kD}=F(N_k^+,w^{-1}R)$.  On the other hand, $\OT$ is also
 a quotient of $R\hot\O_{\Div_k^+}$, which is the ring of symmetric
 power series in $k$ variables over $R$; a symmetric power series $p$
 corresponds to the function
 $\um\mapsto p(y_\um):=p(y_{n_0},\ldots,y_{n_{k-1}})$.

 If $\um\in N^k$ we put $e_\um=e_{n_0}\ot\ldots\ot e_{n_{k-1}}$, so
 these elements form a basis for $w^{-1}\O_D^{\ot k}$ over $w^{-1}R$.
 Similarly, the set $\{\mu_k(e_\um)\st\um\in N_k^+\}$ is a basis for
 $w^{-1}\lm^k\OD$.  Using the previous paragraph we see that
 $p.\mu_k(e_\um)=p(y_\um)\mu_k(e_\um)$, which tells us the
 $\OT$-module structure on $w^{-1}\lm^k\OD$.

 We next analyse $w^{-1}\O_{D'}$.  This is a quotient of the ring
 \[ w^{-1}\OT\ot\OD=F(N_k^+,w^{-1}\OD)=
     \prod_{\um\in N_k^+}w^{-1}R\{e_i\st i<d\}.
 \]
 It is not hard to check that the relevant ideal is a product of terms
 $I_\um$, where $I_\um$ is spanned by the elements $e_i$ that do not
 lie in the list $e_{n_0},\ldots,e_{n_{k-1}}$.  Thus
 \begin{align*}
  w^{-1}\O_{D'} &= \prod_\um w^{-1}R\{e_{n_j}\st j<k\} \\
  w^{-1}\lm^k_T\O_{D'} &=
    \prod_\um w^{-1}R . e_{n_0}\Smash\ldots\Smash e_{n_{k-1}}.
 \end{align*}
 Let $e'_\um$ be the element of this module whose $\um$'th component
 is $e_{n_0}\Smash\ldots\Smash e_{n_{k-1}}$, and whose other
 components are zero.  Clearly $\{e'_\um\st\um\in N_k^+\}$ is a basis
 for $w^{-1}\lm^k_T\O_{D'}$ over $w^{-1}R$.  As a symmetric power
 series $p$ corresponds to the function $\um\mapsto p(y_\um)$ and
 $e'_\um$ is concentrated in the $\um$'th factor we have
 $p.e'_\um=p(y_\um)e'_\um$.  It is also easy to see that
 $\phi(\mu_k(e_\um))=e'_\um$, and it follows that $\phi$ is
 $\OT$-linear as claimed.
\end{proof}

We next give a formula for $\phi$ in terms of suitable bases of
$\lm^k_S\OD$ and $\lm^k_{\Sub_r}\O_{D'}$.  (This could be used to give
an alternative proof that $\phi$ is an isomorphism.)
\begin{proposition}
 Suppose we have an element
 $x^{\bt_0}\Smash\ldots\Smash x^{\bt_{k-1}}\in\lm^k_S\OD$, where
 $0\leq\bt_0<\ldots<\bt_{d-k-1}$.  Let $\gm_0,\ldots,\gm_{k-1}$ be the
 elements of $\{0,\ldots,d-1\}\sm\{\bt_0,\ldots,\bt_{d-k-1}\}$, listed
 in increasing order.  Then 
 \[ \phi(x^{\bt_0}\Smash\ldots\Smash x^{\bt_{k-1}}) = 
     \pm x^0\Smash\ldots\Smash x^{k-1} \,.\, 
     \det(c_{k+i-\gm_j})_{0\leq i,j<d-k},
 \] 
 where the elements $c_i$ are the usual parameters of the divisor
 $D'$.
\end{proposition}
\begin{proof}
 For any increasing sequence $\al_0<\ldots<\al_{n-1}$ we write
 $x(\al)=x^{\al_0}\Smash\ldots\Smash x^{\al_{n-1}}$.  We also write
 $e'=x(0,1,\ldots,k-1)$ and $e=x(0,1,\ldots,d-1)$, and we put
 $T=\Sub_k(D)$. 
 
 We certainly have $\phi(x(\bt))=b_\bt e'$ for some $b_\bt\in\OT$.  To
 analyse these elements, put $J'=\ker(\OT\ot\OD\xra{}\O_{D'}$, which
 is freely generated over $\OT$ by $\{x^if_{D'}(x)\st i<d-k\}$.
 Consider the element
 \[ a = f_{D'}\Smash xf_{D'}\Smash\ldots x^{d-k-1}f_{D'} \in
         \lm^{d-k}J' \subset \OT\ot\lm^{d-k}\OD.
 \]
 This clearly annihilates $J'\subset\OT\ot\lm^1\OD$, so
 multiplication by $a$ induces a map
 $\lm^k\O_{D'}\xra{}\OT\ot\lm^d\OD$.  As $f_{D'}$ is monic of degree
 $k$, we see that $e'a=e$.  It follows that
 $x(\bt)a=b_\bt e'a=b_\bt e$.  

 On the other hand, we can expand $a$ in the form
 $a=\sum_\gm a_\gm x(\gm)$, where $\gm$ runs over sequences
 $0\leq\gm_0<\ldots<\gm_{d-k-1}<d$.  We have $x(\bt)x(\gm)=\pm e$ if
 $\bt$ and $\gm$ are related as in the statement of the proposition,
 and $x(\bt)x(\gm)=0$ otherwise.  It follows that
 $x(\bt)a=\pm a_\gm e$, and thus that $b_\bt=\pm a_\gm$.  

 Let $A_\gm$ be the matrix whose $(i,j)$'th entry is the coefficient
 of $x^{\gm_j}$ in $x^if_{D'}(x)$; it is then clear that
 $a_\gm=\det(A_\gm)$.  On the other hand, we have
 $x^if_{D'}(x)=\sum_mc_mx^{k+i-m}$, so $(A_\gm)_{ij}=c_{k+i-\gm_j}$,
 and the proposition follows.
\end{proof}

\section{Thom spectra of adjoint bundles}
\label{sec-adjoint}

The following proposition is an immediate consequence of
Theorem~\ref{thm-E-Fk} and its proof (the first statement is just the
case $k=d$ of the second statement).
\begin{proposition}\label{prop-Gku}
 Let $V$ be a $d$-dimensional bundle over a space $X$.  Then there are
 natural isomorphisms 
 \begin{align*}
  \tE^0\Sg^{-d}X^{\uu(V)} &= \lm^d_{E^0X}E^0PV \\
  \tE^0\Sg^{-k}G_kV^\uu   &= \lm^k_{E^0X}E^0PV
                             \text{ for $0\leq k\leq d$. \qed}  
 \end{align*}
\end{proposition}

\begin{remark}
 Note that the proposition gives two different descriptions of the
 module $\tE^0\Sg^{-k}G_kV^\uu$: the first statement with $X$ replaced
 by $G_k(V)$ and $V$ by $T$ gives
 \[ \tE^0\Sg^{-k}G_kV^\uu=\lm^k_{E^0G_kV}E^0PT, \]
 whereas the second
 statement gives 
 \[ \tE^0\Sg^{-k}G_k(V)^\uu=\lm^k_{E^0X}E^0PV. \]
 We leave it to the reader to check that these two descriptions are
 related by the isomorphism
 $\phi\:\lm^k_S\OD\xra{}\lm^k_{\Sub_k(D)}\O_{D'}$ of
 Proposition~\ref{prop-phi-iso}.
\end{remark}

In the present section we examine the isomorphisms of
Proposition~\ref{prop-Gku} more carefully.  We will construct a
diagram as follows, whose effect in cohomology will be identified with
the diagram in Proposition~\ref{prop-cross}.
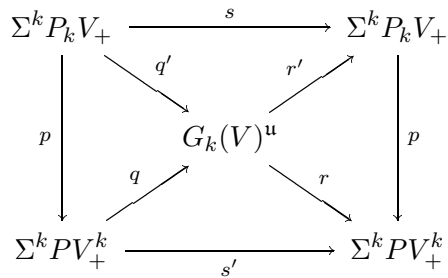
\begin{diag}
 \node{\Sg^k P_kV_+}
  \arrow[2]{e,t}{s}
  \arrow{se,t}{q'}
  \arrow[2]{s,l}{p}
 \node[2]{\Sg^k P_kV_+}
  \arrow[2]{s,r}{p} \\
 \node[2]{G_k(V)^\uu}
  \arrow{ne,t}{r'} 
  \arrow{se,t}{r} \\
 \node{\Sg^k PV^k_+}
  \arrow{ne,t}{q}
  \arrow[2]{e,b}{s'}
 \node[2]{\Sg^k PV^k_+}
\end{diag}%
Here $PV^k$ means the fibre product
\[ PV^k = PV\tm_X\ldots\tm_XPV = 
   \{(x,L_0,\ldots,L_{k-1})\st x\in X\;,\;
        L_0,\ldots,L_{k-1}\in PV_x\}.
\]
Write $Q_kU(V)=F_kU(V)/F_{k-1}U(V)$, so that 
$Q_kU(V)\simeq G_kV^\uu$.  As the filtration of $U(V)$ is
multiplicative, the multiplication $U(V)^k\xra{}U(V)$ induces a map
$(Q_1U(V))^{(k)}\xra{}Q_kU(V)$, or equivalently
$\Sg^kPV_+\xra{}G_k(V)^{\uu}$.  This is the map $q$ in the diagram.

Recall that $P_kV$ is the set of points
$(x,L_0,\ldots,L_{k-1})\in PV^k$ such that the lines $L_i$ are
mutually orthogonal.  The map $p\:P_kV\xra{}PV^k$ is just the
inclusion.  We also have a map $P_kV\xra{}G_kV$ sending $(x,\un{L})$
to $(x,\bigoplus_iL_i)$, and we note that $\uu(\bigoplus_iL_i)$
contains $\bigoplus_i\uu(L_i)$.  Moreover, when $L$ is one-dimensional
there is a canonical isomorphism $\uu(L)\simeq i\R\simeq\R$, so
$\bigoplus_i\uu(L_i)\simeq\R^k$, so we get an inclusion
$\Sg^kP_kV_+\xra{}G_kV^\uu$, which we call $q'$.  It is not hard to
see that this is the same as $qp$, so the left hand triangle commutes
on the nose.

We next define the map $r'\:G_kV^\uu\xra{}\Sg^kP_kV_+$ by a
Pontrjagin-Thom construction.  Let $N'_0\subset\R^k$ be the set of
sequences $(t_0,\ldots,t_{k-1})$ such that $t_0<\ldots<t_{k-1}$, and
let $N'$ be the space of triples $(x,W,\al)$ where $W\in G_kV_x$ and
$\al\in\uu(W)$ and $\al$ has $k$ distinct eigenvalues.  This is easily
seen to be an open subspace of the total space of the bundle $\uu$
over $G_kV$.  Given such a triple, we note that the eigenvalues of
$\al$ are purely imaginary, so we can write them as
$it_0,\ldots,it_{k-1}$ with $t_0<\ldots<t_{k-1}$.  We also put
$L_j=\ker(\al-it_j)$, so the spaces $L_j$ are one-dimensional and
mutually orthogonal, and their direct sum is $W$.  Using this we see
that the map $qp\:\Sg^kP_kV_+\xra{}G_kV^\uu$ induces a homeomorphism
$N'_0\tm P_kV\xra{}N'$, and this gives a collapse map
\[ G_kV^\uu \xra{} N'\cup\{\infty\} \simeq 
    (N'_0\tm P_kV)\cup\{\infty\} \simeq
    (N'_0\cup\{\infty\})\Smash P_kV_+.
\]
On the other hand, the inclusion $N'_0\xra{}\R^k$ gives a collapse map
$S^k\xra{}N'_0\cup\{\infty\}$ which is a homotopy equivalence; after
composing with the inverse of this, we obtain a map
$G_kV^\uu\xra{}\Sg^kP_kV_+$, which we denote by $r'$.

We now define a map $r\:G_kV^\uu\xra{}PV^k_+$.  We first mimic
Lemma~\ref{lem-diagonal} and define maps $m_j\:U(1)\xra{}U(1)$ (for
$0\leq j< k$) by 
\[ m_j(e^{i\tht}) = \begin{cases}
    e^{ik\tht} & \text{ if } j/k\leq\tht/2\pi\leq (j+1)/k \\
    1          & \text{ otherwise }
   \end{cases}
\]
(where $\tht$ is assumed to be in the interval $[0,2\pi]$).  We then
define $\dl'_k\:U(V)\xra{}U(V)^k$ by
$\dl'_k(g)=(m_0(g),\ldots,m_{k-1}(g))$.  This is homotopic to the
diagonal and preserves filtrations so it induces a map
$G_kV^\uu=Q_kU(V)\xra{}Q_k(U(V)^k)$.  The target of this map is the
wedge of all the spaces
$Q_{l_0}U(V)\Smash\ldots\Smash Q_{l_{k-1}}U(V)$ for which
$\sum_il_i=k$.  We can thus project down to the factor
$Q_1U(V)\Smash\ldots\Smash Q_1U(V)=\Sg^kPV^k_+$ to get a map
$G_kV^\uu\xra{}\Sg^kPV^k_+$, which we call $r$.

It is not hard to recover the following more explicit description of
$r$.  Recall that we have a homeomorphism 
\[ \gm\: U(1)\xra{} \R\cup\{\infty\} \]
given by $\gm(z)=(z+1)(z-1)^{-1}/i$ and $\gm^{-1}(t)=(it+1)/(it-1)$.
One checks that $\gm(e^{i\tht})=-\cot(\tht/2)$, which is a strictly
increasing function of $\tht$ for $0<\tht<2\pi$.  
Let $A_j$ denote the arc $\{e^{i\tht}\st j/k<\tht/2\pi<(j+1)/k\}$, so
$\gm A_j$ is the interval $(-\cot(\pi j/k),-\cot(\pi(j+1)/k))$.  We
also define $\ov{m}_j=\gm m_j\gm^{-1}$, which can be regarded as a
homeomorphism $\gm A_j\cup\{\infty\}\xra{}\R\cup\{\infty\}$, homotopy
inverse to the evident collapse map in the opposite direction.  If we
put $N_0=\prod_j\gm A_j\subset N'_0\subset\R^k$ then the maps
$\ov{m}_j$ combine to give a homeomorphism
$\ov{m}\:N_0\cup\{\infty\}\xra{}\R^k\cup\{\infty\}$, which is again
homotopy inverse to the evident collapse map in the opposite
direction.  Now let $N\subset N'$ be the space of triples $(x,W,\al)$
such that $\al/i$ has precisely one eigenvalue in $\gm A_j$ for each
$j$.  If $(x,W,\al)\in N$ and $t_j$ is the eigenvalue in $\gm A_j$ and
$L_j=\ker(\al-it_j)$ then we find that $\un{t}\in N_0$ and $\un{L}\in
P_kV_x$ and $r(x,W,\al)=(\ov{m}(\un{t}),x,\un{L})$.  On the other
hand, if $(x,W,\al)\not\in N$ we find that $r(x,W,\al)=\infty$.  It
follows that $r$ is constructed in the same way as $r'$, except that
$N'$ and $N'_0$ are replaced by the smaller sets $N$ and $N_0$.  The
projections $N'\cup\{\infty\}\xra{}N\cup\{\infty\}$ and
$N'_0\cup\{\infty\}\xra{}N_0\cup\{\infty\}$ are homotopy equivalences,
and it follows that $r$ is homotopic to $pr'$.  This shows that the
right hand triangle in our diagram commutes up to homotopy.

We now consider the composite $s=r'q'\:\Sg^kP_kV_+\xra{}\Sg^kP_kV_+$,
which is essentially obtained by collapsing out the complement of
$(q')^{-1}(N')$.  There is an evident action of the symmetric group
$\Sg_k$ on the space $\Sg^kP_kV_+$, given by 
\[ \sg.(\un{t},\un{L}) = 
    (t_{\sg^{-1}(0)},\ldots,t_{\sg^{-1}(k-1)},
     L_{\sg^{-1}(0)},\ldots,L_{\sg^{-1}(k-1)}).
\]
One checks that $(q')^{-1}(N')=\coprod_\sg\sg.(N'_0\tm P_kV)$, and
using this one can see that $s$ is just the trace map
$\tr_{\Sg_k}=\sum_{\sg\in\Sg_k}\sg$.

Finally, we define 
\[ s'=rq\:\Sg^kPV^k_+\xra{}\Sg^kPV^k_+. \]
We can also define $\tr_{\Sg_k}\:\Sg^kPV^k_+\xra{}\Sg^kPV^k_+$; we
suspect that this is \emph{not} the same as $s'$, although we will see
shortly that it induces the same map in cohomology.

We now apply the functor $\tE^k(-)=\tE^0\Sg^{-k}(-)$ to our diagram of
spaces and write $D=D(V)$ to get the following diagram:
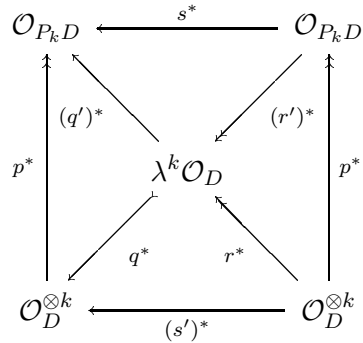
\begin{diag}
 \node{\O_{P_kD}}
 \node[2]{\O_{P_kD}}
  \arrow[2]{w,t}{s^*} 
  \arrow{sw,b,A}{(r')^*} \\  
 \node[2]{\lm^k\OD} 
  \arrow{nw,b}{(q')^*} 
  \arrow{sw,b,V}{q^*} \\
 \node{\O_D^{\ot k}} 
  \arrow[2]{n,l,A}{p^*}
 \node[2]{\O_D^{\ot k}} 
  \arrow[2]{n,r,A}{p^*}
  \arrow{nw,b,A}{r^*}
  \arrow[2]{w,b}{(s')^*}
\end{diag}%
The map $p^*\:\O_D^{\ot k}=\O_{D^k}\xra{}\O_{P_kD}$ is the same as
considered previously; this is the definition of our identification of
$(P_kV)_E$ with $P_kD$.  It follows from the Hopf algebra isomorphism
of Theorem~\ref{thm-EUV} that $r^*=\mu_k$ and $q^*=\psi_k$, and thus
that $(s')^*=\psi_k\mu_k=\alt_k$.  As $\mu_k$ factors uniquely through
$p^*$ we must have $(r')^*=\mu'_k$.  As $q'=pq$ and
$\psi'_k=p^*\psi_k$ we have $(q')^*=\psi'_k$.  Finally, we know that
$s=\tr_{\Sg_k}$ and any permutation $\sg\in\Sg_k$ acts on the sphere
$S^k$ with degree equal to its signature so it follows that
$s^*=\alt_k\:E^0P_kV\xra{}E^0P_kV$.

\section{Fibrewise loop groups}
\label{sec-loops}

We conclude the main part of this paper by studying the fibrewise loop
space $\Om_XU(V)$ and thereby providing a topological realisation of
the diagram in Proposition~\ref{prop-invariant}.

First, the group structure on $U(V)$ gives a group structure on
$\Om_XU(V)$.  We also have $\Om_XU(V)\simeq\Om^2_XBU(V)$, and there is
a canonical homotopy showing that a double loop space is
homotopy-commutative, so the proof goes through to show that
$\Om_XU(V)$ is fibrewise homotopy commutative.

We next recall certain subspaces of $\Om_XU(V)$ which have been
considered by a number of previous authors --- we will mostly refer to
Crabb's exposition~\cite{cr:ssu}, which cites on Mitchell's
paper~\cite{mi:fls} and (apparently unpublished) work of Mahowald and
Richter. 

Let $V$ be a vector space.  Any finite Laurent series
$f(z)=\sum_{i=-N}^Na_iz^i$ with coefficients $a_i\in\End(V)$ can be
regarded as a map $U(1)\xra{}\End(V)$ which we can compose with the
standard homeomorphism $\gm^{-1}\:S^1\xra{}U(1)$ to get a map
$\hat{f}\:S^1\xra{}\End(V)$.  We write $\Oml U(V)$ for the space of
based loops $u\:S^1\xra{}U(V)$ that have the form $u=\hat{f}$ for some
finite Laurent series $f$, and call this the space of Laurent loops.
Similarly, we write $\Omp U(V)$ for the space of loops that have the
form $\hat{f}$ for some polynomial $f$.

If $u=\hat{f}$ is a Laurent loop we have
$f(z)^{-1}=f(z)^*=\sum_ia_i^*z^{-i}$ which is again a finite Laurent
series.  Using this we see that that $\Oml U(V)$ is a subgroup of 
$\Om U(V)$ (but $\Omp U(V)$ is merely a submonoid).  We also find that
the function $d(z)=\det(f(z))$ is a finite Laurent series in
$\C[z,z^{-1}]$ satisfying $d(z)\ov{d}(\ov{z})=1$ and 
$d(1)=1$; it follows easily that $d(z)=z^n$ for some integer $n$,
called the degree of $u$.

\begin{definition}
 \begin{itemize}
  \item[\rm(a)] We write $S_kV$ for the space of polynomial loops of
   degree $k$ on $U(V)$.
  \item[\rm(b)] The product structure on $\Om U(V)$ induces maps
   $S_kV\tm S_lV\xra{}S_{k+l}V$, which we call $\mu_{kl}$.
  \item[\rm(c)] Given $W\in G_kV$ and $z\in U(1)$ we have a polynomial
   $z\pi_W+(1-\pi_W)\in\End(V)[z]$ giving rise to a based loop in
   $U(V)$ which we call $\phi_k(W)$.  This defines a map
   $\phi_k\:G_kV\xra{}S_kV$.  It is not hard to show that
   $\phi_1\:PV\xra{}S_1V$ is a homeomorphism.
  \item[\rm(d)] By combining $\phi_1$ with the product map we get a map
   $\nu_k\:PV^k_X\xra{}S_kV$. 
 \end{itemize}
 If $V$ is a bundle rather than a vector space, we make all these
 definitions fibrewise in the obvious way.
\end{definition}

Note that $\nu_k$ induces a map $E^*S_kV\xra{}(E^*PV)^{\ot k}$, where
the tensor product is taken over $E^*X$.  We write $\Sym^k(E^*PV)$ for
the submodule invariant under the action of $\Sg_k$.
\begin{proposition}
 $\nu_k$ induces an isomorphism $E^*S_kV=\Sym^k(E^*PV)$, and
 thus an isomorphism $D(V)^k/\Sg_k\xra{}(S_kV)_E$.
\end{proposition}
\begin{proof}
 Put $d=\dim(V)$ and let $A$ be the set of lists $\al=(\al_i\st i<k)$
 with $0\leq\al_i<d$.  We have 
 \[ E^*PV^k_X=(E^*PV)^{\ot k}=E^*\psb{x_i\st i<k}/(f_V(x_i)\st i<k),
 \] 
 and the set $\{x^\al\st\al\in A\}$ is a basis for this ring over
 $E^*X$.  Put $A_+=\{\al\in A\st \al_0\leq\ldots\leq\al_{k-1}\}$ and
 $M_+=\bigoplus_{\al\in A_+}E^*X$, and let $\pi\:E^*PV^k_X\xra{}M_+$
 be the obvious projection.  This clearly induces an isomorphism
 $\Sym^kE^*PV\xra{}M_+$.  

 We take as our basic input (proved by Mitchell~\cite{mi:fls}) 
 the fact that when $X$ is a point, the map $\nu_k$ induces an
 isomorphism $(H_*PV)^{\ot k}_{\Sg_k}\xra{}H_*S_kV$.  In particular,
 this means that $H_*S_kV$ is a finitely generated free Abelian group,
 concentrated in even degrees.  By duality we see that
 $H^*S_kV=\Sym^kH^*PV$, and thus that the map
 $\pi\nu_k^*\:H^*S_kV\xra{}\bigoplus_{A_+}H^*X$ is an isomorphism.
 Using an Atiyah-Hirzebruch spectral sequence we see that
 $\pi\nu_k^*\:E^*S_kV\xra{}\bigoplus_{A_+}E^*X$ is an isomorphism for
 any $E$.

 Now let $X$ be arbitrary.  If $V$ is trivialisable with fibre $V_0$
 then $S_kV=X\tm S_kV_0$ and it follows from the above that
 $\pi\nu_k^*$ is an isomorphism.  If $V$ is not trivialisable, we can
 still give $X$ a cell structure such that the restriction to any
 closed cell is trivialisable, and then use Mayer-Vietoris sequences,
 the five lemma, and the Milnor sequence to see that $\pi\nu_k^*$ is
 an isomorphism.

 We next claim that the maps 
 \[ \mu_{kl}^*\:E^*S_{k+l}V\xra{}E^*S_kV\ot_{E^*X}E^*S_lV \]
 give rise to a cocommutative coproduct.  To see this, let $C(V)$
 denote the following diagram:
 \begin{diag}
  \node{S_kV\tm S_lV}
   \arrow{s,l}{\text{twist}}
   \arrow{e,t}{\mu_{kl}}
  \node{S_{k+l}V}
   \arrow{s,r}{1} \\
  \node{S_lV\tm S_kV}
   \arrow{e,b}{\mu_{lk}}
  \node{S_{k+l}V.}
 \end{diag}%
The claim is that the diagram $E^*C(V)$ commutes.  Let
 $i_0,i_1\:V\xra{}V^2$ be the two inclusions.  The map $i_0$
 induces a map $E^*C(V^2)\xra{}E^*C(V)$, and it follows easily from
 our previous discussion that this is surjective.  It will thus be
 enough to show that the two ways round $E^*C(V^2)$ become the same
 when composed with the map
 \[ (S_k(i_0)\tm S_l(i_0))^*\:
     E^*(S_k(V^2)\tm_X S_l(V^2))\xra{}E^*(S_kV\tm_XS_lV).
 \]
 It is standard that $i_0$ is homotopic to $i_1$ through linear
 isometries, so $S_l(i_0)$ is fibre-homotopic to $S_l(i_1)$.
 Similarly, the identity map of $S_{k+l}(V^2)$ is homotopic to
 $S_{k+l}(\text{twist})$.  It is thus enough to check that the two
 composites $S_kV\tm S_lV\xra{}S_{k+l}(V^2)$ in the following diagram
 are the same:
 \begin{diag}
  \node{S_k(V)\tm S_l(V)}
   \arrow{e,t}{\scriptscriptstyle{S_k(i_0)\tm S_k(i_1)}}
  \node{S_k(V^2)\tm S_l(V^2)}
   \arrow{s,l}{\text{twist}}
   \arrow{e,t}{\mu_{kl}}
  \node{S_{k+l}(V^2)}
   \arrow{s,r}{\scriptscriptstyle{S_{k+l}(\text{twist})}} \\
  \node[2]{S_l(V^2)\tm S_k(V^2)}
   \arrow{e,b}{\mu_{lk}}
  \node{S_{k+l}(V^2).}
 \end{diag}%
This is easy to see directly.

 We now see that the map 
 \[ \nu_k^*\:E^*S_kV\xra{}(E^*PV)^{\ot k} \]
 factors through $\Sym^k(E^*PV)$.  As the map
 $\pi\:\Sym^k(E^*PV)\xra{}M_+$ and its composite with $\nu_k^*$ are
 both isomorphisms, we deduce that
 $\nu_k^*\:E^*S_kV\xra{}\Sym^k(E^*PV)$ is an isomorphism as claimed.
\end{proof}

\begin{corollary}
 The formal scheme $(\Om_XU(V))_E$ is the free commutative formal
 group over $X_E$ generated by the divisor $DV$.
\end{corollary}
\begin{proof}
 We refer to~\cite[Section 6.2]{st:fsfg} for background on free
 commutative formal groups; the results there mostly state that the
 obvious methods for constructing such objects work as expected under
 some mild hypotheses.  Given a formal scheme $T$ over a formal scheme
 $S$, we use the following notation:
 \begin{itemize}
  \item[\rm(a)] $M^+T$ is the free commutative monoid over $S$ generated
   by $T$.  This is characterised by the fact that monoid homomorphisms
   from $M^+T$ to any monoid $H$ over $S$ biject with maps $T\xra{}H$
   of schemes over $S$.  It is clear that if there exists an $M^+T$
   with this property, then it is unique up to canonical isomorphism.
   Similar remarks apply to our other definitions.  In reasonable cases
   we can construct the colimit $\coprod_kT^k_S/\Sg_k$ and this works
   as $M^+T$; see~\cite[Proposition 6.8]{st:fsfg} for technicalities.
  \item[\rm(b)] $MT$ is the free commutative group over $S$ generated by
   $T$.
  \item[\rm(c)] If $T$ has a specified section $z\:S\xra{}T$, then $N^+T$
   is the free commutative monoid scheme generated by the based scheme
   $T$, so homomorphisms from $N^+T$ to $H$ biject with maps
   $T\xra{}H$ such that the composite $S\xra{z}T\xra{}H$ is zero.  In
   reasonable cases $N^+T$ can be constructed as
   $\colim_kT^k_S/\Sg_k$. 
  \item[\rm(d)] If $T$ has a specified section we also write $NT$ for the
   free commutative group over $S$ generated by the based scheme $T$.
 \end{itemize}
 The one surprise in the theory is that often $NT=N^+T$; this is
 analogous to the fact that a graded connected Hopf algebra
 automatically has an antipode.  It is easy to check that
 $MT=\Z\tm NT$, where $\Z$ is regarded as a discrete group scheme in
 an obvious way.

 We first suppose that $V$ has a one-dimensional summand, so
 $V=L\op W$ for some bundles $L$ and $W$ with $\dim(L)=1$.  Note that
 for each $x\in X$ there is a canonical isomorphism $\C\xra{}\End(L)$
 giving $U(1)\tm X\simeq U(L)$.  This gives
 an evident inclusion $\al\:U(1)\tm X\xra{}U(V)$ with
 $\det\circ\al=1$.  We define $\bt\:U(V)\xra{}SU(V)$ by
 $\bt(g)=\al(\det(g))^{-1}g$, and note that $\bt(\al(z)g)=\bt(g)$ for
 all $z$.

 We also have an evident map $z\:X\simeq PL\xra{}PV\simeq S_1V$
 splitting the projection $PV\xra{}X$.  Left multiplication by $z$
 gives a map $i_k\:S_kV\xra{}S_{k+1}V$.  We also define
 $i_k\:S_kV\xra{}\Om_XSU(V)$ to be the restriction of $\Om_X\bt$ to
 $S_kV\subset\Om_XU(V)$.  Using the fact that $\bt(\al(z)g)=\bt(g)$ we
 see that $j_{k+1}i_k=j_k$.  Thus, if we define $S_\infty V$ to be the
 homotopy colimit of the spaces $S_kV$, we get a map
 $j_\infty\:S_\infty V\xra{}\Om_XSU(V)$ of spaces over $X$.  Using the
 usual bases for $E^*S_kV=\Sym^k(E^*PV)$ we find that the maps
 $i_k^*\:\Sym^{k+1}(E^*PV)\xra{}\Sym^k(E^*PV)$ are surjective.  It
 follows using the Milnor sequence that
 $E^*S_\infty V=\invlim_k\Sym^k(E^*PV)$ and thus that
 $(S_\infty V)_E=\colim_kDV^k/\Sg_k$.  We claim that this is the same
 as $N^+DV$; this is clear modulo some categorical technicalities,
 which are covered in~\cite[Section 6.2]{st:fsfg}.  In the case where
 $X$ is a point, it is well-known and easy to check (by calculation in
 ordinary homology) that the map $S_\infty V\xra{}\Om_XSU(V)$ is a weak
 equivalence.  In the general case we have a map between fibre bundles
 that is a weak equivalence on each fibre; it follows easily that the
 map is itself a weak equivalence, and thus that $\Om_XSU(V)_E=N^+DV$.
 On the other hand, as $\Om_XSU(V)$ is actually a group bundle, we see
 that $\Om_XSU(V)_E$ is a formal group scheme, so $N^+DV=NDV$.

 We now turn to the groups $\Om_XU(V)$.  We define $\Z_X=\Z\tm X$,
 viewed as a bundle of groups over $X$ in the obvious way.  This can
 be identified with $\Om_X(U(1)\tm X)$ so the determinant map gives
 rise to a homomorphism $\dl\:\Om_XU(V)\xra{}\Z_X$.  Given
 $(n,x)\in\Z_X$ we have a homomorphism $U(1)\xra{}U(V_x)$ given by
 $z\mapsto\al(z^n)$.  This construction gives us a map
 $\sg\:\Z_X\xra{}\Om_XU(V)$ with $\dl\sg=1$ and thus a splitting
 $\Om_XU(V)\simeq\Z\tm\Om_XSU(V)$ and thus an isomorphism
 $\Om_XU(V)_E=\Z\tm NDV=MDV$.  One can check that the various uses of
 the map $\al$ cancel out and that the standard inclusion
 $DV\xra{}MDV$ is implicitly identified with the map coming from the
 inclusion $PV=S_1V\xra{}\Om_XU(V)$.  This proves the corollary in the
 case where $V$ has a one-dimensional summand.

 Now suppose that $V$ does not have such a summand.  We have an
 evident coequaliser diagram $PV\tm_XPV\arrow{e,=>}PV\xra{}X$, giving
 rise to a coequaliser diagram $DV\tm_{X_E}DV\xra{}DV\xra{}X_E$ of
 schemes over $X_E$, in which the map $DV\xra{}X_E$ is faithfully
 flat.  The pullback of $V$ to $PV$ has a tautological one-dimensional
 summand, which implies that $(PV\tm_X\Om_XU(V))_E$ has the required
 universal property in the category of formal group schemes over
 $PV_E$.  Similar remarks apply to $PV\tm_XPV\tm_X\Om_XU(V)$.  It
 follows by a descent argument that $\Om_XU(V)$ itself has the
 required universal property, as one sees easily
 from~\cite[Proposition~2.76 and Remark~4.52]{st:fsfg}.
\end{proof}

We next recall the standard line bundle over $S_kV$, which we will
call $T$; see~\cite{cr:ssu,mi:fls} for more details.  Write $A=\C[z]$
and $K=\C[z,z^{-1}]$.  A point of $S_kV$ has the form $(x,\hat{f})$
for some $f\in\End_\C(V_x)[z]\simeq\End_A(A\ot V_x)$.  Multiplication
by $f$ defines a surjective endomorphism $m(f)$ of $(K/A)\ot V$, and
we define $T_{(x,f)}$ to be the kernel of this endomorphism.  One can
check that this always has dimension $k$ over $\C$ and that we get a
vector bundle.  This is classified by a map
$\tau_k\:S_kV\xra{}BU(k)\tm X$ of spaces over $X$.  It is easy to see
that the restriction of $T$ to $G_kV\subset S_kV$ is just the
tautological bundle.

There are evident short exact sequences 
\[ \ker(m(g)) \arrow{e,V} \ker(m(fg)) 
              \arrow{e,t,A}{m(g)} \ker(m(f)), 
\]
which can be split using the inner products to give isomorphisms
$\mu_{kl}^*T\simeq\pi_0^*T\op\pi_1^*T$ over $S_kV\tm_XS_lV$.  This
means that the map $\tau\:\coprod_kS_kV\xra{}(\coprod_kBU(k))\tm X$ is
a homomorphism of $H$-spaces over $X$.

We now have a diagram of spaces as follows:
\begin{diag}
 \node{P_kV} 
  \arrow{s,l}{q'}
  \arrow{e,t}{p} 
 \node{PV^k_X}
  \arrow{s,l}{\nu_k}
  \arrow{e,t}{\tau_1^k}
 \node{(\CPi)^k\tm X}
  \arrow{s} \\
 \node{G_kV}
  \arrow{e,b}{\phi_k}
 \node{S_kV}
  \arrow{e,b}{\tau_k}
 \node{BU(k)\tm X.}
\end{diag}%
It is easy to identify the corresponding diagram of schemes with the
diagram of Proposition~\ref{prop-invariant}.

\appendix

\section{Appendix : Functional calculus}
\label{apx-functional}

In this appendix we briefly recall some basic facts about functional
calculus for normal operators.  An endomorphism $\al$ of a vector
space $V$ is \emph{normal} if it commutes with its adjoint.  For us
the relevant examples are Hermitian operators (with $\al^*=\al$),
anti-Hermitian operators (with $\al^*=-\al$) and unitary operators
(with $\al^*=\al^{-1}$).  

For any operator $\al$ and any $\lm\in\C$ we have
\[ \ker(\al-\lm)^\perp=\img((\al-\lm)^*)=\img(\al^*-\ov{\lm}). \]
If $\al$ is normal we deduce that $\ker(\al-\lm)^\perp$ is preserved
by $\al$, and it follows easily that $V$ is the orthogonal direct sum
of the eigenspaces of $\al$.  It follows in turn that the operator
norm of $\al$ (defined by $\|\al\|=\sup\{\|\al(v)\|\;:\;\|v\|=1\}$) is
just the same as the spectral radius (defined as the maximum absolute
value of the eigenvalues of $\al$).

Now let $X$ be a subset of $\C$ containing the eigenvalues of $\al$,
and let $f\:X\xra{}\C$ be a continuous function.  We define $f(\al)$
to be the endomorphism of $V$ that has eigenvalue $f(\lm)$ on the
space $\ker(\al-\lm)$.  From this definition it is clear that the
following equations are valid whenever they make sense:
\begin{align*}
 c(\al) &= c.1_V \text{ if $c$ is constant } \\
 \text{id}(\al) &= \al \\
 \text{Re}(\al) &= (\al+\al^*)/2 \\
 \text{Im}(\al) &= (\al-\al^*)/(2i) \\ 
 (f+g)(\al) &= f(\al) + g(\al) \\
 (fg)(\al)  &= f(\al)g(\al) \\
 \ov{f}(\al) &= f(\al)^* \\
 (f\circ g)(\al) &= f(g(\al)) \\
 \|f(\al)\| &\leq \sup_{x\in X}|f(x)|.
\end{align*}

The continuity properties of $f(\al)$ are less clear from our
definition.  However, they are provided by the following result.  
\begin{proposition}
 Let $X$ be a closed subset of $\C$, and $V$ a vector space.  Let
 $N(X,V)$ be the set of normal operators on $V$ whose eigenvalues lie
 in $X$, and let $C(X,\C)$ be the set of continuous functions from $X$
 to $\C$ (with the topology of uniform convergence on compact sets).
 Define a function $E\:C(X,\C)\tm N(X,V)\xra{}\End(V)$ by
 $E(f,\al)=f(\al)$.  Then $E$ is continuous.
\end{proposition}
\begin{proof}
 Let $A$ be the set of functions $f\in C(X,\C)$ for which the function
 $\al\mapsto f(\al)$ is continuous.  Using the above algebraic
 properties, we see that $A$ is a subalgebra of $C(X,\C)$ containing
 the functions $z\mapsto\text{Re}(z)$ and $z\mapsto\text{Im}(z)$.  By
 the Stone-Weierstrass theorem, it is dense in $C(X,\C)$.  Now suppose
 we have $f\in C(X,\C)$, $\al\in N(X,V)$ and $\ep>0$.  Put
 $Y=\{x\in X\st |x|\leq \|\al\|+1\}$, which is compact.  As $A$ is
 dense we can choose $p\in A$ with $|f-p|<\ep/4$ on $Y$.  As $p\in A$
 can choose $\dl$ such that $\|p(\bt)-p(\al)\|<\ep/4$ whenever
 $\|\bt-\al\|<\dl$.  We may also assume that $\dl<1$, which means that
 when $\|\bt-\al\|<\dl$ we have $\bt\in Y$.  Now if $|f-g|<\ep/4$ on
 $Y$ and $\|\al-\bt\|<\dl$ then
 \begin{align*}
  \|f(\al)-g(\bt)\| \leq & 
    \|f(\al)-p(\al)\| + \|p(\al)-p(\bt)\| + \\
  & \|p(\bt)-f(\bt)\| + \|f(\bt)-g(\bt)\| \\
  < & \ep/4 + \ep/4 + \ep/4 + \ep/4 = \ep,
 \end{align*}
 as required.
\end{proof}

The following proposition is an elementary exercise in linear algebra.
\begin{proposition}\label{prop-square}
 Let $\al\:V\xra{}W$ be a linear map.  Then $\al^*\al$ and $\al\al^*$
 are self-adjoint endomorphisms of $V$ and $W$ with nonnegative
 eigenvalues.  For each $t>0$ the map $\al$ gives an isomorphism of
 $\ker(\al^*\al-t)$ with $\ker(\al\al^*-t)$, so the nonzero
 eigenvalues of $\al^*\al$ and their multiplicities are the same as
 those of $\al\al^*$.  If $f\:[0,\infty)\xra{}\R$ then
 $\al\circ f(\al^*\al)=f(\al\al^*)\circ\al$. \qed
\end{proposition}

\begin{definition}\label{defn-ek}
 We write $w(V)=\{\al\in\End(V)\st\al^*=\al\}$ (the space of
 self-adjoint endomorphisms of $V$).  If $\al\in w(V)$ then the
 eigenvalues of $\al$ are real, so we can list them in descending
 order, repeated according to multiplicity.  We write $e_k(\al)$ for
 the $k$'th element in this list, so $e_1(\al)\geq\ldots\geq e_n(\al)$
 and $\det(t-\al)=\prod_k(t-e_k(\al))$.
\end{definition}

We will need the following standard result:
\begin{proposition}\label{prop-ek-cont}
 The functions $e_k\:w(V)\xra{}\R$ are continuous.
\end{proposition}
\begin{proof}
 Let $\gm$ be a simple closed curve in $\C$ and let $m$ be an integer.
 Let $U$ be the set of endomorphisms of $V$ that have precisely $m$
 eigenvalues (counted according to multiplicity) inside $\gm$, and no
 eigenvalues on $\gm$.  A standard argument with Rouch\'e's theorem
 shows that $U$ is open in $\End(V)$.

 Given real numbers $r\leq R$, consider the rectangular contour
 $\gm_{r,R}$ with corners at $r\pm i$ and $R\pm i$.  Clearly
 $e_k(\al)>r$ iff $\al$ has at least $k$ eigenvalues inside
 $\gm_{r,R}$ for some $R$.  It follows that $\{\al\st e_k(\al)>r\}$ is
 open, as is $\{\al\st e_k(\al)<r\}$ by a similar argument.  This
 implies that $e_k$ is continuous.
\end{proof}

\Addresses\recd

\end{document}